\pdfoutput=1
\documentclass[11pt]{article}
\usepackage[margin=1in]{geometry}
\usepackage{mathtools}
\usepackage{enumitem}
\usepackage{amsmath}
\usepackage{amssymb}
\usepackage{amsthm}
\usepackage{subcaption}
\usepackage{etoolbox}
\usepackage[hang,flushmargin]{footmisc}
\usepackage{microtype}
\usepackage{titlesec}
\usepackage[longnamesfirst]{natbib}
\usepackage{xcolor}
\usepackage[hypertexnames=false]{hyperref}
\usepackage[norefs,nocites]{refcheck}

\hypersetup{colorlinks=true,linkcolor=blue,urlcolor=blue,citecolor=blue}
\DeclareMathOperator*{\argmin}{arg\,min}
\DeclareMathOperator{\rank}{rank}
\DeclareMathOperator{\Tr}{Tr}
\DeclareMathOperator{\Var}{Var}
\renewcommand{\P}{\ensuremath{\mathbb{P}}}
\newcommand{\E}{\ensuremath{\mathbb{E}}}
\newcommand{\Q}{\ensuremath{\mathbb{Q}}}
\newcommand{\cF}{\ensuremath{\mathcal{F}}}
\newcommand{\rF}{\ensuremath{\mathrm{F}}}
\newcommand{\Bias}{\ensuremath{\mathrm{Bias}\hspace*{0.2mm}}}
\newcommand{\Unif}{\ensuremath{\mathrm{Unif}\hspace*{0.2mm}}}
\newcommand{\Cov}{\ensuremath{\mathrm{Cov}\hspace*{0.2mm}}}
\newcommand{\cA}{\ensuremath{\mathcal{A}}}
\newcommand{\cB}{\ensuremath{\mathcal{B}}}
\newcommand{\cR}{\ensuremath{\mathcal{R}}}

\newcommand{\cC}{\ensuremath{\mathcal{C}}}
\newcommand{\bX}{\ensuremath{\mathbf{X}}}
\newcommand{\bY}{\ensuremath{\mathbf{Y}}}
\newcommand{\bW}{\ensuremath{\mathbf{W}}}
\newcommand{\cH}{\ensuremath{\mathcal{H}}}
\newcommand{\N}{\ensuremath{\mathbb{N}}}
\newcommand{\cN}{\ensuremath{\mathcal{N}}}
\newcommand{\I}{\ensuremath{\mathbb{I}}}
\newcommand{\cX}{\ensuremath{\mathcal{X}}}
\newcommand{\cW}{\ensuremath{\mathcal{W}}}
\newcommand{\R}{\ensuremath{\mathbb{R}}}
\newcommand{\T}{\ensuremath{\mathsf{T}}}
\newcommand{\Z}{\ensuremath{\mathbb{Z}}}

\newlist{inlineroman}{enumerate*}{1}
\setlist[inlineroman]{afterlabel=~,label=(\roman*)}
\newcommand{\vvvert}{{\vert\kern-0.25ex\vert\kern-0.25ex\vert}}
\newcommand{\bigvvvert}{{\big\vert\kern-0.35ex\big\vert\kern-0.35ex\big\vert}}
\newcommand{\Bigvvvert}{{\Big\vert\kern-0.3ex\Big\vert\kern-0.3ex\Big\vert}}
\newcommand{\biggvvvert}{{\bigg\vert\kern-0.3ex\bigg\vert\kern-0.3ex\bigg\vert}}
\newcommand{\diff}[1]{\,\mathrm{d}#1}

\titleformat{\paragraph}[hang]{\bfseries\upshape}{}{0pt}{}[]
\titlespacing*{\paragraph}{0pt}{6pt}{0pt}
\newcounter{proofparagraphcounter}
\newcommand{\proofparagraph}[1]{
  \refstepcounter{proofparagraphcounter}%
\paragraph{Part \theproofparagraphcounter : #1}}%

\newtheoremstyle{break}{\topsep}{\topsep}{\itshape}{}{\bfseries}{}{\newline}%
{\thmname{#1}\thmnumber{ #2}%
\thmnote{ \rm(#3)\addcontentsline{toc}{subsubsection}{{\it#1 #2}\/: #3}}}

\newtheoremstyle{breakplainproof}{\topsep}{\topsep}{}{}{\bfseries}{}{\newline}%
{\thmname{#1}\thmnumber{ #2}%
\thmnote{ \rm(#3)\addcontentsline{toc}{subsubsection}{{\it#1}\/: #3}}}

\theoremstyle{break}
\newtheorem{theorem}{Theorem}[section]
\newtheorem{lemma}{Lemma}[section]

\newtheorem{corollary}{Corollary}[section]
\newtheorem{proposition}{Proposition}[section]
\newtheorem{remark}{Remark}[section]

\theoremstyle{breakplainproof}

\newtheorem*{myproof}{Proof}
\AtBeginEnvironment{myproof}{\setcounter{proofparagraphcounter}{0}}%
\AtEndEnvironment{myproof}{\null\hfill$\square$}%

\newcommand{\myref}[1]{\iftoggle{aos}{\ref*{#1}}{\ref{#1}}}

\newtoggle{aos}
\togglefalse{aos}

\makeatletter%
\begin{document}

\title{Yurinskii's Coupling for Martingales}

\author{
  Matias D.\ Cattaneo\textsuperscript{1}
  \and Ricardo P.\ Masini\textsuperscript{2}
  \and William G.\ Underwood\textsuperscript{3*}
}

\maketitle

\footnotetext[1]{
  Department of Operations Research
  and Financial Engineering,
  Princeton University
}
\footnotetext[2]{
  Department of Statistics,
  University of California, Davis
}
\footnotetext[3]{
  Statistical Laboratory,
  University of Cambridge
}
\let\thefootnote\relax
\footnotetext[1]{
  \textsuperscript{*}Corresponding author:
  \href{mailto:wgu21@cam.ac.uk}{\texttt{wgu21@cam.ac.uk}}
}
\addtocounter{footnote}{-1}\let\thefootnote\svthefootnote

\setcounter{page}{0}\thispagestyle{empty}

\begin{abstract}
Yurinskii's coupling is a popular theoretical tool for non-asymptotic
distributional analysis in mathematical statistics and applied probability,
offering a Gaussian strong approximation with an explicit error bound under
easily verifiable conditions.
Originally stated in $\ell_2$-norm for sums of independent random vectors, it
has recently been extended both to the $\ell_p$-norm,
for $1 \leq p \leq \infty$, and to vector-valued martingales in
$\ell_2$-norm, under some strong conditions.
We present as our main result a Yurinskii coupling for approximate
martingales in $\ell_p$-norm, under substantially weaker conditions than
those previously imposed. Our formulation further allows for the coupling
variable to follow a more general Gaussian mixture distribution, and we
provide a novel third-order coupling method which gives tighter
approximations in certain settings.
We specialize our main result to mixingales, martingales, and independent
data, and derive uniform Gaussian mixture strong approximations for
martingale empirical processes. Applications to nonparametric
partitioning-based and local polynomial regression procedures are provided,
alongside central limit theorems for high-dimensional martingale vectors.
 \end{abstract}

\vspace*{10mm}

\clearpage
\pagebreak

\tableofcontents

\pagebreak

\section{Introduction}

Yurinskii's coupling \citep{yurinskii1978error} has proven to be an important
theoretical tool for developing non-asymptotic distributional approximations in
mathematical statistics and applied probability. For a sum $S$ of $n$
independent zero-mean $d$-dimensional random vectors, this coupling technique
constructs (on a suitably enlarged probability space)
a zero-mean $d$-dimensional
Gaussian vector $T$ which has the same covariance matrix as $S$ and
which is close
to $S$ in probability, bounding the discrepancy $\|S-T\|$ as a function of $n$,
$d$, the choice of norm, and some features of the underlying distribution.
See, for example, \citet[Chapter 10]{pollard2002user} for a textbook
introduction, and \citet{csorgo1981strong} and
\citet{Lindvall_1992_Book} for background references.

When compared to other coupling approaches, such as the celebrated Hungarian
construction \citep{komlos1975approximation} or Zaitsev's coupling
\citep{zaitsev1987estimates,zaitsev1987gaussian}, Yurinskii's approach stands
out for its simplicity, robustness, and wider applicability, while also offering
tighter couplings in some applications (see below for more discussion and
examples). These features have led many scholars to use Yurinskii's coupling to
study the distributional properties of high-dimensional
statistical procedures in
a variety of settings, often with the end goal of developing uncertainty
quantification or hypothesis testing methods. For example, in recent years,
Yurinskii's coupling has been used
to construct Gaussian approximations for the suprema of empirical processes
\citep{chernozhukov2014gaussian};
to establish distribution theory for non-Donsker stochastic $t$-processes
generated in nonparametric series regression \citep{belloni2015some};
to prove distributional approximations for
high-dimensional $\ell_p$-norms \citep{biau2015high};
to develop distribution theory for vector-valued martingales
\citep{belloni2018high,li2020uniform};
to derive a law of the iterated logarithm for stochastic gradient descent
optimization methods \citep{anastasiou2019normal};
to establish uniform distributional results for nonparametric
high-dimensional quantile processes \citep{belloni2019conditional};
to develop distribution theory for non-Donsker stochastic
$t$-processes generated in partitioning-based series regression
\citep{cattaneo2020large};
to deduce Bernstein--von Mises theorems in high-dimensional settings
\citep{ray2021bernstein};
and to develop distribution theory for non-Donsker U-processes based on dyadic
network data \citep{cattaneo2024uniform}.
There are also many other early applications of Yurinskii's coupling:
\citet{dudley1983invariance} and \citet{dehling1983limit} establish invariance
principles for Banach space-valued random variables, and \citet{lecam1988} and
\citet{sheehy1992uniform} obtain uniform Donsker results for empirical
processes, to name just a few.

This paper presents a new Yurinskii coupling which encompasses and
improves upon all of the results previously available in the
literature, offering four new primary features:
\begin{enumerate}[label=(\roman*),leftmargin=*]
  \item
    \label{it:contribution_approximate_martingale}
    It applies to vector-valued
    \textit{approximate martingale} data.
  \item
    \label{it:contribution_gaussian_mixture}
    It allows for a \textit{Gaussian mixture} coupling distribution.
  \item
    \label{it:contribution_degeneracy}
    It imposes \textit{no restrictions on degeneracy} of the
    data covariance matrix.
  \item
    \label{it:contribution_third_order}
    It establishes a \textit{third-order} coupling to
    improve the approximation in certain situations.
\end{enumerate}

Closest to our work are
the recent paper by \citet{li2020uniform}
and the unpublished manuscript by \citet{belloni2018high},
which both investigated distribution theory for martingale
data using Yurinskii's coupling and related methods. Specifically,
\citet{li2020uniform}
established a Gaussian $\ell_2$-norm Yurinskii coupling for
mixingales and martingales under the assumption that the covariance structure
has a minimum eigenvalue bounded away from zero. As formally demonstrated in
this paper (see Section~\ref{sec:kde}), such eigenvalue assumptions can be
prohibitively strong in practically relevant applications. In contrast, our
Yurinskii coupling does not impose any restrictions on covariance degeneracy
\ref{it:contribution_degeneracy}, in addition to offering several other new
features not present in \citet{li2020uniform}, including
\ref{it:contribution_approximate_martingale},
\ref{it:contribution_gaussian_mixture}, \ref{it:contribution_third_order}, and
applicability to general $\ell_p$-norms. In addition, we correct a slight
technical inaccuracy in their proof relating to the derivation of bounds in
probability (see Remark \ref{rem:coupling_bounds_probability}).

\citet{belloni2018high} did not establish a Yurinskii coupling for
martingales, but rather a central limit theorem for smooth functions
of high-dimensional martingales
using the celebrated second-order Lindeberg method
\citep[see][and references therein]{chatterjee2006generalization},
explicitly accounting for covariance
degeneracy. As a consequence, their result could be leveraged to deduce a
Yurinskii coupling for martingales with additional, non-trivial technical work
(see
  \iftoggle{aos}{the supplementary material
  \citep{cattaneo2025yurinskiisupplement}}{Appendix~\ref{sec:proofs}}
for details). Nevertheless, a
Yurinskii coupling derived from \citet{belloni2018high} would not feature
\ref{it:contribution_approximate_martingale},
\ref{it:contribution_gaussian_mixture}, \ref{it:contribution_third_order}, or
general $\ell_p$-norms, as our results do. We discuss
further the connections between our work and the related literature in the
upcoming sections, both when introducing our main theoretical results and when
presenting examples and statistical applications.

The most general coupling result of this paper (Theorem~\ref{thm:sa_dependent})
is presented in Section~\ref{sec:main_results}, where we also specialize it to
a slightly weaker yet more user-friendly formulation
(Proposition~\ref{pro:sa_simplified}). Our Yurinskii coupling for approximate
martingales is a strict generalization of all previous Yurinskii couplings
available in the literature, offering a Gaussian mixture strong approximation
for approximate martingale vectors in $\ell_p$-norm, with an improved rate of
approximation when the third moments of the data are negligible, making no
assumptions on the spectrum of the data covariance matrix. A key technical
innovation underlying the proof of Theorem~\ref{thm:sa_dependent} is that we
explicitly account for the possibility that the minimum eigenvalue of the
variance may be zero, or that its lower bound may be unknown, with the argument
proceeding using a carefully tailored regularization. Establishing a coupling
to a Gaussian mixture distribution is achieved by an appropriate conditioning
argument, leveraging a conditional version of Strassen's theorem
\citetext{\citealp[Theorem~B.2]{chen2020jackknife};
\citealp[Theorem~4]{monrad1991nearby}},
along with some related technical work
detailed in
\iftoggle{aos}{the supplementary material
\citep{cattaneo2025yurinskiisupplement}}{Appendix~\ref{sec:proofs}}.
A third-order coupling is obtained via a
modification of a standard smoothing technique for Borel sets from classical
versions of Yurinskii's coupling
\iftoggle{aos}{%
  (see Lemma~SA.2 in
    the supplementary material
\citep{cattaneo2025yurinskiisupplement})}{%
  (see Lemma~\ref{lem:smooth_approximation}
in the appendix)},
enabling improved approximation errors
whenever third moments are negligible.

In Proposition~\ref{pro:sa_simplified}, we explicitly tune the parameters of
the aforementioned regularization to obtain a simpler, parameter-free version
of Yurinskii's coupling for approximate martingales, again offering Gaussian
mixture coupling distributions and an improved third-order approximation.
This specialization of our main result takes an agnostic approach to potential
singularities in the data covariance matrix and, as such, may be improved in
specific applications where additional knowledge of the covariance structure is
available. Section~\ref{sec:main_results} also presents some further
refinements when additional structure is imposed, deriving Yurinskii couplings
for mixingales, martingales, and independent data as
Corollaries~\ref{cor:sa_mixingale}, \ref{cor:sa_martingale}, and
\ref{cor:sa_indep}, respectively. We take the opportunity to discuss and correct
in Remark~\ref{rem:coupling_bounds_probability} a technical issue which is
often neglected \citep{pollard2002user, li2020uniform} when using Yurinskii's
coupling to derive bounds in probability. Section~\ref{sec:factor} presents a
stylized example portraying the relevance of our main technical results in the
context of canonical factor models, illustrating the importance of each of our
new Yurinskii coupling features \ref{it:contribution_approximate_martingale}--%
\ref{it:contribution_third_order}.

Section~\ref{sec:emp_proc} considers a substantive application of our main
results: strong approximation of martingale empirical processes. We begin with
the motivating example of canonical kernel density estimation, demonstrating
how Yurinskii's coupling can be applied, and showing in
Lemma~\ref{lem:kde_eigenvalue} why it is essential that we do not place any
conditions on the minimum eigenvalue of the variance matrix
\ref{it:contribution_degeneracy}.
We then present a general-purpose strong
approximation for martingale empirical processes in
Proposition~\ref{pro:emp_proc}, combining classical results in the empirical
process literature \citep{van1996weak} with our coupling from
Corollary~\ref{cor:sa_martingale}. This statement appears to be the first of
its kind for martingale data, and when specialized to independent
(and not necessarily identically distributed) data, it is
shown to be superior to the best known %
comparable strong approximation result
available in the literature \citep{berthet2006revisiting}. Our improvement
comes from using Yurinskii's coupling for the $\ell_\infty$-norm, where
\citet{berthet2006revisiting} apply Zaitsev's coupling
\citep{zaitsev1987estimates, zaitsev1987gaussian} with the larger
$\ell_2$-norm.

Section~\ref{sec:nonparametric} further illustrates the applicability of our
results through two examples in nonparametric regression estimation. Firstly,
we deduce strong approximations for partitioning-based least squares series
estimators with time series data, applying Corollary~\ref{cor:sa_martingale}
directly and additionally imposing only a mild mixing condition on the
regressors. We show that our Yurinskii coupling for martingale vectors delivers
the same distributional approximation rate as the best %
known result for
independent data, and discuss how this can be leveraged to yield a feasible
statistical inference procedure. We also show that if the residuals have
vanishing conditional third moment, an improved rate of Gaussian approximation
can be established. Secondly, we deduce a strong approximation for local
polynomial estimators with time series data,
using our result on martingale empirical processes
(Proposition~\ref{pro:emp_proc}) and again imposing a mixing assumption.
Appealing to empirical process theory is essential here as, in contrast with
series estimators, local polynomials do not possess certain additive
separability properties. The bandwidth restrictions we require are relatively
mild, and, as far as we know, they have not been improved upon even with
independent data.

Section \ref{sec:conclusion} concludes the paper.
Appendix~\ref{sec:high_dim_clt} demonstrates
how our coupling results can be used to derive
distributional Gaussian approximations
(central limit theorems)
for possibly high-dimensional martingale vectors
(Proposition~\ref{pro:clt}).
This result complements a recent literature on probability
and statistics studying
the same problem but with independent data
\citep[see][and references
therein]{buzun2022strong,lopes2022central,%
chernozhukov2023nearly,kock2024remark}.
We also present a version of this result employing
a covariance estimator (Proposition~\ref{pro:bootstrap}),
enabling the construction of valid high-dimensional
confidence sets via a Gaussian multiplier bootstrap.
\iftoggle{aos}{}{%
  Finally we present
  some further results on applications of our theory to
  deriving distributional approximations for
  $\ell_p$-norms of high-dimensional martingale vectors
  in Appendix~\ref{sec:lp}.
}

All proofs are collected in
\iftoggle{aos}{the supplementary material
\citep{cattaneo2025yurinskiisupplement}}{Appendix~\ref{sec:proofs}},
where we also include other technical lemmas of
potential independent interest%
\iftoggle{aos}{%
  , alongside
  some further results on
  distributional approximations for
  $\ell_p$-norms of high-dimensional martingale vectors%
}{}.

\subsection{Notation}

We write $\|x\|_p$ for $p\in[1,\infty]$ to denote the $\ell_p$-norm
if $x$ is a
(possibly random) vector or the induced operator
$\ell_p$--$\ell_p$-norm if $x$
is a matrix. For $X$ a real-valued random variable and an Orlicz function
$\psi$, we use $\vvvert X \vvvert_\psi$ to denote the Orlicz $\psi$-norm
\citep[Section~2.2]{van1996weak} and $\vvvert X \vvvert_p$ for the $L^p(\P)$
norm where $p\in [1,\infty]$. For a matrix $M$, we write
$\|M\|_{\max}$ for the
maximum absolute entry and $\|M\|_\rF$ for the Frobenius norm. We denote
positive semi-definiteness by $M \succeq 0$ and write $I_d$ for the $d \times
d$ identity matrix.

For scalar sequences $x_n$ and $y_n$, we write $x_n \lesssim y_n$ if there
exists a positive constant $C$ such that $|x_n| \leq C |y_n|$ for sufficiently
large $n$. We write $x_n \asymp y_n$ to indicate both $x_n \lesssim y_n$ and
$y_n \lesssim x_n$. Similarly, for random variables $X_n$ and $Y_n$, we write
$X_n \lesssim_\P Y_n$ if for every $\varepsilon > 0$ there exists a positive
constant $C$ such that $\P(|X_n| \geq C |Y_n|) \leq \varepsilon$, and write
$X_n \to_\P X$ for limits in probability. For real numbers $a$ and $b$ we use
$a \vee b = \max\{a,b\}$. We write $\kappa \in \N^d$ for a multi-index, where
$d \in \N = \{0, 1, 2, \ldots\}$, and define
$|\kappa| = \sum_{j=1}^d \kappa_j$, along with
$\kappa! = \prod_{j=1}^{d} \kappa_j !$, and
$x^\kappa = \prod_{j=1}^d x_j^{\kappa_j}$ for $x \in \R^d$.

Since our results concern couplings, some statements must be made on a new or
enlarged probability space. We omit the details of this for clarity of
notation, but technicalities are handled by the Vorob'ev--Berkes--Philipp
Theorem \citep[Theorem~1.1.10]{dudley1999uniform}.

\section{Main results}
\label{sec:main_results}

We begin with our most general result: an $\ell_p$-norm Yurinskii
coupling for a
sum of vector-valued approximate martingale differences to a Gaussian
mixture-distributed random vector. The general result is presented in
Theorem~\ref{thm:sa_dependent}, while
Proposition~\ref{pro:sa_simplified} gives
a simplified and slightly weaker version which is easier to use in many
applications. We then further specialize
Proposition~\ref{pro:sa_simplified} to
three scenarios with successively stronger assumptions, namely mixingales,
martingales, and independent data, in Corollaries~\ref{cor:sa_mixingale},
\ref{cor:sa_martingale}, and \ref{cor:sa_indep} respectively. In each case we
allow for possibly random quadratic variations (cf.\ mixing convergence),
thereby establishing Gaussian mixture couplings in the general setting. In
Remark~\ref{rem:coupling_bounds_probability} we comment on and
correct an often
overlooked technicality relating to the derivation of bounds in probability
from Yurinskii's coupling. As a first illustration of the power of our
generalized $\ell_p$-norm Yurinskii coupling, we present in
Section~\ref{sec:factor} a simple factor model example relating to
all three of
the aforementioned scenarios, discussing further how our contributions are
related to the existing literature.

\begin{theorem}[Strong approximation for vector-valued approximate
  martingales]%
  \label{thm:sa_dependent}

  Take a complete probability space with a countably generated
  filtration $\cH_0, \ldots, \cH_n$
  for some $n \geq 1$, supporting the
  $\R^d$-valued square-integrable random vectors
  $X_1, \ldots, X_n$.
  Let $S = \sum_{i=1}^n X_i$ and
  define
  \begin{align*}
    \tilde X_i
    &=
    \sum_{r=1}^n \big(\E[X_{r} \mid \cH_{i}] - \E[X_{r} \mid \cH_{i-1}]\big)
    & &\text{and}
    &U
    &=
    \sum_{i=1}^{n} \big(
      X_i - \E[ X_i \mid \cH_n]
      + \E[ X_i \mid \cH_0 ]
    \big).
  \end{align*}
  Let $V_i = \Var[\tilde X_i \mid \cH_{i-1}]$ and
  define $\Omega = \sum_{i=1}^n V_i - \Sigma$
  where $\Sigma$ is an almost surely positive semi-definite $\cH_0$-measurable
  $d \times d$ random matrix.
  Then, for each $\eta > 0$ and $p \in [1,\infty]$,
  there exists, on an enlarged probability space,
  an $\R^d$-valued random vector $T$ with
  $T \mid \cH_0 \sim \cN(0, \Sigma)$ such that
  \begin{align}
    \label{eq:sa_dependent}
    \P\big(\|S-T\|_p > 6\eta\big)
    &\leq
    \inf_{t>0}
    \left\{
      2 \P\big( \|Z\|_p > t \big)
      + \min\left\{
        \frac{\beta_{p,2} t^2}{\eta^3},
        \frac{\beta_{p,3} t^3}{\eta^4}
        + \frac{\pi_3 t^3}{\eta^3}
      \right\}
    \right\} \nonumber \\
    &\quad+
    \inf_{M \succeq 0}
    \Big\{ 2 \P\big(\Omega \npreceq M\big) + \delta_p(M,\eta)
    + \varepsilon_p(M, \eta)\Big\}
    +\P\big(\|U\|_p>\eta\big),
  \end{align}
  where $Z, Z_1,\dots ,Z_n$ are i.i.d.\ standard Gaussian random variables on
  $\R^d$ independent of $\cH_n$,
  the second infimum is taken over all positive semi-definite
  $d \times d$ non-random matrices $M$,
  \begin{align*}
    \beta_{p,k}
    &=
    \sum_{i=1}^n \E\left[\| \tilde X_i \|^k_2 \| \tilde X_i \|_p
    + \|V_i^{1/2} Z_i \|^k_2 \|V_i^{1/2} Z_i \|_p \right],
    &\pi_3
    &=
    \sum_{i=1}^{n}
    \sum_{|\kappa| = 3}
    \E \Big[ \big|
      \E [ \tilde X_i^\kappa \mid \cH_{i-1} ]
    \big| \Big]
  \end{align*}
  for $k \in \{2, 3\}$, with $\pi_3 = \infty$ if the associated
  conditional expectation does not exist, and with
  \begin{align*}
    \delta_p(M,\eta)
    &=
    \P\left(
      \big\|\big((\Sigma +M)^{1/2}- \Sigma^{1/2}\big) Z\big\|_p
      \geq \eta
    \right), \\
    \varepsilon_p(M, \eta)
    &=
    \P\left(\big\| (M - \Omega)^{1/2} Z \big\|_p\geq \eta, \
    \Omega \preceq M\right).
  \end{align*}
\end{theorem}

This theorem offers four novel contributions to the literature on
coupling theory and strong approximation,
as discussed in the introduction.
Firstly \ref{it:contribution_approximate_martingale}, it allows for
approximate
vector-valued martingales, with the variables $\tilde X_i$ forming martingale
differences with respect to $\cH_i$ by construction, and $U$ quantifying the
associated martingale approximation error. Such martingale approximation
techniques for sequences of dependent random vectors are well established and
have been used in a range of scenarios: see, for example,
\citet{wu2004martingale}, \citet{wu2005nonlinear},
\citet{dedecker2007weak}, \citet{zhao2008martingale},
\citet{peligrad2010conditional}, \citet{atchade2014martingale},
\citet{cuny2014martingale}, \citet{magda2018martingale}, and references
therein. In Section~\ref{sec:mixingales} we demonstrate how this approximation
can be established in practice by restricting our general theorem to the
special case of mixingales, while the upcoming example in
Section~\ref{sec:factor} provides an illustration in the context of
auto-regressive factor models.

Secondly \ref{it:contribution_gaussian_mixture},
Theorem~\ref{thm:sa_dependent} allows for the
resulting coupling variable $T$
to follow a multivariate Gaussian distribution only conditionally,
and thus we offer a useful analog of mixing convergence in the context
of strong approximation.
To be more precise, the random matrix $\sum_{i=1}^{n} V_i$
is the quadratic variation of the constructed martingale
$\sum_{i=1}^n \tilde X_i$, and we approximate it using the $\cH_0$-measurable
random matrix $\Sigma$. This yields the coupling variable
$T \mid \cH_0 \sim \cN(0, \Sigma)$, which can alternatively be written as
$T=\Sigma^{1/2} Z$ with $Z \sim \cN(0,I_d)$ independent of $\cH_0$.
The errors in this quadratic variation
approximation are accounted for by the terms
$\P(\Omega \npreceq M)$, $\delta_p(M, \eta)$ and $\varepsilon_p(M, \eta)$,
utilizing a regularization argument through the free matrix parameter $M$.
If a non-random $\Sigma$ is used, then $T$ is unconditionally Gaussian,
and one can take $\cH_0$ to be the trivial $\sigma$-algebra.
As demonstrated in our proof, our approach to establishing a
mixing approximation is different from naively taking an unconditional version
of Yurinskii's coupling and applying
it conditionally on $\cH_0$, which will not deliver the same coupling as in
Theorem~\ref{thm:sa_dependent} for a few reasons.
To begin with, we explicitly indicate in the
conditions of Theorem~\ref{thm:sa_dependent} where conditioning is required.
Next, our error of approximation is given unconditionally,
involving only marginal expectations and probabilities.
Finally, we provide a rigorous account of the construction of the
conditionally Gaussian coupling variable $T$ via a conditional version
of Strassen's theorem \citetext{\citealp[Theorem~B.2]{chen2020jackknife};
\citealp[Theorem~4]{monrad1991nearby}}.
Section~\ref{sec:martingales}
illustrates how a strong approximation akin to
mixing convergence can arise when the data
forms an exact martingale, and Section~\ref{sec:factor} gives a simple example
relating to factor modeling in statistics and data science.

As a third contribution to the literature
\ref{it:contribution_degeneracy}, and
of particular importance for applications,
Theorem~\ref{thm:sa_dependent} makes
no requirements on the minimum eigenvalue of the quadratic variation of the
approximating martingale sequence. Instead, our proof technique employs a
careful regularization scheme designed to account for any such exact or
approximate rank degeneracy in $\Sigma$. This capability is
fundamental in some
applications, a fact which we illustrate in Section \ref{sec:kde} by
demonstrating the significant improvements in strong approximation errors
delivered by Theorem~\ref{thm:sa_dependent} relative to those obtained using
prior results in the literature.

Finally \ref{it:contribution_third_order},
Theorem \ref{thm:sa_dependent} gives
a third-order strong approximation alongside the usual second-order
version considered in all prior literature.
More precisely, we observe that an analog of the term
$\beta_{p,2}$ is present in the
classical Yurinskii coupling and comes from a Lindeberg
telescoping sum argument,
replacing random variables by Gaussians with the same mean
and variance to match the first and second moments.
Whenever the third conditional moments of $\tilde X_i$ are negligible
(quantified by $\pi_3$), this moment-matching argument can be extended to
third-order terms, giving a new quantity $\beta_{p,3}$.
At this level of generality, it is not possible to obtain explicit
bounds on $\pi_3$ because we make no assumptions on the
relationship between the
data $X_i$ and the $\sigma$-algebras $\cH_i$
(and therefore the variables $\tilde X_i$ resulting from the
martingale approximation).
However, if $X_1, \ldots, X_n$ form martingale differences
with respect to $\cH_0, \ldots, \cH_n$,
then $\tilde X_i = X_i$ almost surely (see Section~\ref{sec:martingales}).
In this setting, assuming that
$\E \big[ X_i^\kappa \mid \cH_{i-1} \big] = 0$
for each multi-index $\kappa$ with $|\kappa| = 3$
(e.g.\ if the data is conditionally symmetrically distributed around zero),
then using $\beta_{p,3}$ rather than $\beta_{p,2}$
can give smaller coupling approximation errors in \eqref{eq:sa_dependent}.
Such a refinement can be viewed as a strong approximation counterpart
to classical Edgeworth expansion methods,
and we illustrate this phenomenon in our
upcoming applications to nonparametric inference
(Section~\ref{sec:nonparametric}).

\subsection{User-friendly formulation of the main result}%

The result in Theorem~\ref{thm:sa_dependent} is given in a somewhat implicit
manner, involving infima over the free parameters $t > 0$ and $M \succeq 0$,
and it is not clear how to compute these in general. In the upcoming
Proposition~\ref{pro:sa_simplified}, we set $M = \nu^2 I_d$ and approximately
optimize over $t > 0$ and $\nu > 0$, resulting in a simplified and slightly
weaker version of our main general result. In specific applications, where
there is additional knowledge of the quadratic variation structure, other
choices of regularization schemes may be more appropriate. Nonetheless, the
choice $M = \nu^2 I_d$ leads to arguably the principal result of our work,
due to its simplicity and utility in statistical applications. For
convenience,
define the functions $\phi_p : \{1, 2, \ldots\} \to \R$,
for $p \in [0, \infty]$, by
\begin{align*}
  \phi_p(d) =
  \begin{cases}
    \sqrt{pd^{2/p} }  & \text{ if } p \in [1,\infty),\\
    \sqrt{2\log 2d}   & \text{ if } p =\infty.
  \end{cases}
\end{align*}
With $Z \sim \cN(0, I_d)$ and $t > 0$,
these functions satisfy
$\P( \|Z\|_p > t ) \leq \E[\|Z\|_p] / t \leq \phi_p(d) / t$
(see \iftoggle{aos}{%
    Lemma~SA.4
    in the supplementary material
  \citep{cattaneo2025yurinskiisupplement}}{%
Lemma~\ref{lem:gaussian_pnorm} in the appendix}).

\begin{proposition}[Simplified strong approximation for
  vector-valued approximate martingales]%
  \label{pro:sa_simplified}

  Assume the setup and notation of Theorem~\ref{thm:sa_dependent}.
  For each $\eta > 0$ and $p \in [1,\infty]$,
  there exists a random vector $T \mid \cH_0 \sim \cN(0, \Sigma)$ satisfying
  \begin{align*}
    \P\big(\|S-T\|_p > \eta\big)
    &\leq
    24 \left(
      \frac{\beta_{p,2} \phi_p(d)^2}{\eta^3}
    \right)^{1/3}
    + 17 \left(
      \frac{\E \left[ \|\Omega\|_2 \right] \phi_p(d)^2}{\eta^2}
    \right)^{1/3}
    +\P\left(\|U\|_p>\frac{\eta}{6}\right).
  \end{align*}
  If further $\pi_3 = 0$, then also
  \begin{align*}
    \P\big(\|S-T\|_p > \eta\big)
    &\leq
    24 \left(
      \frac{\beta_{p,3} \phi_p(d)^3}{\eta^4}
    \right)^{1/4}
    + 17 \left(
      \frac{\E \left[ \|\Omega\|_2 \right] \phi_p(d)^2}{\eta^2}
    \right)^{1/3}
    +\P\left(\|U\|_p>\frac{\eta}{6}\right).
  \end{align*}
\end{proposition}

Proposition~\ref{pro:sa_simplified} makes clear the potential benefit of a
third-order coupling when $\pi_3 = 0$, as in this case the bound features
$\beta_{p,3}^{1/4}$ rather than $\beta_{p,2}^{1/3}$. If $\pi_3$ is small but
non-zero, an analogous result can easily be derived by adjusting the optimal
choices of $t$ and $\nu$, but we omit this for clarity of notation. In
applications (see Section~\ref{sec:series}), this reduction of the
exponent can
provide a significant improvement in terms of the dependence of the bound on
the sample size $n$, the dimension $d$, and other problem-specific quantities.
When using our results for strong approximation, it is usual to set
$p = \infty$ to bound the maximum discrepancy over the entries of a vector (to
construct uniform confidence sets, for example). In this setting, we have that
$\phi_\infty(d) = \sqrt{2 \log 2d}$ has a sub-Gaussian slow-growing dependence
on the dimension. The remaining term depends on $\E[\|\Omega\|_2]$
and requires
that the matrix $\Sigma$ be a good approximation of $\sum_{i=1}^{n}
V_i$, while
remaining $\cH_0$-measurable. In some applications (such as factor modeling;
see Section~\ref{sec:factor}), it can be shown that the quadratic variation
$\sum_{i=1}^n V_i$ remains random and $\cH_0$-measurable even in
large samples,
giving a natural choice for $\Sigma$.

In the next few sections, we continue to refine
Proposition~\ref{pro:sa_simplified}, presenting a sequence of results with
increasingly strict assumptions on the dependence structure of the data $X_i$.
These allow us to demonstrate the broad applicability of our main results,
providing more explicit bounds in settings which are likely to be of special
interest. In particular, we consider mixingales, martingales, and independent
data, comparing our derived results with those in the existing literature.

\subsection{Mixingales}\label{sec:mixingales}

In our first refinement, we provide a natural method for bounding the
martingale approximation error term $U$. Suppose that $X_i$ form an
$\ell_p$-mixingale in $L^1(\P)$ in the sense that there exist non-negative
$c_1, \ldots, c_n$ and $\zeta_0, \ldots, \zeta_n$ such that for all
$1 \leq i \leq n$ and $0 \leq r \leq i$,
\begin{align}
  \label{eq:mixingale_1}
  \E \left[ \left\|
    \E \left[ X_i \mid \cH_{i-r} \right]
  \right\|_p \right]
  &\leq
  c_i \zeta_r,
\end{align}
and for all $1 \leq i \leq n$ and $0 \leq r \leq n-i$,
\begin{align}
  \label{eq:mixingale_2}
  \E \left[ \big\|
    X_i - \E \big[ X_i \mid \cH_{i+r} \big]
  \big\|_p \right]
  &\leq
  c_i \zeta_{r+1}.
\end{align}
These conditions are satisfied, for example, if $X_i$ are integrable strongly
$\alpha$-mixing random variables \citep{mcleish1975invariance}, or
if $X_i$ are
generated by an auto-regressive or auto-regressive moving average process (see
Section~\ref{sec:factor}), among many other possibilities
\citep{bradley2005basic}. Then, in the notation of
Theorem~\ref{thm:sa_dependent}, we have by Markov's inequality that
\begin{align*}
  \P \left( \|U\|_p > \frac{\eta}{6} \right)
  &\leq
  \frac{6}{\eta}
  \sum_{i=1}^{n}
  \E \left[
    \big\|
    X_i - \E \left[ X_i \mid \cH_n \right]
    \big\|_p
    + \big\|
    \E \left[ X_i \mid \cH_0 \right]
    \big\|_p
  \right]
  \leq \frac{\zeta}{\eta},
\end{align*}
with $\zeta = 6 \sum_{i=1}^{n} c_i (\zeta_{i} + \zeta_{n-i+1})$.
Combining Proposition~\ref{pro:sa_simplified} with this
martingale error bound yields the following result for mixingales.
\begin{corollary}[Strong approximation for vector-valued mixingales]%
  \label{cor:sa_mixingale}

  Assume the setup and notation of Theorem~\ref{thm:sa_dependent}, and suppose
  that the mixingale conditions \eqref{eq:mixingale_1} and
  \eqref{eq:mixingale_2} hold. For each $\eta > 0$ and $p \in
  [1,\infty]$ there
  exists a random vector $T \mid \cH_0 \sim \cN(0, \Sigma)$ satisfying
  \begin{align*}
    \P\big(\|S-T\|_p > \eta\big)
    &\leq
    24 \left(
      \frac{\beta_{p,2} \phi_p(d)^2}{\eta^3}
    \right)^{1/3}
    + 17 \left(
      \frac{\E \left[ \|\Omega\|_2 \right] \phi_p(d)^2}{\eta^2}
    \right)^{1/3}
    + \frac{\zeta}{\eta}.
  \end{align*}
  If further $\pi_3 = 0$ then
  \begin{align*}
    \P\big(\|S-T\|_p > \eta\big)
    &\leq
    24 \left(
      \frac{\beta_{p,3} \phi_p(d)^3}{\eta^4}
    \right)^{1/4}
    + 17 \left(
      \frac{\E \left[ \|\Omega\|_2 \right] \phi_p(d)^2}{\eta^2}
    \right)^{1/3}
    + \frac{\zeta}{\eta}.
  \end{align*}
\end{corollary}

The closest antecedent to Corollary~\ref{cor:sa_mixingale} is found in
\citet[Theorem~4]{li2020uniform}, who also considered Yurinskii's coupling for
mixingales. Our result improves on this work in the following manner: it
removes any requirements on the minimum eigenvalue of the quadratic variation
of the mixingale sequence; it allows for general $\ell_p$-norms with
$p\in[1,\infty]$; it establishes a coupling to a multivariate Gaussian
mixture distribution in general; and it permits third-order couplings
(when $\pi_3=0$). These improvements have important practical implications as
demonstrated in Section~\ref{sec:factor} and Section~\ref{sec:nonparametric},
where significantly better coupling approximation
errors are demonstrated for a variety of statistical applications. On the
technical side, our result is rigorously established using a conditional
version of Strassen's theorem, a
carefully crafted
regularization argument, and a third-order Lindeberg method.
Furthermore (Remark~\ref{rem:coupling_bounds_probability}), we
clarify a technical issue in \citet{li2020uniform} surrounding the derivation
of valid probability bounds for $\|S-T\|_p$.

Corollary~\ref{cor:sa_mixingale} focused on mixingales for simplicity, but, as
previously discussed, any method for constructing a martingale approximation
$\tilde X_i$ and bounding the resulting error $U$ could be used instead in
Proposition~\ref{pro:sa_simplified} to derive a similar result.

\subsection{Martingales}\label{sec:martingales}

For our second refinement, suppose that
$X_i$ form martingale differences with respect to $\cH_i$.
In this case, $\E[X_i \mid \cH_n] = X_i$ and $\E[X_i \mid \cH_0] = 0$,
so $U = 0$, and the martingale approximation error term vanishes.
Applying Proposition~\ref{pro:sa_simplified} in this setting
directly yields the following result.
\begin{corollary}[Strong approximation for vector-valued martingales]%
  \label{cor:sa_martingale}

  With the setup and notation of Theorem~\ref{thm:sa_dependent}, suppose $X_i$
  is $\cH_i$-measurable with $\E[X_i \mid \cH_{i-1}] = 0$ for
  $1 \leq i \leq n$. Then, for each $\eta > 0$ and $p \in
  [1,\infty]$, there is
  a random vector $T \mid \cH_0 \sim \cN(0, \Sigma)$ with
  \begin{align}
    \label{eq:sa_martingale_order_2}
    \P\big(\|S-T\|_p > \eta\big)
    &\leq
    24 \left(
      \frac{\beta_{p,2} \phi_p(d)^2}{\eta^3}
    \right)^{1/3}
    + 17 \left(
      \frac{\E \left[ \|\Omega\|_2 \right] \phi_p(d)^2}{\eta^2}
    \right)^{1/3}.
  \end{align}
  If further $\pi_3 = 0$ then
  \begin{align}
    \label{eq:sa_martingale_order_3}
    \P\big(\|S-T\|_p > \eta\big)
    &\leq
    24 \left(
      \frac{\beta_{p,3} \phi_p(d)^3}{\eta^4}
    \right)^{1/4}
    + 17 \left(
      \frac{\E \left[ \|\Omega\|_2 \right] \phi_p(d)^2}{\eta^2}
    \right)^{1/3}.
  \end{align}
\end{corollary}

The closest antecedents to Corollary~\ref{cor:sa_martingale} are
\citet{belloni2018high} and \citet{li2020uniform}, who also (implicitly or
explicitly) considered Yurinskii's coupling for martingales. More
specifically,
\citet[Theorem 1]{li2020uniform} established an explicit
$\ell_2$-norm Yurinskii coupling
for martingales under a strong assumption on the minimum eigenvalue of the
martingale quadratic variation, while \citet[Theorem 2.1]{belloni2018high}
established a
central limit theorem for vector-valued martingale sequences employing the
standard second-order Lindeberg method. As such, their proof could be
adapted to deduce a Yurinskii coupling for martingales with the help of a
conditional version of Strassen's theorem
and some additional nontrivial technical work.

Corollary~\ref{cor:sa_martingale} improves over this prior work as follows.
With respect to \citet{li2020uniform}, our result establishes an $\ell_p$-norm
Gaussian mixture Yurinskii coupling for martingales without any
requirements on
the minimum eigenvalue of the martingale quadratic variation, and permits a
third-order coupling if $\pi_3=0$. The first probability bound
\eqref{eq:sa_martingale_order_2} in
Corollary~\ref{cor:sa_martingale} gives the
same rate of strong approximation as that in Theorem~1 of
\citet{li2020uniform}
when $p=2$, with non-random $\Sigma$, and when the eigenvalues of a normalized
version of $\Sigma$ are bounded away from zero. In Section~\ref{sec:kde} we
demonstrate the crucial importance of removing this eigenvalue lower bound
restriction in applications involving nonparametric kernel
estimators, while in
Section~\ref{sec:series} we demonstrate how the availability of a third-order
coupling \eqref{eq:sa_martingale_order_3} can give improved
approximation rates
in applications involving nonparametric series estimators with conditionally
symmetrically distributed residual errors. Finally, our technical
work improves
on \citet{li2020uniform} in two respects: (i) we employ a conditional version
of Strassen's theorem (see
  \iftoggle{aos}{%
    Lemma~SA.1 in
    the supplementary material
  \citep{cattaneo2025yurinskiisupplement}}{Lemma~\ref{lem:strassen}
in the appendix})
to appropriately handle the conditioning arguments; and (ii) we deduce valid
probability bounds for $\|S-T\|_p$, as the following
Remark~\ref{rem:coupling_bounds_probability} makes clear.

\begin{remark}[Yurinskii's coupling and bounds in probability]
  \label{rem:coupling_bounds_probability}
  Given a sequence of random vectors $S_n$, Yurinskii's method provides a
  coupling in the following form: for each $n$ and any $\eta > 0$,
  there exists
  a random vector $T_n$ with $\P\big(\|S_n - T_n\| > \eta\big) < r_n(\eta)$,
  where $r_n(\eta)$ is the approximation error. Crucially, each coupling
  variable $T_n$ is a function of the desired approximation level $\eta$ and,
  as such, deducing bounds in probability on $\|S_n - T_n\|$ requires some
  extra care. One option is to select a sequence $R_n \to \infty$
  and note that
  $\P\big(\|S_n - T_n\| > r_n^{-1}(1 / R_n)\big) < 1 / R_n \to 0$ and hence
  $\|S_n - T_n\| \lesssim_\P r_n^{-1}(1 / R_n)$. In this case,
  $T_n$ depends on
  the choice of $R_n$, which can in turn typically be chosen to diverge slowly
  enough to cause no issues in applications.
\end{remark}

Technicalities akin to those outlined in
Remark~\ref{rem:coupling_bounds_probability} have been both addressed and
neglected alike in the prior literature. \citet[Chapter 10.4, Example
16]{pollard2002user} apparently misses this subtlety, providing an
inaccurate bound in probability based on the Yurinskii coupling.
\citet{li2020uniform} seem to make the same mistake in the proof of their
Lemma~A2, which invalidates the conclusion of their Theorem~1. In contrast,
\citet{belloni2015some} and \citet{belloni2019conditional} directly provide
bounds in $o_\P$ instead of $O_\P$, circumventing these issues in a manner
similar to our approach involving a diverging sequence $R_n$.

To see how this phenomenon applies to our main results, observe that the
second-order martingale coupling given as \eqref{eq:sa_martingale_order_2} in
Corollary~\ref{cor:sa_martingale} implies that for any $R_n \to \infty$,
\begin{align*}
  \|S - T\|_p
  \lesssim_\P
  \beta_{p,2}^{1/3}
  \phi_p(d)^{2/3} R_n
  + \E[\|\Omega\|_2]^{1/2}
  \phi_p(d) R_n.
\end{align*}
This bound is comparable to that obtained by \citet[Theorem~1]{li2020uniform}
with $p=2$, albeit with their formulation missing the $R_n$ correction terms.
In Section~\ref{sec:series} we discuss further their (amended) result, in the
setting of nonparametric series estimation. Our approach using
$p = \infty$ obtains superior distributional approximation rates, alongside
exhibiting various other improvements such as the aforementioned third-order
coupling.

Turning to the comparison with \citet{belloni2018high}, our
Corollary~\ref{cor:sa_martingale} again offers the same improvements, with the
only exception being that the authors did account for the implications of a
possibly vanishing minimum eigenvalue. However, their results exclusively
concern high-dimensional central limit theorems for vector-valued martingales,
and therefore while their findings
could in principle enable the derivation of a result similar to our
Corollary~\ref{cor:sa_martingale}, this would require additional
technical work
on their behalf in multiple ways
(see
  \iftoggle{aos}{the supplementary material
\citep{cattaneo2025yurinskiisupplement}}{Appendix~\ref{sec:proofs}}):
(i) a correct application of a conditional
version of Strassen's theorem
\iftoggle{aos}{%
  (Lemma~SA.1 in
    the supplementary material
\citep{cattaneo2025yurinskiisupplement})}{(Lemma~\ref{lem:strassen}
in the appendix)};
(ii) the development of a third-order Borel set smoothing technique and
associated $\ell_p$-norm moment control
\iftoggle{aos}{%
(Lemmas~SA.2, SA.3, and SA.4)}{(Lemmas~\ref{lem:smooth_approximation},
\ref{lem:gaussian_useful}, and \ref{lem:gaussian_pnorm})};
(iii) a careful truncation scheme to account for
$\Omega\npreceq0$; and (iv) a valid third-order Lindeberg argument
\iftoggle{aos}{%
(Lemma~SA.8)}{(Lemma~\ref{lem:sa_martingale})};
among others.

\subsection{Independence}

As a final refinement, suppose that $X_i$ are independent and
zero-mean conditionally on $\cH_0$,
and take $\cH_i$ to be the filtration
generated by $X_1, \ldots, X_i$ and $\cH_0$ for $1 \leq i \leq n$.
Then, taking $\Sigma = \sum_{i=1}^n V_i$
gives $\Omega = 0$, and hence Corollary~\ref{cor:sa_martingale}
immediately yields the following result.
\begin{corollary}[Strong approximation for sums of independent vectors]%
  \label{cor:sa_indep}

  Assume the setup of Theorem~\ref{thm:sa_dependent},
  and suppose $X_i$ are independent given $\cH_0$,
  with $\E[X_i \mid \cH_0] = 0$.
  Then, for each $\eta > 0$ and $p \in [1,\infty]$,
  with $\Sigma = \sum_{i=1}^n V_i$,
  there exists $T \mid \cH_0 \sim \cN(0, \Sigma)$ satisfying
  \begin{align}
    \label{eq:sa_indep_order_2}
    \P\big(\|S-T\|_p > \eta\big)
    &\leq
    24 \left(
      \frac{\beta_{p,2} \phi_p(d)^2}{\eta^3}
    \right)^{1/3}.
  \end{align}
  If further $\pi_3 = 0$ then
  \begin{align*}
    \P\big(\|S-T\|_p > \eta\big)
    &\leq
    24 \left(
      \frac{\beta_{p,3} \phi_p(d)^3}{\eta^4}
    \right)^{1/4}.
  \end{align*}
\end{corollary}

Taking $\cH_0$ to be trivial, the first inequality \eqref{eq:sa_indep_order_2}
in Corollary~\ref{cor:sa_indep} provides an $\ell_p$-norm approximation
analogous to that presented in \cite{belloni2019conditional}. By further
restricting to $p=2$, we recover the original Yurinskii coupling as presented
in \citet[Theorem~1]{lecam1988} and \citet[Theorem~10]{pollard2002user}. Thus,
in the independent data setting, our result improves on prior work as follows:
(i) it establishes a coupling to a multivariate Gaussian mixture distribution;
and (ii) it permits a third-order coupling if $\pi_3=0$.

\subsection{Stylized example: factor modeling}
\label{sec:factor}

In this section, we present a simple statistical example of how our
improvements over prior coupling results can have important theoretical and
practical implications. Consider the stylized factor model
\begin{align*}
  X_i = L f_i + \varepsilon_i, \qquad 1 \leq i \leq n,
\end{align*}
with random variables $L$ taking values in $\R^{d \times m}$, $f_i$ in $\R^m$,
and $\varepsilon_i$ in $\R^d$. We interpret $f_i$ as a latent factor variable
and $L$ as a random factor loading, with
independent (idiosyncratic) disturbances
$(\varepsilon_1, \ldots, \varepsilon_n)$.
See \citet{fan2020statistical}, and references therein, for a
textbook review of factor analysis in statistics and econometrics.

We employ the above factor model to give a first illustration of the
applicability of our main result Theorem~\ref{thm:sa_dependent}, the
user-friendly Proposition~\ref{pro:sa_simplified}, and their specialized
Corollaries~\ref{cor:sa_mixingale}--\ref{cor:sa_indep}. We
consider three different sets of conditions to demonstrate the
applicability of
each of our corollaries for mixingales, martingales, and independent data,
respectively. We assume throughout that each
$\varepsilon_i$ is zero-mean and finite variance, and
that $(\varepsilon_1, \ldots, \varepsilon_n)$ is independent
of $L$ and $(f_1, \ldots, f_n)$. Let $\cH_i$ be the $\sigma$-algebra generated
by $L$, $(f_1, \ldots, f_i)$ and $(\varepsilon_1, \ldots,
\varepsilon_i)$, with
$\cH_0$ the $\sigma$-algebra generated by $L$ alone.

\begin{enumerate}[label=(\roman*)]
  \item \emph{Independent data}.
    Suppose that the factors $(f_1, \ldots,
    f_n)$ are independent conditional on $L$ and satisfy
    $\E [ f_i \mid L ] = 0$.
    Then, since $X_i$ are independent conditional on $\cH_0$ and with
    $\E [ X_i \mid \cH_0 ] = \E [ L f_i + \varepsilon_i \mid L ] = 0$,
    we can apply Corollary~\ref{cor:sa_indep} to $\sum_{i=1}^n X_i$.
    In general, we will obtain a coupling variable which has the Gaussian
    mixture distribution $T \mid \cH_0 \sim \cN(0, \Sigma)$ where
    $\Sigma= \sum_{i=1}^n (L\Var[f_i \mid L]L^\T +\Var[\varepsilon_i])$.
    In the special case where $L$ is non-random
    and $\cH_0$ is trivial, the coupling is Gaussian. Furthermore,
    if $f_i\mid L$ and $\varepsilon_i$ are symmetric about zero
    and bounded almost surely, then $\pi_3=0$, and the coupling is improved.

  \item \emph{Martingales}.
    Suppose instead that we assume only a martingale
    condition on the latent factor variables so that
    $\E \left[ f_i \mid L, f_1, \ldots, f_{i-1} \right] = 0$.
    Then $\E [ X_i \mid \cH_{i-1} ]
    = L\, \E \left[ f_i \mid \cH_{i-1} \right] = 0$
    and Corollary~\ref{cor:sa_martingale} is applicable to $\sum_{i=1}^n X_i$.
    The preceding comments on Gaussian mixture distributions
    and third-order couplings continue to apply.

  \item \emph{Mixingales}.
    Finally, assume that the factors follow the
    auto-regressive model $f_i = A f_{i-1} + u_i$ where
    $A \in \R^{m \times m}$ is non-random and $(u_1, \ldots, u_n)$ are
    zero-mean, independent, and independent of
    $(\varepsilon_1, \ldots, \varepsilon_n)$.
    Then $\E \left[ f_i \mid f_0 \right] = A^i f_0$, so taking
    $p \in [1, \infty]$ we see that
    $\E \big[ \| \E [ f_i \mid f_0 ] \|_p \big]
    = \E \big[ \| A^i f_0 \|_p \big] \leq \|A\|_p^i\,\E [ \|f_0\|_p ]$,
    and that clearly $f_i - \E [ f_i \mid \cH_n ] = 0$.
    Thus, whenever $\|A\|_p < 1$, the geometric sum formula implies that
    the mixingale result from Corollary~\ref{cor:sa_mixingale} applies to
    $\sum_{i=1}^n X_i$. The conclusions on Gaussian mixture distributions
    and third-order couplings parallel the previous cases.
\end{enumerate}

This simple application to factor modeling gives a preliminary illustration of
the power of our main results, encompassing settings which could
not be handled
by employing Yurinskii couplings available in the existing literature. Even
with independent data, we offer new Yurinskii couplings to Gaussian mixture
distributions (due to the presence of the common random factor loading $L$),
which could be further improved whenever the factors and residuals possess
symmetric (conditional) distributions. Furthermore, our results do not impose
any restrictions on the minimum eigenvalue of $\Sigma$, thereby allowing for
more general factor structures. These improvements are maintained in the
martingale, mixingale, and weakly dependent stationary data settings.

\section{Strong approximation for martingale empirical processes}%
\label{sec:emp_proc}

In this section, we demonstrate how our main results can be applied
to some more
substantive problems in statistics. Having until this point studied only
finite-dimensional (albeit potentially high-dimensional) random
vectors, we now
turn our attention to infinite-dimensional stochastic processes. Specifically,
we consider empirical processes of the form \[S(f) = \sum_{i=1}^{n} f(X_i),
\qquad f \in \cF,\] with $\cF$ a problem-specific class of real-valued
functions, where for each $f \in \cF$, the variables
$f(X_1), \ldots, f(X_n)$ form martingale differences with
respect to an appropriate filtration. We construct (conditionally) Gaussian
processes $T(f)$ for which upper bounds on the uniform coupling error
$\sup_{f \in \cF} |S(f) - T(f)|$ are precisely quantified. We control the
complexity of $\cF$ using metric entropy under Orlicz norms.

The novel strong approximation results which we present concern the entire
martingale empirical process $(S(f):f \in \cF)$, as opposed to just the scalar
supremum of the empirical process, $\sup_{f \in \cF} |S(f)|$. This distinction
has been carefully noted by \citet{chernozhukov2014gaussian}, who studied
Gaussian approximation of empirical process suprema in the independent data
setting and wrote (p.\ 1565): ``A related but different problem is that of
approximating \textit{whole} empirical processes by a sequence of Gaussian
processes in the sup-norm. This problem is more difficult than
[approximating the supremum of the empirical process].''
Indeed, the results we establish in
this section are for strong approximations of entire empirical processes by
sequences of Gaussian mixture processes in supremum norm, when the data
has a martingale difference structure (cf.\ Corollary
\ref{cor:sa_martingale}).
Our results can be further generalized to \emph{approximate} martingale
empirical processes (including mixingale empirical processes;
cf.\ Corollary \ref{cor:sa_mixingale}), but
to reduce notation and the technical burden
we do not consider this extension.

\subsection{Motivating example: kernel density estimation}
\label{sec:kde}

We begin with a brief study of a canonical example of an empirical process
which is non-Donsker (thus precluding the use of uniform central limit
theorems) due to the presence of a function class whose complexity increases
with the sample size: the kernel density estimator with i.i.d.\ scalar data.
We give an overview of our general strategy for
strong approximation of stochastic processes
via discretization, and show explicitly in Lemma~\ref{lem:kde_eigenvalue}
how it is crucial
that we do not impose lower bounds on the eigenvalues of the discretized
covariance matrix. Detailed calculations for this section are
relegated to
\iftoggle{aos}{the supplementary material
\citep{cattaneo2025yurinskiisupplement}}{Appendix~\ref{sec:proofs}}
for conciseness.

Let $X_1, \ldots, X_n$ be i.i.d.\ $\Unif[0,1]$, take
$K(x) = \frac{1}{\sqrt{2 \pi}} e^{-x^2/2}$ the Gaussian kernel and let
$h \in (0,1]$ be a bandwidth. Then, for $a \in (0,1/4]$ and
$x \in \cX = [a, 1-a]$ to avoid boundary issues, the kernel density estimator
of the true density function $g(x) = 1$ is
\begin{align*}
  \hat g(x)
  &=
  \frac{1}{n}
  \sum_{i=1}^{n}
  K_h( X_i - x),
  \qquad K_h(u) = \frac{1}{h} K\left( \frac{u}{h} \right).
\end{align*}
Consider establishing a strong approximation for the process
$(\hat g(x)-\E [ \hat g(x) ] : x\in\cX)$
which is, upon rescaling, non-Donsker whenever
the bandwidth decreases to zero in large samples.
To match notation with the upcoming
general result for empirical processes, set
$f_x(u) = \frac{1}{n} (K_h( u - x) - \E[K_h( X_i - x)])$
so $S(x) \vcentcolon= S(f_x) = \hat g(x)-\E [ \hat g(x) ]$.
The next step is standard: a
mesh separates the local oscillations of the processes from
the finite-dimensional coupling. %
For $\delta \in (0,1/2)$, set
$N = \left\lfloor 1 + \frac{1 - 2a}{\delta} \right\rfloor$
and $\cX_\delta = (a + (j-1)\delta : 1 \leq j \leq N)$.
Letting $T(x)$ be the approximating stochastic
process to be constructed, consider the following decomposition:
\begin{align*}
  \sup_{x \in \cX}
  \big|S(x) - T(x)\big|
  &\leq
  \iftoggle{aos}{\!\!}{}
  \sup_{|x-x'| \leq \delta}
  \iftoggle{aos}{\!\!}{}
  \big|S(x) - S(x') \big|
  + \max_{x \in \cX_\delta}
  \big|S(x) - T(x)\big|
  + \iftoggle{aos}{\!\!}{}
  \sup_{|x-x'| \leq \delta}
  \iftoggle{aos}{\!\!}{}
  \big|T(x) - T(x')\big|.
\end{align*}
Writing $S(\cX_\delta)$ for
$\big(S(x) : x \in \cX_\delta\big)\in \mathbb{R}^N$,
and noting that this is a sum of i.i.d.\ random vectors, we apply
Corollary~\ref{cor:sa_indep} as
$\max_{x \in \cX_\delta} |S(x) - T(x)|
= \| S(\cX_\delta) - T(\cX_\delta) \|_\infty$.
We thus obtain that, for each
$\eta > 0$, there exists a Gaussian vector
$T(\cX_\delta)$ with the same covariance matrix as $S(\cX_\delta)$
satisfying
\begin{align*}
  \P\left(
    \|S(\cX_\delta) - T(\cX_\delta)\|_\infty > \eta
  \right)
  &\leq
  31 \left(
    \frac{N \log 2 N}{\eta^3 n^2 h^2}
  \right)^{1/3}
\end{align*}
assuming that $1/h \geq \log 2 N$.
By the Vorob'ev--Berkes--Philipp theorem
\citep[Theorem~1.1.10]{dudley1999uniform},
$T(\cX_\delta)$ extends to a Gaussian process $T(x)$
defined for all $x \in \cX$ and with the same covariance structure
as $S(x)$.

Next, it is not difficult to show by chaining with the Bernstein--Orlicz and
sub-Gaussian norms respectively \citep[Section~2.2]{van1996weak} that if
$\log(N/h) \lesssim \log n$ and $n h \gtrsim \log n$,
\begin{align*}
  \sup_{|x-x'| \leq \delta}
  \big\|S(x) - S(x') \big\|_\infty
  &\lesssim_\P
  \delta
  \sqrt{\frac{\log n}{n h^3}}, \quad\text{and}\quad
  \sup_{|x-x'| \leq \delta}
  \big\|T(x) - T(x')\big\|_\infty
  \lesssim_\P
  \delta
  \sqrt{\frac{\log n}{n h^3}}.
\end{align*}
Finally, for any sequence $R_n\to\infty$
(Remark~\ref{rem:coupling_bounds_probability}),
the resulting bound on the coupling error is
\begin{align*}
  \sup_{x \in \cX}
  \big| S(x) - T(x) \big| \lesssim_\P
  \left( \frac{N \log 2N}{n^2 h^2} \right)^{1/3} R_n
  + \delta \sqrt{\frac{\log n}{n h^3}},
\end{align*}
where the mesh size $\delta$ is then optimized to obtain the tightest
possible strong approximation.
In particular, since $N \lesssim 1/\delta$, setting
$\delta \asymp n^{-1/8} h^{5/8} (\log n)^{-1/8}$ yields
\begin{align*}
  \sup_{x \in \cX}
  \big| S(x) - T(x) \big| \lesssim_\P
  \left( \frac{(\log n)^3}{n^5 h^7} \right)^{1/8} R_n
\end{align*}
which, after standardization by $\sqrt{n h}$,
vanishes whenever $R_n (\log n)^3 / (n h^3) \to 0$. This is a more
stringent assumption on the bandwidth
$h$ than $(\log n) / (n h) \to 0$ imposed by \citet{gine2004kernel}
and \citet{cattaneo2024strong} when employing a Hungarian
construction \citep{komlos1975approximation}, or $(\log n)^6 / (n h) \to 0$
imposed by \citet{chernozhukov2014gaussian} when studying in
particular the Kolmogorov--Smirnov distance between the scalar
suprema. The difference in side restrictions is a result of the
specific assumptions imposed and coupling approaches used; see
Section~\ref{sec:local_poly} for related discussion.

The discretization strategy outlined above is at the core of the
proof strategy
for our upcoming Proposition~\ref{pro:emp_proc}. Since we will consider
martingale empirical processes, our proof will rely on
Corollary~\ref{cor:sa_martingale}, which, unlike the martingale Yurinskii
coupling established by \citet{li2020uniform}, does not require a lower bound
on the minimum eigenvalue of $\Sigma$. Using the simple kernel density example
just discussed, we now demonstrate precisely the crucial importance
of removing
such eigenvalue conditions. The following Lemma~\ref{lem:kde_eigenvalue} shows
that the discretized covariance matrix $\Sigma = n h\Var[S(\cX_\delta)]$ has
exponentially small eigenvalues, which in turn will negatively affect the
strong approximation bound if the \citet{li2020uniform} coupling were to be
used instead of the results in this paper.

\begin{figure}[ht]
  \centering
  \begin{subfigure}{0.45\textwidth}
    \centering
    \includegraphics[width=\textwidth]{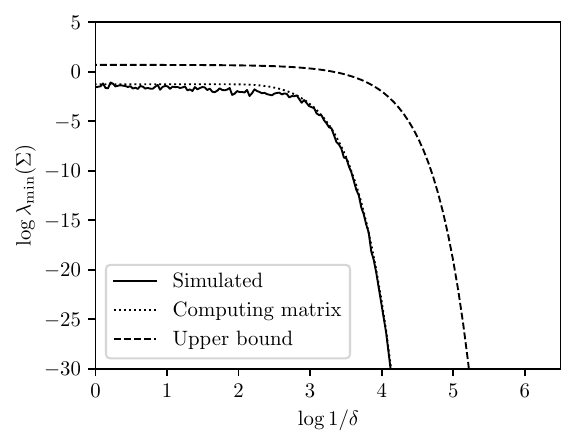}
    \caption{Bandwidth $h = 0.03$}
  \end{subfigure}
  \begin{subfigure}{0.45\textwidth}
    \centering
    \includegraphics[width=\textwidth]{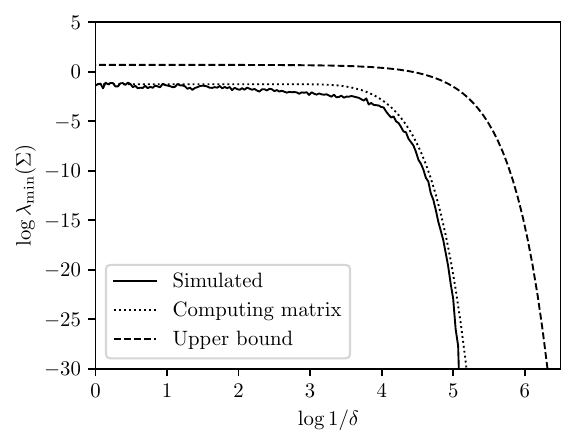}
    \caption{Bandwidth $h = 0.01$}
  \end{subfigure}
  \caption{
    Upper bounds on the minimum eigenvalue of the discretized covariance
    matrix in kernel density estimation,
    with $n=100$ and $a = 0.2$.
    Simulated: the kernel density estimator is simulated,
    resampling the data $100$ times
    to estimate its covariance.
    Computing matrix: the minimum eigenvalue of the limiting covariance
    matrix $\Sigma$ is computed explicitly.
    Upper bound: the bound derived in
    Lemma~\ref{lem:kde_eigenvalue}
    is shown.
  }
  \label{fig:min_eig}
\end{figure}

\begin{lemma}[Minimum eigenvalue of a
  kernel density estimator covariance matrix]%
  \label{lem:kde_eigenvalue}
  The minimum eigenvalue of
  $\Sigma=n h\Var[S(\cX_\delta)] \in \R^{N \times N}$
  satisfies the upper bound
  \begin{align*}
    \lambda_{\min}(\Sigma)
    &\leq
    2 e^{-h^2/\delta^2}
    + \frac{h}{\pi a \delta}
    e^{-a^2 / h^2}.
  \end{align*}
\end{lemma}
Figure~\ref{fig:min_eig} shows how the upper bound in Lemma
\ref{lem:kde_eigenvalue} captures the behavior of the simulated minimum
eigenvalue of $\Sigma$. In particular, the smallest eigenvalue decays
exponentially fast in the discretization level $\delta$ and the bandwidth $h$.
As seen in the calculations above, the coupling rate depends on $\delta / h$,
while the bias will generally depend on $h$, implying that both $\delta$ and
$h$ must converge to zero to ensure valid statistical inference. In general,
this will lead to $\Sigma$ possessing extremely small eigenvalues, rendering
strong approximation approaches such as that of \citet{li2020uniform}
ineffective in such scenarios.

The discussion in this section focuses on the strong approximation of the
centered process $\hat g(x)-\E [ \hat g(x) ]$. In practice, the goal is often
rather to approximate $\hat g(x)- g(x)$. The difference
between these is captured by the smoothing bias $\E [ \hat g(x) ] - g(x)$,
which is straightforward to control with
$\sup_{x \in \cX} \big| \E [ \hat g(x) ] - g(x) \big|
\lesssim \frac{h}{a} e^{-a^2 / (2 h^2)}$.
See Section \ref{sec:nonparametric} for further
discussion.

\subsection{General result for martingale empirical processes}

We now give our general result on a strong approximation for
martingale empirical processes, obtained by applying
the first result \eqref{eq:sa_martingale_order_2} in
Corollary~\ref{cor:sa_martingale} with $p=\infty$
to a discretization of the empirical process,
as in Section~\ref{sec:kde}.
We then control the increments in the stochastic processes
using chaining with Orlicz norms,
but note that other tools are available,
including generalized entropy with bracketing \citep{geer2000empirical}
and sequential symmetrization \citep{rakhlin2015sequential}.

A class of functions is said to be \emph{pointwise measurable}
if it contains a countable subclass which is dense under
the pointwise convergence topology.
For a finite class $\cF$, write
$\cF(x) = \big(f(x) : f \in \cF\big)$.
Define the set of Orlicz functions
\begin{align*}
  \Psi
  &=
  \Big\{
    \psi: [0, \infty) \to [0, \infty)
    \text{ convex increasing, }
    \psi(0) = 0,\
    \limsup_{x,y \to \infty} \tfrac{\psi(x) \psi(y)}{\psi(C x y)} < \infty
    \text{ for } C > 0
  \Big\}
\end{align*}
and, for real-valued $Y$, the Orlicz norm
$\vvvert Y \vvvert_\psi
= \inf
\left\{ C > 0:
  \E \left[ \psi(|Y|/C) \leq 1 \right]
\right\}$
as in \citet[Section~2.2]{van1996weak}.

\begin{proposition}[Strong approximation for martingale empirical processes]%
  \label{pro:emp_proc}

  Let $X_i$ be random variables for $1 \leq i \leq n$ taking values in a
  measurable space $\cX$, and $\cF$ be a pointwise measurable class of
  functions from $\cX$ to $\R$. Let $\cH_0, \ldots, \cH_n$ be a
  filtration such
  that each $X_i$ is $\cH_i$-measurable, with $\cH_0$ the trivial
  $\sigma$-algebra, and suppose that $\E[f(X_i) \mid \cH_{i-1}] = 0$ for all
  $f \in \cF$. Define $S(f) = \sum_{i=1}^n f(X_i)$ for $f\in\cF$ and let
  $\Sigma: \cF \times \cF \to \R$ be an almost surely positive semi-definite
  $\cH_0$-measurable random function. Suppose that for a non-random
  metric $d$ on $\cF$, constant $L$ and $\psi \in \Psi$,
  \begin{align}%
    \label{eq:emp_proc_var}
    \Sigma(f,f) - 2\Sigma(f,f') + \Sigma(f',f')
    + \bigvvvert S(f) - S(f') \bigvvvert_\psi^2
    &\leq L^2 d(f,f')^2 \quad \text{a.s.}
  \end{align}
  Then for each $\eta > 0$ there is a
  process $T(f)$ indexed by $f\in\cF$
  which, conditional on $\cH_0$, is zero-mean and Gaussian,
  satisfying
  $\E\big[ T(f) T(f') \mid \cH_0 \big] = \Sigma(f,f')$
  for all $f, f' \in \cF$, and for all $t > 0$ has
  \begin{align*}
    &\P\left(
      \sup_{f \in \cF}
      \big| S(f) - T(f) \big|
      \geq C_\psi(t + \eta)
    \right)
    \leq
    C_\psi
    \inf_{\delta > 0}
    \inf_{\cF_\delta}
    \Bigg\{
      \frac{\beta_\delta^{1/3} (\log 2 |\cF_\delta|)^{1/3}}{\eta } \\
      &\qquad\quad+
      \left(\frac{\sqrt{\log 2 |\cF_\delta|}
      \sqrt{\E\left[\|\Omega_\delta\|_2\right]}}{\eta }\right)^{2/3}
      + \psi\left(\frac{t}{L J_\psi(\delta)}\right)^{-1}
      + \exp\left(\frac{-t^2}{L^2 J_2(\delta)^2}\right)
    \Bigg\},
  \end{align*}
  where $\cF_\delta$ is any finite $\delta$-cover of $(\cF,d)$
  and $C_\psi$ is a constant depending only on $\psi$, with
  \begin{align*}
    \beta_\delta
    &= \sum_{i=1}^n
    \E\left[ \|\cF_\delta(X_i)\|^2_2\|\cF_\delta(X_i)\|_\infty
      + \|V_i(\cF_\delta)^{1/2}Z_i\|^2_2
    \|V_i(\cF_\delta)^{1/2}Z_i\|_\infty \right],
  \end{align*}
  \vspace*{-6mm}
  \begin{align*}
    V_i(\cF_\delta)
    &=
    \E\big[\cF_\delta(X_i) \cF_\delta(X_i)^\T \mid \cH_{i-1} \big],
    &
    \Omega_\delta
    &=
    \sum_{i=1}^n V_i(\cF_\delta) - \Sigma(\cF_\delta), \\
    J_\psi(\delta)
    &=
    \int_0^\delta \psi^{-1}\big( N_\varepsilon \big)
    \diff{\varepsilon}
    + \delta \psi^{-1} \big( N_\delta^2 \big),
    &
    J_2(\delta)
    &= \int_0^\delta \sqrt{\log N_\varepsilon}
    \diff{\varepsilon},
  \end{align*}
  where $N_\delta = N(\delta, \cF, d)$
  is the $\delta$-covering number of $(\cF, d)$
  and $Z_i$ are i.i.d.\ $\cN\big(0, I_{|\cF_\delta|}\big)$
  independent of $\cH_n$.
  If $\cF_\delta$ is a minimal $\delta$-cover
  of $(\cF, d)$, then $|\cF_\delta| = N_\delta$.
\end{proposition}

Proposition~\ref{pro:emp_proc}
is given in a rather general form to accommodate a range of different
settings and applications.
In particular, consider the following well-known Orlicz functions.
\begin{description}

  \item[Polynomial:]
    $\psi(x) = x^a$ for $a \geq 2$
    has $\vvvert X \vvvert_2 \leq \vvvert X \vvvert_\psi$ and
    $\sqrt{\log x} \leq \sqrt{a} \psi^{-1}(x)$.

  \item[Exponential:]
    $\psi(x) = \exp(x^a) - 1$ for $a \in [1,2]$
    has $\vvvert X \vvvert_2 \leq 2\vvvert X \vvvert_\psi$ and
    $\sqrt{\log x} \leq \psi^{-1}(x)$.

  \item[Bernstein:]
    $\psi(x) = \exp\hspace*{-0.8mm}
    \Big(\hspace*{-0.6mm}
      \Big(\frac{\sqrt{1+2ax}-1}{a}\Big)^{\hspace*{-0.6mm}2}
    \Big)-1$
    for $a > 0$ has
    $\vvvert X \vvvert_2 \leq (1+a)\vvvert X \vvvert_\psi$ and
    $\sqrt{\log x}~\leq~\psi^{-1}(x)$.

\end{description}
For these Orlicz functions and when $\Sigma(f, f') = \Cov[S(f), S(f')]$ is
non-random, the terms involving $\Sigma$ in \eqref{eq:emp_proc_var} can be
controlled by the Orlicz $\psi$-norm term; similarly, $J_2$ is bounded by
$J_\psi$. Further, $C_\psi$ can be replaced by a universal constant $C$ which
does not depend on the parameter $a$. See Section~2.2 in \citet{van1996weak}
for details. If the conditional third moments of $\cF_\delta(X_i)$
given $\cH_{i-1}$ are
all zero (if $f$ and $X_i$ are appropriately symmetric, for example), then the
second inequality in Corollary~\ref{cor:sa_martingale} can be
applied to obtain
a tighter coupling inequality; the details of this are omitted for
brevity, and
the proof would proceed in exactly the same manner.

In general, however, Proposition~\ref{pro:emp_proc} allows for a random
covariance function, yielding a coupling to a stochastic process that is
Gaussian only conditionally. Such a process can equivalently be
formally viewed as a
mixture of Gaussian processes, writing $T=\Sigma^{1/2} Z$ with an operator
square root and where $Z$ is a Gaussian white noise on $\cF$ independent of
$\cH_0$. This extension is in contrast with much of the existing strong
approximation and empirical process literature, which tends to focus on
couplings and weak convergence results with marginally Gaussian processes.

A similar approach was taken by \citet{berthet2006revisiting}, who used a
Gaussian coupling due to
\citet{zaitsev1987estimates,zaitsev1987gaussian} along
with a discretization method to obtain strong approximations for empirical
processes with independent data. They handled fluctuations in the stochastic
processes with uniform $L^2$ covering numbers and bracketing numbers where we
opt instead for chaining with Orlicz norms. Our version using the (martingale)
Yurinskii coupling can improve upon theirs in approximation rate even for
independent data under certain circumstances, as follows. Suppose the setup of
Proposition~1 in \citet{berthet2006revisiting}; that is, $X_1,
\ldots, X_n$ are
i.i.d.\ and $\sup_{\cF} \|f\|_\infty \leq M$, with the VC-type assumption
$\sup_\Q N(\varepsilon, \cF, d_\Q) \leq c_0 \varepsilon^{-\nu_0}$ where
$d_\Q(f,f')^2 = \E_\Q\big[(f-f')^2\big]$ for a measure $\Q$ on $\cX$ and
$M, c_0, \nu_0$ are constants. Then, using uniform $L^2$ covering numbers
rather than Orlicz norm chaining in our Proposition~\ref{pro:emp_proc}
gives the following.
Firstly as $X_i$ are i.i.d.\ we take $\Sigma(f, f') = \Cov[S(f), S(f')]$ so
$\Omega_\delta = 0$. Let $\cF_\delta$ be a minimal $\delta$-cover of
$(\cF, d_\P)$ with cardinality $N_\delta \lesssim \delta^{-\nu_0}$ where
$\delta \to 0$. It is not difficult to show that
$\beta_\delta \lesssim n \delta^{-\nu_0} \sqrt{\log(1/\delta)}$.
Theorem~2.2.8 and Theorem~2.14.1 in \citet{van1996weak} give
\begin{align*}
  \E\left[
    \sup_{d_\P(f,f') \leq \delta}
    \Big(
      |S(f) - S(f')|
      + |T(f) - T(f')|
    \Big)
  \right]
  &\lesssim
  \sup_\Q
  \int_0^\delta
  \sqrt{n \log N(\varepsilon, \cF, d_\Q)}
  \diff{\varepsilon}
  \iftoggle{aos}{\\ &\lesssim}{\lesssim}
  \delta \sqrt{n\log(1/\delta)},
\end{align*}
where we used the VC-type property to bound the entropy integral.
So by our Proposition~\ref{pro:emp_proc},
for any sequence $R_n \to \infty$
(see Remark~\ref{rem:coupling_bounds_probability}),
\begin{align*}
  \sup_{f \in \cF}
  \big| S(f) - T(f) \big|
  &\lesssim_\P
  n^{1/3} \delta^{-\nu_0/3}
  \sqrt{\log(1/\delta)} R_n
  + \delta \sqrt{n\log(1/\delta)}
  \lesssim_\P
  n^{\frac{2+\nu_0}{6+2\nu_0}}
  \sqrt{\log n} R_n,
\end{align*}
where we minimized over $\delta$ in the last step.
\citet[Proposition~1]{berthet2006revisiting} achieved
\begin{align*}
  \sup_{f \in \cF}
  \big| S(f) - T(f) \big|
  &\lesssim_\P
  n^{\frac{5\nu_0}{4+10\nu_0}}
  (\log n)^{\frac{4+5\nu_0}{4+10\nu_0}},
\end{align*}
showing that our approach achieves a better approximation rate whenever
$\nu_0 > 4/3$. In particular, our method is superior in richer
function classes
with larger VC-type dimension. For example, if $\cF$ is smoothly parametrized
by $\theta \in \Theta \subseteq \R^d$ where $\Theta$ contains an
open set, then
$\nu_0 > 4/3$ corresponds to $d \geq 2$ and our rate is better as soon as the
parameter space is more than one-dimensional. The difference in approximation
rate is due to Zaitsev's coupling having better dependence on the sample size
but worse dependence on the dimension. In particular, Zaitsev's coupling is
stated only in $\ell_2$-norm and hence
\citet[Equation~5.3]{berthet2006revisiting} are compelled to use
the inequality
$\|\cdot\|_\infty \leq \|\cdot\|_2$ in the coupling step, a bound which is
loose when the dimension of the vectors (here on the order of
$\delta^{-\nu_0}$) is even moderately large. We use the fact that our version
of Yurinskii's coupling applies directly to the supremum norm, giving sharper
dependence on the dimension.

In Section~\ref{sec:local_poly} we apply Proposition~\ref{pro:emp_proc} to
obtain strong approximations for local polynomial estimators in the
nonparametric regression setting. In contrast with the series
estimators of the
upcoming Section~\ref{sec:series}, local polynomial estimators are
not linearly
separable and hence cannot be analyzed directly using the finite-dimensional
Corollary~\ref{cor:sa_martingale}.

\section{Applications to nonparametric regression}
\label{sec:nonparametric}

We illustrate the applicability of our previous strong approximation results
with two substantial and classical examples in nonparametric regression
estimation. Firstly, we present an analysis of partitioning-based series
estimators, in which we can apply the finite-dimensional result of
Corollary~\ref{cor:sa_martingale} directly
due to an intrinsic linear separability property. Secondly, we consider local
polynomial estimators, this time using the stochastic process
formulation in Proposition~\ref{pro:emp_proc} due to
the presence of a non-linearly separable martingale empirical process.

\subsection{Partitioning-based series estimators}
\label{sec:series}

Partitioning-based least squares methods are essential tools for
estimation and
inference in nonparametric regression, encompassing splines, piecewise
polynomials, compactly supported wavelets and decision trees as special cases.
See \citet{cattaneo2020large} for further details and references throughout
this section. We illustrate the usefulness of
Corollary~\ref{cor:sa_martingale}
by deriving a Gaussian strong approximation for partitioning series estimators
based on multivariate martingale data. Proposition~\ref{pro:series} shows how
we achieve the best known rate of strong approximation for independent data by
imposing an additional mild $\alpha$-mixing condition to control the time
series dependence of the regressors.

Consider the nonparametric regression setup with martingale difference
residuals defined by $Y_i = \mu(W_i) + \varepsilon_i$ for $ 1 \leq i \leq n$
where the regressors $W_i$ have compact connected support $\cW
\subseteq \R^m$,
$\cH_i$ is the $\sigma$-algebra generated by
$(W_1, \ldots, W_{i+1}, \varepsilon_1, \ldots, \varepsilon_i)$,
$\E[\varepsilon_i \mid \cH_{i-1}] = 0$ and $\mu: \cW \to \R$ is the estimand.
Let $p(w)$ be a $k$-dimensional vector of bounded basis functions on $\cW$
which are locally supported on a quasi-uniform partition
\citep[Assumption~2]{cattaneo2020large}. Under minimal regularity conditions,
the least-squares partitioning-based series estimator is
$\hat\mu(w) = p(w)^{\T} \hat H^{-1} \sum_{i=1}^n p(W_i) Y_i$
with $\hat H = \sum_{i=1}^n p(W_i) p(W_i)^\T$.
The approximation power of the estimator $\hat\mu(w)$ derives from letting
$k\to\infty$ as $n\to\infty$. The assumptions made on $p(w)$ are
mild enough to
accommodate splines, wavelets, piecewise polynomials, and certain types of
decision trees. For such a tree, $p(w)$ is comprised of indicator functions
over $k$ axis-aligned rectangles forming a partition of $\cW$ (a Haar basis),
provided that the partitions are constructed using independent data
(e.g., with sample splitting).

Our goal is to approximate the law of the stochastic process
$(\hat\mu(w)-\mu(w):w\in\cW)$, which upon rescaling is typically not
asymptotically tight as $k \to \infty$ and thus does not converge weakly.
Nevertheless, exploiting the intrinsic linearity of the estimator
$\hat\mu(w)$,
we can apply Corollary~\ref{cor:sa_martingale} directly to
construct a Gaussian
strong approximation. Specifically, we write
\begin{align*}
  \hat\mu(w) - \mu(w)
  &= p(w)^\T H^{-1} S
  + p(w)^\T \big(\hat H^{-1} - H^{-1}\big) S
  + \Bias(w),
\end{align*}
where $H= \sum_{i=1}^n \E\left[p(W_i) p(W_i)^\T\right]$
is the expected outer product matrix, $S = \sum_{i=1}^n p(W_i) \varepsilon_i$
is the score vector, and
$\Bias(w) = p(w)^{\T} \hat H^{-1}\sum_{i=1}^n p(W_i) \mu(W_i) - \mu(w)$.
Imposing some mild time series restrictions and assuming stationarity for
simplicity, it is not difficult to show (see
  \iftoggle{aos}{the supplementary material
\citep{cattaneo2025yurinskiisupplement}}{Appendix~\ref{sec:proofs}})
that $\|\hat H - H\|_1 \lesssim_\P \sqrt{n k}$ and
$\sup_{w\in\cW} |\Bias(w)| \lesssim_\P k^{-\gamma}$
for some $\gamma>0$, depending on the specific structure of the basis
functions, the dimension $m$ of the regressors, and the smoothness of the
regression function $\mu$. Thus, it remains to study the $k$-dimensional
zero-mean martingale $S$ by applying Corollary~\ref{cor:sa_martingale} with
$X_i=p(W_i) \varepsilon_i$. Controlling the convergence of the quadratic
variation term $\E[\|\Omega\|_2]$ also requires some time series dependence
assumptions; we impose an $\alpha$-mixing condition on $(W_1,
\ldots, W_n)$ for
illustration \citep{bradley2005basic}.

\begin{proposition}[Strong approximation for partitioning series estimators]%
  \label{pro:series}
  Consider the nonparametric regression setup described above
  and further assume the following:
  \begin{enumerate}[label=(\roman*)]

    \item
      $(W_i, \varepsilon_i)_{1 \leq i \leq n}$
      is strictly stationary.

    \item
      $W_1, \ldots, W_n$ is $\alpha$-mixing with mixing coefficients
      satisfying $\sum_{j=1}^\infty \alpha(j) < \infty$.

    \item
      $W_i$ has a Lebesgue density on $\cW$
      which is bounded above and away from zero.

    \item
      $\E\big[|\varepsilon_i|^3 \big] < \infty$
      and
      $\E\big[\varepsilon_i^2 \mid \cH_{i-1}\big]=\sigma^2(W_i)$
      is bounded away from zero.

    \item
      $p(w)$ forms a basis with $k$ features satisfying
      Assumptions~2 and~3 in \citet{cattaneo2020large}.

  \end{enumerate}
  Then, for any sequence $R_n \to \infty$,
  there is a zero-mean Gaussian process
  $G(w)$ indexed on $\cW$
  with $\Var[G(w)] \asymp\frac{k}{n}$
  satisfying
  $\Cov[G(w), G(w')]
  = \Cov[p(w)^\T H^{-1} S,\, p(w')^\T H^{-1} S]$
  and
  \begin{align*}
    \sup_{w \in \cW}
    \left| \hat\mu(w) - \mu(w) - G(w) \right|
    &\lesssim_\P
    \sqrt{\frac{k}{n}}
    \left( \frac{k^3 (\log k)^3}{n} \right)^{1/6} R_n
    + \sup_{w \in \cW} |\Bias(w)|
  \end{align*}
  assuming the number of basis functions satisfies $k^3 / n \to 0$.
  If further $\E \left[ \varepsilon_i^3 \mid \cH_{i-1} \right] = 0$ then
  \begin{align*}
    \sup_{w \in \cW}
    \left| \hat\mu(w) - \mu(w) - G(w) \right|
    &\lesssim_\P
    \sqrt{\frac{k}{n}}
    \left( \frac{k^3 (\log k)^2}{n} \right)^{1/4} R_n
    + \sup_{w \in \cW} |\Bias(w)|.
  \end{align*}
\end{proposition}

The core of the proof of Proposition~\ref{pro:series} involves applying
Corollary~\ref{cor:sa_martingale} with $S = \sum_{i=1}^n p(W_i) \varepsilon_i$
and $p=\infty$ to construct $T \sim \cN\big(0, \Var[S]\big)$ such that
$\|S - T \|_\infty$ is small, and then setting $G(w) = p(w)^\T H^{-1} T$. So
long as the bias can be appropriately controlled, this result allows for
uniform inference procedures such as uniform confidence bands or shape
specification testing. The condition $k^3 / n \to 0$ is the same (up to logs)
as that imposed by \citet{cattaneo2020large} for i.i.d. data, which gives the
best known strong approximation rate for this problem. Thus,
Proposition~\ref{pro:series} gives the same best approximation rate, without
requiring any extra restrictions, for $\alpha$-mixing time series data.

Our results improve substantially on \citet[Theorem~1]{li2020uniform}: using
the notation of our Corollary~\ref{cor:sa_martingale}, and with any sequence
$R_n \to \infty$, a valid (see Remark~\ref{rem:coupling_bounds_probability})
version of their martingale Yurinskii coupling is
\begin{align*}
  \|S-T\|_2
  \lesssim_\P
  d^{1/2} r^{1/2}_n
  + (B_n d)^{1/3} R_n,
\end{align*}
where $B_n = \sum_{i=1}^n \E[\|X_i\|_2^3]$ and $r_n$ is a term controlling the
convergence of the quadratic variation, playing a similar role to our
term $\E[\|\Omega\|_2]$. Under the assumptions of our
Proposition~\ref{pro:series}, applying this
result with $S = \sum_{i=1}^n p(W_i) \varepsilon_i$ yields a rate no better
than $\|S-T\|_2 \lesssim_\P (n k)^{1/3} R_n$. As such, they attain a rate of
strong approximation no faster than
\begin{align*}
  \sup_{w \in \cW}
  \left| \hat\mu(w) - \mu(w) - G(w) \right|
  &\lesssim_\P
  \sqrt{\frac{k}{n}}
  \left( \frac{k^5}{n} \right)^{1/6} R_n
  + \sup_{w \in \cW} |\Bias(w)|.
\end{align*}
Hence, for this approach to yield a valid strong approximation, the number of
basis functions must satisfy $k^5/n \to 0$, a more restrictive assumption than
our $k^3 / n \to 0$ (up to logs). This difference is due to
\citet{li2020uniform} using the $\ell_2$-norm version of Yurinskii's coupling
rather than the more recently established $\ell_\infty$-norm version. Further,
our approach allows for an improved rate of distributional approximation
whenever the residuals have zero conditional third moment.

To illustrate the statistical applicability of
Proposition~\ref{pro:series}, consider constructing a feasible uniform
confidence band for the regression function $\mu$, using standardization and
Studentization for statistical power improvements. We assume throughout that
the bias is negligible. Proposition~\ref{pro:series} and
anti-concentration for
Gaussian suprema \citep[Corollary~2.1]{chernozhukov2014anti} yield
a distributional approximation for the supremum statistic whenever
$k^3(\log n)^6 / n \to 0$, giving
\begin{align*}
  \sup_{t \in \R}
  \left|
  \P\left(
    \sup_{w \in \cW}
    \left|
    \frac{\hat\mu(w)-\mu(w)}{\sqrt{\rho(w,w)}}
    \right| \leq t
  \right)
  -
  \P\left(
    \sup_{w \in \cW}
    \left|
    \frac{G(w)}{\sqrt{\rho(w,w)}}
    \right| \leq t
  \right)
  \right|
  &\to 0,
\end{align*}
where $\rho(w,w') = \E[G(w)G(w')]$. Furthermore, using a Gaussian--Gaussian
comparison result \citep[Lemma~3.1]{chernozhukov2013gaussian} and
anti-concentration again, it is not difficult to show (see the proof of
Proposition~\ref{pro:series}) that with $\bW = (W_1, \ldots, W_n)$ and
$\bY = (Y_1, \ldots, Y_n)$,
\begin{align*}
  \sup_{t \in \R}
  \left|
  \P\left(
    \sup_{w \in \cW}
    \left|
    \frac{\hat\mu(w)-\mu(w)}{\sqrt{\hat\rho(w,w)}}
    \right| \leq t
  \right)
  - \P\left(
    \sup_{w \in \cW}
    \left|
    \frac{\hat G(w)}{\sqrt{\hat\rho(w,w)}}
    \right| \leq t \biggm| \bW, \bY
  \right)
  \right|
  &\to_\P 0,
\end{align*}
where $\hat G(w)$ is a zero-mean Gaussian process
conditional on $\bW$ and $\bY$ with conditional covariance function
$\hat\rho(w,w')
=\E\big[\hat G(w) \hat G(w') \mid \bW, \bY \big]
= p(w)^\T \hat H^{-1} \widehat{\Var}[S] \hat H^{-1}p(w')$
for some estimator $\widehat{\Var}[S]$ satisfying
$\frac{k (\log n)^2}{n}
\big\|\widehat{\Var}[S]-\Var[S]\big\|_2 \to_\P 0$.
For example, one could use the plug-in estimator
$\widehat{\Var}[S]=\sum_{i=1}^n p(W_i) p(W_i)^\T \hat{\sigma}^2(W_i)$
where $\hat{\sigma}^2(w)$ satisfies
$(\log n)^2 \sup_{w \in \cW}
|\hat{\sigma}^2(w)-\sigma^2(w)| \to_\P 0$.
This leads to the following feasible and asymptotically valid
$100(1-\tau)\%$
uniform confidence band for partitioning-based series estimators
based on martingale data.

\begin{proposition}[Feasible uniform confidence bands for partitioning
  series estimators]%
  \label{pro:series_feasible}
  Assume the setup as described above. Then
  \begin{align*}
    \P\Big(
      \mu(w) \in
      \Big[
        \hat\mu(w) \pm \hat q(\tau)
        \sqrt{\hat\rho(w,w)}
      \Big]
      \ \text{for all }
    w \in \cW \Big)
    \to 1-\tau,
  \end{align*}
  where
  \vspace*{-2mm}
  \begin{align*}
    \hat{q}(\tau)
    &=
    \inf
    \left\{
      t \in \R:
      \P\left(
        \sup_{w \in \cW}
        \left|
        \frac{\hat G(w)}{\sqrt{\hat\rho(w,w)}}
        \right|
        \leq t
        \Bigm| \bW, \bY
      \right)
      \geq \tau
    \right\}
  \end{align*}
  is the conditional quantile of the supremum of the Studentized Gaussian
  process. This can be estimated by resampling the conditional law of
  $\hat G(w) \mid \bW, \bY$ with a discretization of $w \in \cW$.
\end{proposition}

\subsection{Local polynomial estimators}
\label{sec:local_poly}

As a second example application we consider nonparametric
regression estimation
with martingale data employing local polynomial methods
\citep{Fan-Gijbels_1996_Book}. In contrast with the partitioning-based series
methods of Section~\ref{sec:series}, local polynomials induce stochastic
processes which are not linearly separable, allowing us to showcase the
empirical process result given in Proposition \ref{pro:emp_proc}.

As before, suppose that
$Y_i = \mu(W_i) + \varepsilon_i$
for $ 1 \leq i \leq n$
where $W_i$ has compact connected support $\cW \subseteq \R^m$,
$\cH_i$ is the $\sigma$-algebra generated by
$(W_1, \ldots, W_{i+1}, \varepsilon_1, \ldots, \varepsilon_i)$,
$\E[\varepsilon_i \mid \cH_{i-1}] = 0$,
and $\mu: \cW \to \R$ is the estimand. Let $K$ be a kernel function on $\R^m$
and $K_h(w) = h^{-m} K(w/h)$
for some bandwidth $h > 0$.
Take $\gamma \geq 0$ and let
$k = (m+\gamma)!/(m!\gamma!)$ be the number of monomials up to order $\gamma$.
Using multi-index notation,
let $p(w)$ be the $k$-dimensional vector
collecting the monomials $w^{\kappa}/\kappa!$
for $0 \leq |\kappa| \leq \gamma$,
and set $p_h(w) = p(w/h)$.
The local polynomial regression estimator of $\mu(w)$ is,
with $e_1 = (1, 0, \ldots, 0)^\T \in \R^k$,
\begin{align*}
  \hat{\mu}(w)
  &=
  e_1^\T\hat{\beta}(w)
  &\text{where}      &
  &\hat{\beta}(w)
  &=
  \argmin_{\beta \in \R^{k}}
  \sum_{i=1}^n
  \left(Y_i - p_h(W_i-w)^\T \beta \right)^2
  K_h(W_i-w).
\end{align*}

Our goal is again to approximate the distribution of the entire stochastic
process, $(\hat{\mu}(w)-\mu(w):w\in\cW)$, which upon rescaling is non-Donsker
if $h \to 0$, and decomposes as follows:
\begin{align*}
  \hat{\mu}(w)-\mu(w)
  &= e_1^\T H(w)^{-1} S(w)
  + e_1^\T \big(\hat H(w)^{-1} - H(w)^{-1}\big) S(w)
  + \Bias(w)
\end{align*}
where
$\hat H(w) = \sum_{i=1}^n K_h(W_i-w) p_h(W_i-w) p_h(W_i-w)^\T$,
$H(w) = \E \big[ \hat H(w) \big]$,
$S(w)= \sum_{i=1}^n K_h(W_i-w) p_h(W_i-w) \varepsilon_i$
and
$\Bias(w) = e_1^\T \hat H(w)^{-1}
\sum_{i=1}^n K_h(W_i-w) p_h(W_i-w) \mu(W_i) - \mu(w)$.
A key distinctive feature of local polynomial regression is that both
$\hat H(w)$ and $S(w)$ are functions of the evaluation point $w\in\cW$;
contrast this with the partitioning-based series estimator discussed in
Section~\ref{sec:series}, for which neither $\hat H$ nor $S$ depend on $w$.
Therefore we use Proposition \ref{pro:emp_proc} to obtain a Gaussian strong
approximation for the martingale empirical process directly.

Under some mild regularity conditions, including stationarity for simplicity
and an $\alpha$-mixing assumption on the time-dependence of the data, we first
show
$\sup_{w\in\cW} \|\hat H(w)-H(w)\|_2
\lesssim_\P \sqrt{n h^{-2m}\log n}$.
Further,
$\sup_{w\in\cW} |\Bias(w)|
\lesssim_\P h^\gamma$
provided that the regression function is sufficiently smooth.
Thus it remains to analyze the martingale empirical process
$\big(e_1^\T H(w)^{-1} S(w) : w\in\cW\big)$
via Proposition \ref{pro:emp_proc} by setting
\begin{align*}
  \cF = \left\{
    (W_i, \varepsilon_i) \mapsto
    e_1^\T H(w)^{-1}
    K_h(W_i-w) p_h(W_i-w) \varepsilon_i
    : w \in \cW
  \right\}.
\end{align*}
With this approach, we obtain the following result.

\begin{proposition}[Strong approximation for local polynomial estimators]%
  \label{pro:local_poly}

  Under the nonparametric regression setup described above,
  assume further that
  \begin{enumerate}[label=(\roman*)]

    \item
      $(W_i, \varepsilon_i)_{1 \leq i \leq n}$
      is strictly stationary.

    \item
      $(W_i, \varepsilon_i)_{1 \leq i \leq n}$
      is $\alpha$-mixing with mixing coefficients
      $\alpha(j) \leq e^{-2 j / C_\alpha}$
      for some $C_\alpha > 0$.

    \item
      $W_i$ has a Lebesgue density on $\cW$
      which is bounded above and away from zero.

    \item
      $\E\big[e^{|\varepsilon_i|/C_\varepsilon}\big] < \infty$
      for $C_\varepsilon > 0$ and
      $\E\left[\varepsilon^2_i \mid \cH_{i-1}\right]=\sigma^2(W_i)$
      is bounded away from zero.

    \item
      $K$ is a non-negative Lipschitz
      compactly supported kernel with
      $\int K(w) \diff{w} < \infty$.

  \end{enumerate}
  Then for any $R_n \to \infty$,
  there is a zero-mean Gaussian process
  $T(w)$ on $\cW$
  with $\Var[T(w)] \asymp\frac{1}{n h^m}$
  satisfying
  $\Cov[T(w), T(w')]
  = \Cov[e_1^\T H(w)^{-1} S(w),\, e_1^\T H(w')^{-1} S(w')]$
  and
  \begin{align*}
    \sup_{w \in \cW}
    \left|\hat \mu(w) - \mu(w) - T(w) \right|
    &\lesssim_\P
    \frac{R_n}{\sqrt{n h^m}}
    \left(
      \frac{(\log n)^{m+4}}{n h^{3m}}
    \right)^{\frac{1}{2m+6}}
    + \sup_{w \in \cW} |\Bias(w)|,
  \end{align*}
  provided that the bandwidth sequence satisfies
  $n h^{3m} \to \infty$.
\end{proposition}

If the residuals further satisfy
$\E \left[ \varepsilon_i^3 \mid \cH_{i-1} \right] = 0$, then
a third-order Yurinskii coupling delivers an improved rate of strong
approximation for Proposition~\ref{pro:local_poly}; this is omitted here for
brevity. For completeness, the proof of Proposition~\ref{pro:local_poly}
verifies that if the regression function $\mu(w)$ is $\gamma$ times
continuously differentiable on $\cW$ then
$\sup_w |\Bias(w)| \lesssim_\P h^\gamma$. Further, the assumption that $p(w)$
is a vector of monomials is unnecessary in general; any collection of bounded
linearly independent functions which exhibit appropriate approximation power
will suffice \citep{eggermont2009maximum}. As such, we can encompass local
splines and wavelets, as well as polynomials, and also choose
whether or not to
include interactions between the regressor variables. The bandwidth
restriction
of $n h^{3m} \to \infty$ is analogous to that imposed in
Proposition~\ref{pro:series} for partitioning-based series estimators, and as
far as we know, has not been improved upon for non-i.i.d.\ data.

Applying an anti-concentration result for Gaussian process suprema, such as
Corollary~2.1 in \citet{chernozhukov2014anti}, allows one to write a
Kolmogorov--Smirnov bound comparing the law of
$\sup_{w \in \cW}|\hat\mu(w) - \mu(w)|$ to that of $\sup_{w \in \cW}|T(w)|$.
With an appropriate covariance estimator, we can further replace $T(w)$ by a
feasible version $\hat T(w)$ or its Studentized counterpart, enabling
procedures for uniform inference analogous to the confidence bands constructed
in Section~\ref{sec:series}. We omit the details of this to conserve space but
note that our assumptions on $W_i$ and $\varepsilon_i$ ensure that
Studentization is possible even when the discretized covariance matrix has
small eigenvalues (Section~\ref{sec:kde}), as we normalize only by
the diagonal
entries.

In this setting of kernel-based local empirical
processes, it is essential that our initial strong approximation result
(Corollary~\ref{cor:sa_martingale}) does not impose a lower bound on the
eigenvalues of the variance matrix $\Sigma$. This effect was demonstrated by
Lemma \ref{lem:kde_eigenvalue} and its surrounding discussion in
Section~\ref{sec:kde}, and as such, the result of \citet{li2020uniform} is
unsuited for this application due to its strong minimum eigenvalue assumption.
Finally, for the special case of i.i.d.\ data,
\citet[Remark~3.1]{chernozhukov2014gaussian} achieve better rates for
approximating the scalar supremum of the $t$-process in
Kolmogorov--Smirnov distance by bypassing the step where we first approximate
the entire stochastic process (see Section~\ref{sec:emp_proc} for a
discussion), while \citet{cattaneo2024strong} obtain better strong
approximations for the entire stochastic process under additional
assumptions via a generalization of the celebrated Hungarian
construction \citep{komlos1975approximation,rio1994local}.

\section{Conclusion}
\label{sec:conclusion}

In this paper we introduced as our main result a new version of Yurinskii's
coupling which strictly generalizes all previously known forms of the result.
Our formulation gave a Gaussian mixture coupling for approximate martingale
vectors in $\ell_p$-norm where $1 \leq p \leq \infty$, with no restrictions on
the minimum eigenvalues of the associated covariance matrices. We further
showed how to obtain an improved approximation whenever third moments of the
data are negligible. We demonstrated the applicability of this main result by
first deriving a user-friendly version, and then specializing it to
mixingales,
martingales, and independent data, illustrating the benefits with a collection
of simple factor models. We then considered the problem of
constructing uniform
strong approximations for martingale empirical processes,
demonstrating how our
new Yurinskii coupling can be employed in a stochastic process setting. As
substantive illustrative applications of our theory to some
well established %
problems in statistical methodology, we showed how to use our coupling results
for both vector-valued and empirical process-valued martingales in developing
uniform inference procedures for partitioning-based series
estimators and local
polynomial models in nonparametric regression.
At each stage we addressed issues of feasibility, compared our work with the
existing literature, and provided implementable statistical inference
procedures.
 
\section*{Acknowledgments}
We thank the Editor, Associate Editor, and several reviewers for
their comments, which led to a much improved version of this paper.
We also thank
Jianqing Fan,
Alexander Giessing,
Boris Hanin,
Michael Jansson,
Jason Klusowski,
Arun Kumar,
Boris Shigida,
and Rae Yu
for comments.
 
\section*{Funding}
The authors gratefully acknowledge financial support from the National Science
Foundation through grant DMS-2210561, and Cattaneo gratefully acknowledges
financial support from the National Science Foundation through grant
SES-2241575 and from the National Institute of Health through grant
R01 GM072611-16.
 
\bibliographystyle{hapalike}
\bibliography{refs}

\clearpage
\appendix

\section{High-dimensional central limit theorems for martingales}%
\label{sec:high_dim_clt}
We present an application of our main results to
central limit theorems for high-dimensional martingale vectors. Our main
contribution in this section is found in the generality of our results,
which are broadly
applicable to martingale data and impose minimal extra assumptions. In exchange
for the scope and breadth of our results, we naturally do not necessarily
achieve state-of-the-art distributional approximation errors in certain special
cases, such as with independent data or when restricting the class of sets over
which the central limit theorem must hold. Extensions of our
results to mixingales and other approximate martingales,
along with third-order refinements and Gaussian mixture coupling distributions,
are possible through methods akin to those used to establish our main results
in Section~\ref{sec:main_results}, but we omit these for succinctness.

Our approach to deriving a high-dimensional martingale central limit theorem
proceeds as follows. Firstly, the upcoming Proposition~\ref{pro:clt} uses our
main result on martingale coupling (Corollary~\ref{cor:sa_martingale}) to
reduce the problem to that of providing anti-concentration results for
high-dimensional Gaussian vectors. We then demonstrate the utility of this
reduction by employing a few such anti-concentration methods from the existing
literature. Proposition~\ref{pro:bootstrap} gives a feasible implementation via
the Gaussian multiplier bootstrap, enabling valid
resampling-based inference using
the resulting conditional Gaussian distribution.
In \iftoggle{aos}{the supplementary material
\citep{cattaneo2025yurinskiisupplement}}{Section~\ref{sec:lp}}
we provide an example application: distributional
approximation for $\ell_p$-norms of high-dimensional martingale vectors
in Kolmogorov--Smirnov distance, relying on recent results
concerning Gaussian perimetric inequalities
\citep[see][and references therein]{nazarov2003maximal,giessing2023anti,%
cattaneo2025sharp}.

We begin with some notation. Assume the setup of
Corollary~\ref{cor:sa_martingale} and suppose $\Sigma$ is
non-random. Let $\cA$ be a class of measurable subsets of
$\R^d$ and take $T \sim \cN(0, \Sigma)$.
For $\eta>0$ and $p \in [1, \infty]$, define the Gaussian perimetric
(anti-concentration) quantity
\begin{align*}
  \Delta_p(\cA, \eta)
  &=
  \sup_{A\in \cA}
  \big\{\P(T\in A_p^\eta\setminus A)
  \vee \P(T\in A \setminus A_p^{-\eta})\big\},
\end{align*}
with $A_p^\eta = \{x \in \R^d : \|x - A\|_p \leq \eta\}$,
$A_p^{-\eta} = \R^d \setminus (\R^d \setminus A)_p^\eta$
and $\|x - A\|_p = \inf_{x' \in A} \|x - x'\|_p$.
This perimetric term allows one to convert coupling results
to central limit theorems as follows.
Denote by $\Gamma_p(\eta)$ the rate of strong approximation attained in
Corollary~\ref{cor:sa_martingale}:
\begin{align*}
  \Gamma_p(\eta)
  &=
  24 \left(
    \frac{\beta_{p,2} \phi_p(d)^2}{\eta^3}
  \right)^{1/3}
  + 17 \left(
    \frac{\E \left[ \|\Omega\|_2 \right] \phi_p(d)^2}{\eta^2}
  \right)^{1/3}.
\end{align*}

\begin{proposition}[High-dimensional central limit theorem for martingales]%
  \label{pro:clt}

  Assume the setup of Corollary~\ref{cor:sa_martingale},
  with $\Sigma$ non-random.
  For a class $\cA$ of measurable subsets of $\R^d$,
  \begin{equation}%
    \label{eq:high_dim_clt}
    \sup_{A\in \cA}
    \big|\P(S\in A) -\P(T\in A)\big|
    \leq \inf_{p \in [1, \infty]} \inf_{\eta>0}
    \big\{\Gamma_p(\eta) + \Delta_p(\cA, \eta) \big\}.
  \end{equation}
\end{proposition}

\begin{myproof}[Proposition~\ref{pro:clt}]

  \iftoggle{aos}{
    This follows directly from Strassen's theorem;
    see Lemma~SA.1 in the supplementary
    materials \citep{cattaneo2025yurinskiisupplement}.
  }{
    This follows from Strassen's theorem (Lemma~\ref{lem:strassen}), but we
    provide a proof for completeness. Note
    \begin{align*}
      \P(S \in A)
      &\leq
      \P(T \in A)
      + \P(T \in A_p^\eta \setminus A)
      + \P(\|S - T\| > \eta)
    \end{align*}
    and applying this to $\R^d \setminus A$ gives
    \begin{align*}
      \P(S\in A)
      &=
      1 - \P(S\in \R^d \setminus A) \\
      &\geq
      1 - \P(T \in \R^d \setminus A)
      - \P(T \in (\R^d \setminus A)_p^\eta \setminus (\R^d \setminus A))
      - \P(\|S - T\| > \eta) \\
      &=
      \P(T \in A)
      - \P(T \in A \setminus A_p^{-\eta})
      - \P(\|S - T\| > \eta).
    \end{align*}
    Since this holds for all $p \in [1, \infty]$,
    \begin{align*}
      \sup_{A\in \cA}
      \big|\P(S\in A) -\P(T\in A)\big|
      &\leq
      \sup_{A \in \cA}
      \big\{\P(T \in A_p^\eta\setminus A)
      \vee \P(T \in A \setminus A_p^{-\eta})\big\} \\
      &\quad+
      \P(\|S - T\| > \eta) \\
      &\leq
      \inf_{p \in [1, \infty]} \inf_{\eta>0}
      \big\{\Gamma_p(\eta) + \Delta_p(\cA, \eta) \big\}.
    \end{align*}%
  }%
\end{myproof}

The term $\Delta_p(\cA, \eta)$ in \eqref{eq:high_dim_clt}
depends on the law of $S$ only through the covariance matrix $\Sigma$, and
can be bounded using a selection of
different results from the literature.
For instance, with
$\cA = \cC = \{A \subseteq \R^d \text{ is convex}\}$,
\citet{nazarov2003maximal} showed
\begin{equation}%
  \label{eq:convex_anticonc}
  \Delta_2(\cC, \eta)
  \asymp
  \eta\sqrt{\|\Sigma^{-1}\|_{\rF}},
\end{equation}
if $\Sigma$ is invertible.
Then Proposition~\ref{pro:clt} with $p=2$
combined with \eqref{eq:convex_anticonc} yields for convex sets
\begin{align*}
  \sup_{A\in \cC}
  \big|\P(S\in A) -\P(T\in A)\big|
  &\lesssim
  \inf_{\eta > 0}
  \left\{
    \left(\frac{\beta_{p,2} d}{\eta^3}\right)^{1/3}
    + \left(\frac{\E[\|\Omega \|_2] d}{\eta^2}\right)^{1/3}
    + \eta \sqrt{\|\Sigma^{-1}\|_\rF}
  \right\}.
\end{align*}

Alternatively, with $\cA = \cR$,
the set of axis-aligned rectangles in $\R^d$, Nazarov
\citep{nazarov2003maximal,chernozhukov2017central} gives
\begin{align}%
  \label{eq:rect_anticonc}
  \Delta_\infty(\cR, \eta)
  \leq \frac{\eta (\sqrt{2\log d} + 2)}{\sigma_{\min}}
\end{align}
whenever $\min_j \, \Sigma_{j j} \geq \sigma_{\min}^2 > 0$.
Proposition~\ref{pro:clt} with $p = \infty$
and \eqref{eq:rect_anticonc} then yields
\begin{align*}%
  &\sup_{A\in \cR}
  \big|\P(S\in A) -\P(T\in A)\big| \\
  &\quad\lesssim
  \inf_{\eta > 0}
  \left\{
    \left(\frac{\beta_{\infty,2} \log 2d}{\eta^3}\right)^{1/3}
    + \left(\frac{\E[\|\Omega \|_2] \log 2d}{\eta^2}\right)^{1/3}
    + \frac{\eta \sqrt{\log 2d}}{\sigma_{\min}}
  \right\}.
\end{align*}
In situations where
$\liminf_n \min_j \, \Sigma_{j j} = 0$,
it may be possible in certain cases to regularize
the minimum variance away from zero and then apply
a Gaussian--Gaussian rectangular approximation result
such as Lemma~2.1 from \citet{chernozhukov2023nearly};
we delegate this to future work.

\begin{remark}[Comparisons with the literature]

  The literature on high-dimensional central limit theorems
  has developed rapidly in recent years
  \citep[see][and references therein]{%
    buzun2022strong,%
    lopes2022central,%
    chernozhukov2023nearly,%
    kock2024remark%
  },
  particularly for the special case of
  sums of independent random vectors
  on rectangular sets $\cR$.
  As a consequence, the results in this appendix are weaker in terms of
  dependence on the dimension than those available in the literature.
  This is an inherent issue due to our approach
  of first considering the class of all Borel sets
  and only afterwards specializing to the smaller class $\cR$. In
  contrast, sharper results in the literature, for example, directly target the
  Kolmogorov--Smirnov distance via Stein's method and Slepian interpolation.
  The main contribution of this section is therefore to obtain
  Gaussian distributional approximations
  for high-dimensional martingale vectors,
  a setting in which alternative proof strategies are not
  available.
\end{remark}

As our final main result,
we present a version of Proposition~\ref{pro:clt} in which the covariance
matrix $\Sigma$ is replaced by an estimator $\hat \Sigma$. This ensures that
the associated conditionally Gaussian vector is feasible and can be resampled,
allowing Monte Carlo quantile estimation via a Gaussian
multiplier bootstrap.

\begin{proposition}[Bootstrap central limit theorem for martingales]%
  \label{pro:bootstrap}

  Assume the setup of Corollary~\ref{cor:sa_martingale},
  with $\Sigma$ non-random,
  and let $\hat \Sigma$ be an $\bX$-measurable random
  $d \times d$ positive semi-definite matrix,
  where $\bX = (X_1, \ldots, X_n)$.
  For a class $\cA$ of measurable subsets of $\R^d$,
  \begin{align*}
    &\sup_{A\in \cA}
    \left|
    \P\big(S \in A\big)
    - \P\big(\hat \Sigma^{1/2} Z \in A \bigm| \bX \big)
    \right| \\[-1mm]
    &\quad\leq
    \inf_{p \in [1,\infty]} \inf_{\eta>0}
    \left\{ \Gamma_p(\eta) + 2 \Delta_p(\cA, \eta)
      + 2d \exp\left(\frac{-\eta^2}
        {2d^{2/p}\big\|\hat \Sigma^{1/2} - \Sigma^{1/2}\big\|_2^2}
      \right)
    \right\},
  \end{align*}
  where $Z \sim \cN(0,I_d)$ is independent of $\bX$.
\end{proposition}

\begin{myproof}[Proposition~\ref{pro:bootstrap}]

  Since $\Sigma^{1/2} Z$ is independent of $\bX$,
  we have
  \iftoggle{aos}{
    $\big|
    \P(S \in A)
    - \P\big(\hat \Sigma^{1/2} Z \in A \bigm| \bX\big)
    \big| \leq
    \big|
    \P(S \in A)
    - \P\big(\Sigma^{1/2} Z \in A\big)
    \big|
    +\big|
    \P\big(\Sigma^{1/2} Z \in A\big)
    - \P\big(\hat \Sigma^{1/2} Z \in A \bigm| \bX\big)
    \big|$.
  }{
    \begin{align*}
      \big|
      \P(S \in A)
      - \P\big(\hat \Sigma^{1/2} Z \in A \bigm| \bX\big)
      \big|
      &\leq
      \big|
      \P(S \in A)
      - \P\big(\Sigma^{1/2} Z \in A\big)
      \big| \\
      &\quad+
      \big|
      \P\big(\Sigma^{1/2} Z \in A\big)
      - \P\big(\hat \Sigma^{1/2} Z \in A \bigm| \bX\big)
      \big|.
    \end{align*}
  }
  The first term is bounded by Proposition~\ref{pro:clt};
  the second by
  \iftoggle{aos}{Lemma~SA.5
  in \citep{cattaneo2025yurinskiisupplement}}{%
  Lemma~\ref{lem:feasible_gaussian}}
  conditional on $\bX$.
  \iftoggle{aos}{}{
    Thus
    \begin{align*}
      &\left|
      \P\big(S \in A\big)
      - \P\left(\hat \Sigma^{1/2} Z \in A \bigm| \bX\right)
      \right| \\
      &\quad\leq
      \Gamma_p(\eta) + \Delta_p(\cA, \eta)
      + \Delta_{p'}(\cA, \eta')
      + 2 d \exp \left( \frac{-\eta'^2}
        {2 d^{2/p'} \big\|\hat\Sigma^{1/2} - \Sigma^{1/2}\big\|_2^2}
      \right)
    \end{align*}
    for all $A \in \cA$,
    any $p, p' \in [1, \infty]$ and each $\eta, \eta' > 0$.
  }
  Taking a supremum over $A$ and infima over
  \iftoggle{aos}{$p$ and $\eta$}{$p = p'$ and $\eta = \eta'$}
  yields the result.
\end{myproof}

A natural choice for $\hat\Sigma$ in certain situations is the sample
covariance matrix $\sum_{i=1}^n X_i X_i^\T$, or a correlation-corrected variant
thereof. In general, whenever $\hat \Sigma$ does not depend on unknown
quantities, one can sample from the law of $\hat T = \hat\Sigma^{1/2} Z$
conditional on $\bX$ to approximate the distribution of $S$.
Proposition~\ref{pro:bootstrap} verifies that this Gaussian multiplier
bootstrap approach is valid whenever $\hat\Sigma$ and $\Sigma$ are sufficiently
close. To this end, Theorem~X.1.1 in \citet{bhatia1997matrix} gives
$\big\|\hat\Sigma^{1/2} - \Sigma^{1/2}\big\|_2
\leq \big\|\hat\Sigma - \Sigma\big\|_2^{1/2}$
and Problem~X.5.5 in the same gives
$\big\|\hat\Sigma^{1/2} - \Sigma^{1/2}\big\|_2
\leq \big\|\Sigma^{-1/2}\big\|_2 \big\|\hat\Sigma - \Sigma\big\|_2$
when $\Sigma$ is invertible. The latter often gives a tighter bound when the
minimum eigenvalue of $\Sigma$ can be bounded away from zero, and consistency
of $\hat \Sigma$ can typically be established using
a range of matrix concentration inequalities.

In \iftoggle{aos}{the supplementary material
\citep{cattaneo2025yurinskiisupplement}}{Section~\ref{sec:lp}}
we apply Proposition~\ref{pro:clt} to the special case
of approximating the distribution of the $\ell_p$-norm of a high-dimensional
martingale. Proposition~\ref{pro:bootstrap} is then used to ensure that
feasible distributional approximations are also available.
 
\subsection{Distributional approximation of martingale
\texorpdfstring{$\ell_p$}{lp}-norms}
\label{sec:lp}
We present some applications of the results
derived in Appendix~\myref{sec:high_dim_clt}.
In certain empirical settings,
including nonparametric significance tests
\citep{lopes2020bootstrapping}
and nearest neighbor search procedures
\citep{biau2015high},
an estimator or test statistic
can be expressed under the null hypothesis
as the $\ell_p$-norm of a zero-mean
(possibly high-dimensional) martingale for some $p \in [1, \infty]$.
In the notation of Corollary~\myref{cor:sa_martingale},
it is therefore of interest to bound Kolmogorov--Smirnov
quantities of the form
\begin{align*}
  \sup_{t \geq 0}
  \big| \P( \|S\|_p \leq t)
  - \P( \|T\|_p \leq t) \big|.
\end{align*}
Let $\cB_p$ be the class of closed $\ell_p$-balls in $\R^d$ centered at the
origin and set
\begin{align*}
  \Delta_p(\eta)
  &\vcentcolon=
  \Delta_p(\cB_p, \eta)
  = \sup_{t \geq 0}
  \P( t < \|T\|_p \leq t + \eta ).
\end{align*}

\begin{proposition}[Distributional approximation of
  martingale $\ell_p$-norms]
  \label{pro:application_lp}

  Assume the setup of Corollary~\myref{cor:sa_martingale},
  with $\Sigma$ non-random. Then for $T \sim \cN(0, \Sigma)$,
  \begin{equation}%
    \label{eq:application_lp}
    \sup_{t \geq 0}
    \big| \P( \|S\|_p \leq t )
    - \P\left( \|T\|_p \leq t \right) \big|
    \leq \inf_{\eta>0}
    \big\{\Gamma_p(\eta) + \Delta_p(\eta) \big\}.
  \end{equation}
\end{proposition}

\begin{myproof}[Proposition~\ref{pro:application_lp}]

  Applying Proposition~\myref{pro:clt}
  with $\cA=\cB_p$ gives
  \begin{align*}
    \sup_{t \geq 0}
    \big| \P( \|S\|_p \leq t )
    - \P\left( \|T\|_p \leq t \right) \big|
    &= \sup_{A\in \cB_p}
    \big|\P(S\in A) -\P(T\in A)\big| \\
    &\leq
    \inf_{\eta>0}
    \big\{\Gamma_p(\eta) + \Delta_p(\cB_p, \eta) \big\}
    \leq
    \inf_{\eta>0}
    \big\{\Gamma_p(\eta) + \Delta_p(\eta) \big\}.
  \end{align*}
\end{myproof}

The right-hand side of
\eqref{eq:application_lp} can be controlled in various ways.
In the case of $p=\infty$,
note that $\ell_\infty$-balls are rectangles so
$\cB_\infty\subseteq \cR$, giving
$\Delta_\infty(\eta) \leq \eta (\sqrt{2\log d} + 2) / \sigma_{\min}$
whenever $\min_j \Sigma_{j j} \geq \sigma_{\min}^2$.
Alternatively, \citet[Theorem~1]{giessing2023anti} provides
$\Delta_\infty(\eta) \lesssim \eta / \sqrt{\Var[\|T\|_\infty] + \eta^2}$.
In fact, by H{\"o}lder duality of $\ell_p$-norms, we can write
$\|T\|_p = \sup_{\|u\|_q \leq 1} u^\T T$ where
$1/p + 1/q = 1$.
Then, applying the Gaussian process anti-concentration result of
\citet[Theorem~2]{giessing2023anti} yields the more general
$\Delta_p(\eta) \lesssim \eta / \sqrt{\Var[\|T\|_p] + \eta^2}$.
Thus, the problem can be reduced to that of obtaining lower bounds
for $\Var\left[\|T\|_p\right]$, with techniques for doing so
discussed, for example, in
\citet[Section~4]{giessing2023anti}.
Note that alongside the $\ell_p$-norms,
other functionals can be analyzed in this manner,
including the maximum statistic
and other order statistics
\citep{kozbur2021dimension,giessing2023anti}.

To conduct inference in this situation, we need to feasibly
approximate the quantiles of $\|T\|_p$.
To that end, take a significance level $\tau\in(0,1)$ and define
\begin{equation*}
  \hat q_p(\tau)
  = \inf \big\{t \in \R:
  \P(\|\hat T\|_p \leq t \mid \bX) \geq \tau \}
  \quad\text{where}\quad
  \hat T \mid \bX \sim \cN(0, \hat\Sigma),
\end{equation*}
with $\hat\Sigma$ any $\bX$-measurable positive semi-definite
estimator of $\Sigma$.
Note that for the canonical estimator $\hat\Sigma = \sum_{i=1}^n X_i X_i^\T$
we can write $\hat T =\sum_{i=1}^n X_i Z_i$ with
$Z_1,\dots,Z_n$ i.i.d.\ standard Gaussian independent of $\bX$,
yielding the Gaussian multiplier bootstrap.
Now assuming
the law of $\|\hat T\|_p \mid \bX$ has no atoms,
we can apply Proposition~\myref{pro:bootstrap}
to see
\begin{align*}
  &\sup_{\tau\in(0,1)}
  \big|\P\left(\|S\|_p \leq \hat q_p(\tau)\right) - \tau \big|
  \leq
  \E\left[
    \sup_{t \geq 0}
    \big|
    \P(\|S\|_p \leq t)
    - \P(\|\hat T\|_p \leq t \mid \bX)
    \big|
  \right] \\
  &\qquad\leq
  \inf_{\eta>0}
  \left\{ \Gamma_p(\eta)
    + 2 \Delta_p(\eta)
    + 2d\, \E\left[
      \exp\left(\frac{-\eta^2}
      {2d^{2/p}\big\|\hat \Sigma^{1/2} - \Sigma^{1/2}\big\|_2^2}\right)
    \right]
  \right\}
\end{align*}
and hence the bootstrap is valid whenever
$\|\hat \Sigma^{1/2} - \Sigma^{1/2}\big\|_2^2$ is sufficiently small. See the
discussion in Appendix~\myref{sec:high_dim_clt}
regarding methods for bounding this object.

\begin{remark}[One-dimensional distributional approximations]
  In our application to distributional approximation of $\ell_p$-norms,
  the object of interest $\|S\|_p$ is a
  one-dimensional functional of the high-dimensional martingale;
  contrast this with the more general Proposition~\myref{pro:clt} which
  directly considers the $d$-dimensional random vector $S$.
  As such, our coupling-based approach may be improved in certain settings
  by applying a more carefully tailored smoothing argument.
  For example, \citet{belloni2018high}
  employ a ``log sum exponential'' bound
  \citep[see also][]{chernozhukov2013gaussian}
  for the maximum statistic
  $\max_{1 \leq j \leq d} S_j$,
  along with a coupling due to \citet{chernozhukov2014gaussian}, to attain
  an improved dependence on the dimension.
  Naturally their approach does not permit the formulation of
  high-dimensional central limit theorems over arbitrary classes of
  Borel sets as in our Proposition~\myref{pro:clt}.
\end{remark}

\section{Proofs of main results}
\label{sec:proofs}

\subsection{Preliminary lemmas}

We give a sequence of preliminary lemmas which are useful for establishing our
main results. Firstly, we present a conditional version of Strassen's theorem
for the $\ell_p$-norm
\citetext{\citealp[Theorem~B.2]{chen2020jackknife};
\citealp[Theorem~4]{monrad1991nearby}}, stated for
completeness as Lemma~\ref{lem:strassen}.

\begin{lemma}[A conditional Strassen theorem for the
  \texorpdfstring{$\ell_p$}{lp}-norm]%
  \label{lem:strassen}
  Let $(\Omega, \cH, \P)$ be a probability space supporting the $\R^d$-valued
  random variable $X$ for some $d \geq 1$. Let $\cH'$ be a countably generated
  sub-$\sigma$-algebra of $\cH$ and suppose there exists a $\Unif[0,1]$ random
  variable on $(\Omega, \cH, \P)$ which is independent of the $\sigma$-algebra
  generated by $X$ and $\cH'$. Consider a regular conditional distribution
  $F(\cdot \mid \cH')$ satisfying the following. Firstly, $F(A \mid \cH')$ is
  an $\cH'$-measurable random variable for all Borel sets $A \in \cB(\R^d)$.
  Secondly, $F(\cdot \mid \cH')(\omega)$ is a Borel probability measure on
  $\R^d$ for all $\omega \in \Omega$. Taking $\eta, \rho > 0$ and
  $p \in [1, \infty]$, with $\E^*$ the outer expectation, if
  \begin{align*}
    \E^* \left[
      \sup_{A \in \cB(\R^d)}
      \Big\{
        \P \big( X \in A \mid \cH' \big)
        - F \big( A_p^\eta \mid \cH' \big)
      \Big\}
    \right]
    \leq \rho,
  \end{align*}
  where $A_p^\eta = \{x \in \R^d : \|x - A\|_p \leq \eta\}$
  and $\|x - A\|_p = \inf_{x' \in A} \|x - x'\|_p$,
  then there exists an $\R^d$-valued random variable $Y$
  with $Y \mid \cH' \sim F(\cdot \mid \cH')$
  and $\P \left( \|X-Y\|_p > \eta \right) \leq \rho$.
\end{lemma}

\begin{myproof}[Lemma~\ref{lem:strassen}]
  By Theorem~B.2 in \citet{chen2020jackknife}, noting that the $\sigma$-algebra
  generated by $Z$ is countably generated and using the metric induced by the
  $\ell_p$-norm.
\end{myproof}

Next, we present in Lemma~\ref{lem:smooth_approximation} an analytic result
concerning the smooth approximation of Borel set indicator functions, similar
to that given in \citet[Lemma~39]{belloni2019conditional}.

\begin{lemma}[Smooth approximation of Borel indicator functions]%
  \label{lem:smooth_approximation}
  Let $A \subseteq \R^d$ be a Borel set and $Z \sim \cN(0, I_d)$.
  For $\sigma, \eta > 0$ and $p \in [1, \infty]$, define
  \begin{align*}
    g_{A\eta}(x)
    &=
    \left( 1 - \frac{\|x-A^\eta\|_p}{\eta} \right) \vee 0
    &                  &\text{and}
    &f_{A\eta\sigma}(x)
    &=
    \E\big[g_{A\eta}(x + \sigma Z) \big].
  \end{align*}
  Then $f$ is infinitely differentiable
  and with $\varepsilon = \P(\|Z\|_p > \eta / \sigma)$,
  for all $k \geq 0$,
  any multi-index $\kappa = (\kappa_1,\dots, \kappa_d)\in\N^d$,
  and all $x,y \in \R^d$,
  we have $|\partial^\kappa f_{A\eta\sigma}(x)| \leq
  \frac{\sqrt{\kappa!}}{\sigma^{|\kappa|}}$ and
  \begin{align*}
    &\Bigg|
    f_{A\eta\sigma}(x+y) - \sum_{|\kappa| = 0}^k
    \frac{1}{\kappa!}
    \partial^\kappa f_{A\eta\sigma}(x)
    y^\kappa
    \Bigg|
    \leq
    \frac{\|y\|_p \|y\|_2^k}{\sigma^k \eta \sqrt{k!}}, \\
    &(1 - \varepsilon) \I\big\{x \in A\big\}
    \leq f_{A\eta\sigma}(x)
    \leq \varepsilon + (1 - \varepsilon)
    \I\big\{x \in A^{3\eta}\big\}.
  \end{align*}
\end{lemma}

\begin{myproof}[Lemma~\ref{lem:smooth_approximation}]
  Drop the subscripts on $g_{A\eta}$ and $f_{A \eta \sigma}$.
  By Taylor's theorem with Lagrange remainder, for a $t \in [0,1]$,
  \begin{align*}
    \Bigg|
    f(x + y)
    - \sum_{|\kappa|=0}^{k}
    \frac{1}{\kappa!}
    \partial^{\kappa} f(x)
    y^\kappa
    \Bigg|
    \leq
    \Bigg|
    \sum_{|\kappa|=k}
    \frac{y^\kappa}{\kappa!}
    \big(
      \partial^{\kappa} f(x + t y)
      - \partial^{\kappa} f(x)
    \big)
    \Bigg|.
  \end{align*}
  Now with $\phi(x) = \frac{1}{\sqrt{2 \pi}} e^{-x^2/2}$,
  \begin{align*}
    f(x)
    &=
    \E\big[g(x + \sigma W) \big]
    =
    \int_{\R^d}
    g(x + \sigma u)
    \prod_{j=1}^{d}
    \phi(u_j)
    \diff u
    =
    \frac{1}{\sigma^d}
    \int_{\R^d}
    g(u)
    \prod_{j=1}^{d}
    \phi \left( \frac{u_j-x_j}{\sigma} \right)
    \diff u
  \end{align*}
  and since the integrand is bounded, we exchange differentiation and
  integration to compute
  \begin{align}
    \nonumber
    \partial^\kappa
    f(x)
    &=
    \left(\frac{-1}{\sigma}\right)^{|\kappa|}
    \frac{1}{\sigma^d}
    \int_{\R^d}
    g(u)
    \prod_{j=1}^{d}
    \partial^{\kappa_j}
    \phi \left( \frac{u_j-x_j}{\sigma} \right)
    \diff u \\
    \nonumber
    &=
    \left( \frac{-1}{\sigma} \right)^{|\kappa|}
    \hspace*{-1mm}
    \int_{\R^d}
    g(x + \sigma u)
    \prod_{j=1}^{d}
    \partial^{\kappa_j}
    \phi(u_j)
    \diff u \\
    \label{eq:smoothing_derivative}
    &=
    \left( \frac{-1}{\sigma} \right)^{|\kappa|}
    \E \Bigg[
      g(x + \sigma Z)
      \prod_{j=1}^{d}
      \frac{\partial^{\kappa_j}\phi(Z_j)}{\phi(Z_j)}
    \Bigg],
  \end{align}
  where $Z \sim \cN(0, I_d)$.
  Recalling that $|g(x)| \leq 1$ and applying the Cauchy--Schwarz inequality,
  \begin{align*}
    \left|
    \partial^\kappa
    f(x)
    \right|
    &\leq
    \frac{1}{\sigma^{|\kappa|}}
    \prod_{j=1}^{d}
    \E \left[
      \left(
        \frac{\partial^{\kappa_j}\phi(Z_j)}{\phi(Z_j)}
      \right)^2
    \right]^{1/2}
    \leq
    \frac{1}{\sigma^{|\kappa|}}
    \prod_{j=1}^{d}
    \sqrt{\kappa_j!}
    =
    \frac{\sqrt{\kappa!}}{\sigma^{|\kappa|}},
  \end{align*}
  as the expected square of the Hermite polynomial of degree
  $\kappa_j$ against the standard Gaussian measure is $\kappa_j!$. By the
  reverse triangle inequality, $|g(x + t y) - g(x)| \leq t \|y\|_p / \eta$,
  so by \eqref{eq:smoothing_derivative},
  \begin{align*}
    &\left|
    \sum_{|\kappa|=k}
    \frac{y^\kappa}{\kappa!}
    \big(
      \partial^{\kappa} f(x + t y)
      - \partial^{\kappa} f(x)
    \big)
    \right| \\
    &\quad=
    \left|
    \sum_{|\kappa|=k}
    \frac{y^\kappa}{\kappa!}
    \frac{1}{\sigma^{|\kappa|}}
    \E \Bigg[
      \big(
        g(x + t y + \sigma Z)
        - g(x + \sigma Z)
      \big)
      \prod_{j=1}^{d}
      \frac{\partial^{\kappa_j}\phi(Z_j)}{\phi(Z_j)}
    \Bigg]
    \right| \\
    &\quad\leq
    \frac{t \|y\|_p}{\sigma^k \eta}
    \, \E \left[
      \Bigg|
      \sum_{|\kappa|=k}
      \frac{y^\kappa}{\kappa!}
      \prod_{j=1}^{d}
      \frac{\partial^{\kappa_j}\phi(Z_j)}{\phi(Z_j)}
      \Bigg|
    \right].
  \end{align*}
  Therefore by the Cauchy--Schwarz inequality,
  \begin{align*}
    &\Bigg(
      \sum_{|\kappa|=k}
      \frac{y^\kappa}{\kappa!}
      \big(
        \partial^{\kappa} f(x + t y)
        - \partial^{\kappa} f(x)
      \big)
    \Bigg)^2
    \leq
    \frac{t^2 \|y\|_p^2}{\sigma^{2k} \eta^2}
    \, \E \left[
      \Bigg(
        \sum_{|\kappa|=k}
        \frac{y^\kappa}{\kappa!}
        \prod_{j=1}^{d}
        \frac{\partial^{\kappa_j} \phi(Z_j)}{\phi(Z_j)}
      \Bigg)^2
    \right] \\
    &\quad=
    \frac{t^2 \|y\|_p^2}{\sigma^{2k} \eta^2}
    \sum_{|\kappa|=k}
    \sum_{|\kappa'|=k}
    \frac{y^{\kappa + \kappa'}}{\kappa! \kappa'!}
    \prod_{j=1}^{d}
    \, \E \left[
      \frac{\partial^{\kappa_j} \phi(Z_j)}{\phi(Z_j)}
      \frac{\partial^{\kappa'_j} \phi(Z_j)}{\phi(Z_j)}
    \right].
  \end{align*}
  Orthogonality of Hermite polynomials gives zero if
  $\kappa_j \neq \kappa'_j$. By the multinomial theorem,
  \begin{align*}
    &\left|
    f(x + y)
    - \sum_{|\kappa|=0}^{k}
    \frac{1}{\kappa!}
    \partial^{\kappa} f(x)
    y^\kappa
    \right|
    \leq
    \frac{\|y\|_p}{\sigma^k \eta}
    \Bigg(
      \sum_{|\kappa|=k}
      \frac{y^{2 \kappa}}{\kappa!}
    \Bigg)^{1/2} \\
    &\quad\leq
    \frac{\|y\|_p}{\sigma^k \eta \sqrt{k!}}
    \Bigg(
      \sum_{|\kappa|=k}
      \frac{k!}{\kappa!}
      y^{2 \kappa}
    \Bigg)^{1/2}
    \leq
    \frac{\|y\|_p \|y\|_2^k}{\sigma^k \eta \sqrt{k!}}.
  \end{align*}
  For the final result, since
  $f(x) = \E \left[ g(x + \sigma Z) \right]$ and
  $\I\big\{x \in A^\eta\big\}\leq g(x)\leq \I\big\{x \in A^{2\eta}\big\}$,
  \begin{align*}
    f(x)
    &\leq
    \P \left( x + \sigma Z \in A^{2 \eta} \right) \\
    &\leq
    \P \left( \|Z\|_p > \frac{\eta}{\sigma} \right)
    + \I \left\{ x \in A^{3 \eta} \right\}
    \P \left( \|Z\|_p \leq \frac{\eta}{\sigma} \right)
    = \varepsilon
    + (1 - \varepsilon) \I \left\{ x \in A^{3 \eta} \right\}
    \hspace*{-1mm}, \\
    f(x)
    &\geq
    \P \left( x + \sigma Z \in A^{\eta} \right)
    \leq
    \I \left\{ x \in A \right\}
    \P \left( \|Z\|_p \leq \frac{\eta}{\sigma} \right)
    = (1 - \varepsilon) \I \left\{ x \in A \right\}.
  \end{align*}
\end{myproof}

We provide a useful Gaussian inequality in Lemma~\ref{lem:gaussian_useful}
which helps bound the $\beta_{\infty,k}$ moment terms appearing in several
places throughout the paper.

\begin{lemma}[A Gaussian inequality]%
  \label{lem:gaussian_useful}

  Let $X \sim \cN(0, \Sigma)$
  where $\sigma_j^2 = \Sigma_{j j} \leq \sigma^2$ for all $1 \leq j \leq d$.
  Then
  \begin{align*}
    \E\left[
      \|X\|_2^2
      \|X\|_\infty
    \right]
    &\leq
    4 \sigma \sqrt{\log 2d}
    \,\sum_{j=1}^d \sigma_j^2
    &&\text{and}
    &\E\left[
      \|X\|_2^3
      \|X\|_\infty
    \right]
    &\leq
    8 \sigma \sqrt{\log 2d}
    \,\bigg( \sum_{j=1}^d \sigma_j^2 \bigg)^{3/2}.
  \end{align*}
\end{lemma}

\begin{myproof}[Lemma~\ref{lem:gaussian_useful}]

  By Cauchy--Schwarz, with $k \in \{2,3\}$, we have
  $\E\left[\|X\|_2^{k} \|X\|_\infty \right]
  \leq \E\big[\|X\|_2^{2k} \big]^{1/2} \E\big[\|X\|_\infty^2 \big]^{1/2}$.
  For the first term, by H{\"o}lder's inequality and the fourth and sixth
  moments of the normal distribution,
  \begin{align*}
    \E\big[\|X\|_2^4 \big]
    &=
    \E\Bigg[
      \bigg(
        \sum_{j=1}^d X_j^2
      \bigg)^2
    \Bigg]
    =
    \sum_{j=1}^d \sum_{k=1}^d
    \E\big[
      X_j^2 X_k^2
    \big]
    \leq
    \bigg(
      \sum_{j=1}^d
      \E\big[X_j^4 \big]^{\frac{1}{2}}
    \bigg)^2
    =
    3 \bigg(
      \sum_{j=1}^d
      \sigma_j^2
    \bigg)^2, \\
    \E\big[\|X\|_2^6 \big]
    &=
    \sum_{j=1}^d \sum_{k=1}^d \sum_{l=1}^d
    \E\big[
      X_j^2 X_k^2 X_l^2
    \big]
    \leq
    \bigg(
      \sum_{j=1}^d
      \E\big[X_j^6 \big]^{\frac{1}{3}}
    \bigg)^3
    =
    15 \bigg(
      \sum_{j=1}^d
      \sigma_j^2
    \bigg)^3.
  \end{align*}
  For the second term, by Jensen's inequality and the $\chi^2$ moment
  generating function,
  \begin{align*}
    \E\big[\|X\|_\infty^2 \big]
    &=
    \E\left[
      \max_{1 \leq j \leq d}
      X_j^2
    \right]
    \leq
    4 \sigma^2
    \log
    \sum_{j=1}^d
    \E\Big[
      e^{X_j^2 / (4\sigma^2)}
    \Big]
    \leq
    4 \sigma^2
    \log
    \sum_{j=1}^d
    \sqrt{2}
    \leq
    4 \sigma^2
    \log 2 d.
  \end{align*}
\end{myproof}

We provide an $\ell_p$-norm tail probability bound for Gaussian variables in
Lemma~\ref{lem:gaussian_pnorm}, motivating the definition of the term
$\phi_p(d)$.

\begin{lemma}[Gaussian \texorpdfstring{$\ell_p$}{lp}-norm bound]%
  \label{lem:gaussian_pnorm}
  Let $X \sim \cN(0, \Sigma)$ where $\Sigma \in \R^{d \times d}$
  is positive semi-definite. Then
  $\E\left[ \|X\|_p \right]
  \leq
  \phi_p(d)
  \max_{1 \leq j \leq d}
  \sqrt{\Sigma_{j j}}$
  where $\phi_p(d) = \sqrt{pd^{2/p} }$ for $p \in [1,\infty)$
  and $\phi_\infty(d) = \sqrt{2\log 2d}$.
\end{lemma}

\begin{myproof}[Lemma~\ref{lem:gaussian_pnorm}]

  For $p \in [1, \infty)$,
  as each $X_j$ is Gaussian, we have
  $\big(\E\big[|X_j|^p\big]\big)^{1/p}
  \leq \sqrt{p\, \E[X_j^2]}
  = \sqrt{p \Sigma_{j j}}$.
  Therefore
  \begin{align*}
    \E\big[\|X\|_p\big]
    &\leq
    \Bigg(\sum_{j=1}^d \E \big[ |X_j|^p \big] \Bigg)^{1/p}
    \leq \Bigg(\sum_{j=1}^d p^{p/2} \Sigma_{j j}^{p/2} \Bigg)^{1/p}
    \leq \sqrt{p d^{2/p}}
    \max_{1\leq j\leq d}
    \sqrt{\Sigma_{j j}}
  \end{align*}
  by Jensen's inequality.
  For $p=\infty$,
  with $\sigma^2 = \max_j \Sigma_{j j}$,
  for $t>0$,
  \begin{align*}
    \E\big[\|X\|_\infty \big]
    &\leq
    t
    \log
    \sum_{j=1}^d
    \E\Big[
      e^{|X_j| / t}
    \Big]
    \leq
    t
    \log
    \sum_{j=1}^d
    \E\Big[
      2 e^{X_j / t}
    \Big] \\
    &\leq t \log \Big(2 d e^{\sigma^2/(2t^2)}\Big)
    \leq t \log 2 d + \frac{\sigma^2}{2t},
  \end{align*}
  again by Jensen's inequality.
  Setting $t = \frac{\sigma}{\sqrt{2 \log 2d}}$ gives
  $\E\big[\|X\|_\infty \big] \leq \sigma \sqrt{2 \log 2d}$.
\end{myproof}

We give a Gaussian--Gaussian $\ell_p$-norm approximation
as Lemma~\ref{lem:feasible_gaussian}, useful for
ensuring approximations remain valid upon substituting
an estimator for the true variance matrix.

\begin{lemma}[Gaussian--Gaussian approximation in
  \texorpdfstring{$\ell_p$}{lp}-norm]%
  \label{lem:feasible_gaussian}

  Let $\Sigma_1, \Sigma_2 \in \R^{d \times d}$ be positive semi-definite
  and take $Z \sim \cN(0, I_d)$.
  For $p \in [1, \infty]$ we have
  \begin{align*}
    \P\left(
      \left\|
      \left(\Sigma_1^{1/2} - \Sigma_2^{1/2}\right) Z
      \right\|_p
      > t
    \right)
    &\leq
    2 d \exp \left(
      \frac{-t^2}
      {2 d^{2/p} \big\|\Sigma_1^{1/2} - \Sigma_2^{1/2}\big\|_2^2}
    \right).
  \end{align*}

\end{lemma}

\begin{myproof}[Lemma~\ref{lem:feasible_gaussian}]

  Let $\Sigma \in \R^{d \times d}$ be positive semi-definite
  and write $\sigma^2_j = \Sigma_{j j} $.
  For $p \in [1, \infty)$ by a union bound and
  Gaussian tail probabilities,
  \begin{align*}
    &\P\left(\big\| \Sigma^{1/2} Z \big\|_p > t \right)
    =
    \P\Bigg(
      \sum_{j=1}^d
      \left|
      \left(
        \Sigma^{1/2} Z
      \right)_j
      \right|^p
    > t^p \Bigg)
    \leq
    \sum_{j=1}^d
    \P\Bigg(
      \left|
      \left(
        \Sigma^{1/2} Z
      \right)_j
      \right|^p
      > \frac{t^p \sigma_j^p}{\|\sigma\|_p^p}
    \Bigg) \\
    &\quad=
    \sum_{j=1}^d
    \P\Bigg(
      \left|
      \sigma_j Z_j
      \right|^p
      > \frac{t^p \sigma_j^p}{\|\sigma\|_p^p}
    \Bigg)
    =
    \sum_{j=1}^d
    \P\left(
      \left| Z_j \right|
      > \frac{t}{\|\sigma\|_p}
    \right)
    \leq
    2 d \,
    \exp\left( \frac{-t^2}{2 \|\sigma\|_p^2} \right).
  \end{align*}
  The same result holds for $p = \infty$ since
  \begin{align*}
    \P\left(\big\| \Sigma^{1/2} Z \big\|_\infty > t \right)
    &=
    \P\left(
      \max_{1 \leq j \leq d}
      \left|
      \left(
        \Sigma^{1/2} Z
      \right)_j
      \right|
    > t \right)
    \leq
    \sum_{j=1}^d
    \P\left(
      \left|
      \left(
        \Sigma^{1/2} Z
      \right)_j
      \right|
      > t
    \right) \\
    &=
    \sum_{j=1}^d
    \P\left(
      \left|
      \sigma_j Z_j
      \right|
      > t
    \right)
    \leq
    2 \sum_{j=1}^d
    \exp\left( \frac{-t^2}{2 \sigma_j^2} \right)
    \leq
    2 d
    \exp\left( \frac{-t^2}{2 \|\sigma\|_\infty^2} \right).
  \end{align*}
  Now we apply this to the matrix
  $\Sigma = \big(\Sigma_1^{1/2} - \Sigma_2^{1/2}\big)^2$.
  For $p \in [1, \infty)$,
  \begin{align*}
    \|\sigma\|_p^p
    &=
    \sum_{j=1}^d (\Sigma_{j j})^{p/2}
    =
    \sum_{j=1}^d
    \Big(\big(\Sigma_1^{1/2} - \Sigma_2^{1/2}\big)^2\Big)_{j j}^{p/2}
    \leq
    d \max_{1 \leq j \leq d}
    \Big(\big(\Sigma_1^{1/2} - \Sigma_2^{1/2}\big)^2\Big)_{j j}^{p/2} \\
    &\leq
    d \, \Big\|\big(\Sigma_1^{1/2} - \Sigma_2^{1/2}\big)^2\Big\|_2^{p/2}
    =
    d \, \big\|\Sigma_1^{1/2} - \Sigma_2^{1/2}\big\|_2^p
  \end{align*}
  Similarly for $p = \infty$ we have
  \begin{align*}
    \|\sigma\|_\infty
    &=
    \max_{1 \leq j \leq d}
    (\Sigma_{j j})^{1/2}
    =
    \max_{1 \leq j \leq d}
    \Big(\big(\Sigma_1^{1/2} - \Sigma_2^{1/2}\big)^2\Big)_{j j}^{1/2}
    \leq
    \big\|\Sigma_1^{1/2} - \Sigma_2^{1/2}\big\|_2.
  \end{align*}
  Thus for all $p \in [1, \infty]$ we have
  $\|\sigma\|_p \leq
  d^{1/p} \big\|\Sigma_1^{1/2} - \Sigma_2^{1/2}\big\|_2$,
  with $d^{1/\infty} = 1$. Hence
  \begin{align*}
    \P\left(
      \left\|
      \left(\Sigma_1^{1/2} - \Sigma_2^{1/2}\right) Z
      \right\|_p
      > t
    \right)
    &\leq
    2 d \exp \left( \frac{-t^2}{2 \|\sigma\|_p^2} \right)
    \leq
    2 d \exp \left(
      \frac{-t^2}
      {2 d^{2/p} \big\|\Sigma_1^{1/2} - \Sigma_2^{1/2}\big\|_2^2}
    \right).
  \end{align*}
\end{myproof}

We also include, for completeness, a variance bound
(Lemma~\ref{lem:variance_mixing})
and an exponential concentration inequality
(Lemma~\ref{lem:exponential_mixing})
for $\alpha$-mixing random variables.

\begin{lemma}[Variance bounds for
  \texorpdfstring{$\alpha$}{alpha}-mixing random variables]
  \label{lem:variance_mixing}

  Let $X_1, \ldots, X_n$ be
  real-valued $\alpha$-mixing random
  variables with mixing coefficients $\alpha(j)$.
  Then
  \begin{enumerate}[label=(\roman*)]

    \item
      \label{eq:variance_mixing_bounded}
      If for constants $M_i$ we have
      $|X_i| \leq M_i$ a.s.\ then
      \begin{align*}
        \Var\left[
          \sum_{i=1}^n X_i
        \right]
        &\leq
        4 \sum_{j=1}^\infty \alpha(j)
        \sum_{i=1}^n M_i^2.
      \end{align*}

    \item
      \label{eq:variance_mixing_exponential}
      If $\alpha(j) \leq e^{-2j / C_\alpha}$ then
      for any $r>2$ there is a constant
      $C_r$ depending only on $r$
      such that
      \begin{align*}
        \Var\left[
          \sum_{i=1}^n X_i
        \right]
        &\leq
        C_r C_\alpha
        \sum_{i=1}^n
        \E\big[|X_i|^r\big]^{2/r}.
      \end{align*}
  \end{enumerate}
\end{lemma}

\begin{myproof}[Lemma~\ref{lem:variance_mixing}]

  Define
  $\alpha^{-1}(t) =
  \inf\{j \in \N : \alpha(j) \leq t\}$
  and
  $Q_i(t) = \inf\{s \in \R : \P(|X_i| > s) \leq t\}$.
  By Corollary~1.1 in \citet{rio2017asymptotic}
  and H{\"o}lder's inequality for $r > 2$,
  \begin{align*}
    \Var\left[
      \sum_{i=1}^n X_i
    \right]
    &\leq
    4 \sum_{i=1}^n
    \int_0^1 \alpha^{-1}(t)
    Q_i(t)^2 \diff{t} \\
    &\leq
    4 \sum_{i=1}^n
    \left(
      \int_0^1 \alpha^{-1}(t)^{\frac{r}{r-2}} \diff{t}
    \right)^{\frac{r-2}{r}}
    \left(
      \int_0^1 |Q_i(t)|^r \diff{t}
    \right)^{\frac{2}{r}}
    \diff{t}.
  \end{align*}
  Now note that if $U \sim \Unif[0,1]$ then
  $Q_i(U)$ has the same distribution as $X_i$.
  Therefore
  \begin{align*}
    \Var\left[
      \sum_{i=1}^n X_i
    \right]
    &\leq
    4
    \left(
      \int_0^1 \alpha^{-1}(t)^{\frac r{r-2}} \diff{t}
    \right)^{\frac{r-2}r}
    \sum_{i=1}^n
    \E[|X_i|^r]^{\frac 2 r}.
  \end{align*}
  If $\alpha(j) \leq e^{-2j/C_\alpha}$ then
  $\alpha^{-1}(t) \leq \frac{-C_\alpha \log t}{2}$
  so, for some constant
  $C_r$ depending only on $r$,
  \begin{align*}
    \Var\left[
      \sum_{i=1}^n X_i
    \right]
    \leq
    2 C_\alpha
    \left(
      \int_0^1 (-\log t)^{\frac r{r-2}} \diff{t}
    \right)^{\frac{r-2} r}
    \sum_{i=1}^n
    \E[|X_i|^r]^{\frac 2 r}
    \leq
    C_r C_\alpha
    \sum_{i=1}^n
    \E[|X_i|^r]^{\frac 2 r}.
  \end{align*}
  Alternatively, if for constants $M_i$ we have
  $|X_i| \leq M_i$ a.s.\ then
  \begin{align*}
    \Var\left[
      \sum_{i=1}^n X_i
    \right]
    &\leq
    4 \int_0^1 \alpha^{-1}(t)
    \diff{t}
    \sum_{i=1}^n M_i^2
    \leq
    4 \sum_{j=1}^\infty \alpha(j)
    \sum_{i=1}^n M_i^2.
  \end{align*}
\end{myproof}

\begin{lemma}[Exponential concentration inequalities for
  \texorpdfstring{$\alpha$}{alpha}-mixing random variables]
  \label{lem:exponential_mixing}

  Let $X_1, \ldots, X_n$ be zero-mean real-valued
  variables with $\alpha$-mixing coefficients
  $\alpha(j) \leq e^{-2 j / C_\alpha}$.

  \begin{enumerate}[label=(\roman*)]

    \item
      \label{eq:exponential_mixing_bounded}
      Suppose $|X_i| \leq M$ a.s.\ for each $1 \leq i \leq n$.
      Then for all $t > 0$ there is a constant $C_1$ with
      \begin{align*}
        \P\left(
          \left|
          \sum_{i=1}^n
          X_i
          \right|
          > C_1 M \big( \sqrt{n t}
          + (\log n)(\log \log n) t \big)
        \right)
        &\leq
        C_1 e^{-t}.
      \end{align*}
    \item
      \label{eq:exponential_mixing_bernstein}
      Suppose further
      $\sum_{j=1}^n |\Cov[X_i, X_j]| \leq \sigma^2$.
      Then for all $t > 0$ there is a constant $C_2$ with
      \begin{align*}
        \P\left(
          \left|
          \sum_{i=1}^n
          X_i
          \right|
          \geq C_2 \big( (\sigma \sqrt n + M) \sqrt t
          + M (\log n)^2 t \big)
        \right)
        &\leq
        C_2 e^{-t}.
      \end{align*}

  \end{enumerate}

\end{lemma}

\begin{myproof}[Lemma~\ref{lem:exponential_mixing}]

  We apply results from \citet{merlevede2009bernstein}, adjusting constants
  where necessary.
  \begin{enumerate}[label=(\roman*)]

    \item
      By Theorem~1 in \citet{merlevede2009bernstein},
      \begin{align*}
        \P\left(
          \left|
          \sum_{i=1}^n
          X_i
          \right|
          > t
        \right)
        &\leq
        \exp\left(
          -\frac{C_1 t^2}{n M^2 + Mt (\log n)(\log\log n)}
        \right).
      \end{align*}
      Replace $t$ by
      $M \sqrt{n t} + M (\log n)(\log \log n) t$.

    \item
      By Theorem~2 in \citet{merlevede2009bernstein},
      \begin{align*}
        \P\left(
          \left|
          \sum_{i=1}^n
          X_i
          \right|
          > t
        \right)
        &\leq
        \exp\left(
          -\frac{C_2 t^2}{n\sigma^2 + M^2 + Mt (\log n)^2}
        \right).
      \end{align*}
      Replace $t$ by
      $\sigma \sqrt n \sqrt t + M \sqrt t + M (\log n)^2 t$.
  \end{enumerate}
\end{myproof}

\subsection{Main results}

To establish Theorem~\myref{thm:sa_dependent}, we first
give the analogous result
for martingales as Lemma~\ref{lem:sa_martingale}. Our approach is similar to
that used in modern versions of Yurinskii's coupling for independent data, as
in Theorem~1 in \citet{lecam1988} and Theorem~10 in Chapter~10 of
\citet{pollard2002user}. The proof of Lemma~\ref{lem:sa_martingale} relies on
constructing a ``modified'' martingale, which is close to the original
martingale, but which has an $\cH_0$-measurable terminal quadratic variation.

\begin{lemma}[Strong approximation for vector-valued martingales]%
  \label{lem:sa_martingale}

  Let $X_1, \ldots, X_n$ be $\R^d$-valued
  square-integrable random vectors
  adapted to a countably generated
  filtration $\cH_0, \ldots, \cH_n$.
  Suppose that
  $\E[X_i \mid \cH_{i-1}] = 0$ for all $1 \leq i \leq n$
  and define the martingale $S = \sum_{i=1}^n X_i$.
  Let $V_i = \Var[X_i \mid \cH_{i-1}]$ and
  $\Omega = \sum_{i=1}^n V_i - \Sigma$
  where $\Sigma$ is a positive semi-definite
  $\cH_0$-measurable $d \times d$ random matrix.
  For each $\eta > 0$ and $p \in [1,\infty]$
  there is $T \mid \cH_0 \sim \cN(0, \Sigma)$ with
  \begin{align*}
    \P\big(\|S-T\|_p > 5\eta\big)
    &\leq
    \inf_{t>0}
    \left\{
      2 \P\big( \|Z\|_p > t \big)
      + \min\left\{
        \frac{\beta_{p,2} t^2}{\eta^3},
        \frac{\beta_{p,3} t^3}{\eta^4}
        + \frac{\pi_3 t^3}{\eta^3}
      \right\}
    \right\} \\
    \nonumber
    &\quad+
    \inf_{M \succeq 0}
    \big\{ 2\gamma(M) + \delta_p(M,\eta)
    + \varepsilon_p(M, \eta)\big\},
  \end{align*}
  where the second infimum is over all positive semi-definite
  $d \times d$ non-random matrices, and
  \begin{align*}
    \beta_{p,k}
    &=
    \sum_{i=1}^n \E\left[\| X_i \|^k_2 \| X_i \|_p
    + \|V_i^{1/2} Z_i \|^k_2 \|V_i^{1/2} Z_i \|_p \right],
    \qquad\gamma(M)
    = \P\big(\Omega \npreceq M\big), \\
    \delta_p(M,\eta)
    &=
    \P\left(
      \big\|\big((\Sigma +M)^{1/2}- \Sigma^{1/2}\big) Z\big\|_p
      \geq \eta
    \right),
    \qquad\pi_3
    =
    \sum_{i=1}^{n+m}
    \sum_{|\kappa| = 3}
    \E \Big[ \big|
      \E \left[ X_i^\kappa \mid \cH_{i-1} \right]
    \big| \Big], \\
    \varepsilon_p(M, \eta)
    &=
    \P\left(\big\| (M - \Omega)^{1/2} Z \big\|_p\geq \eta, \
    \Omega \preceq M\right),
  \end{align*}
  for $k \in \{2,3\}$, with $Z, Z_1,\dots ,Z_n$ i.i.d.\ standard Gaussian
  on $\R^d$ independent of $\cH_n$.
\end{lemma}

\begin{myproof}[Lemma~\ref{lem:sa_martingale}]
  \iftoggle{aos}{\phantom{a}}{}
  \proofparagraph{constructing a modified martingale}

  Take $M \succeq 0$ a fixed positive semi-definite
  $d \times d$ matrix.
  We start by constructing a new martingale based on $S$
  whose quadratic variation is $\Sigma + M$.
  Take $m \geq 1$ and define
  \begin{align*}
    H_k
    &=
    \Sigma
    + M
    - \sum_{i=1}^{k} V_i,
    \qquad\qquad\qquad\qquad\tau
    =
    \sup \big\{ k\in\{0,1,\dots,n\} : H_k \succeq 0 \big\}, \\
    \tilde X_i
    &=
    X_i\I\{i \leq \tau\}
    + \frac{1}{\sqrt{m}} H_\tau^{1/2} Z_i\I\{n+1 \leq i \leq n+m\},
    \qquad\qquad\tilde S
    =
    \sum_{i=1}^{n+m} \tilde X_i,
  \end{align*}
  where $Z_{n+1}, \ldots, Z_{n+m}$ is an i.i.d.\
  sequence of standard Gaussian vectors in $\R^d$
  independent of $\cH_n$,
  noting that $H_0 = \Sigma + M \succeq 0$ a.s.
  Define the filtration
  $\tilde \cH_0, \ldots, \tilde \cH_{n+m}$,
  where $\tilde \cH_i = \cH_i$ for $0 \leq i \leq n$
  and is the $\sigma$-algebra generated by
  $\cH_n$ and $Z_{n+1}, \dots, Z_{i}$ for $n+1 \leq i\leq n+m$.
  Observe that $\tau$ is a stopping time with respect to $\tilde\cH_i$
  because $H_{i+1} - H_i = -V_{i+1} \preceq 0$ almost surely,
  so $\{\tau \leq i\} = \{H_{i+1} \nsucceq 0\}$ for $0\leq i<n$.
  This depends only on $V_1, \dots, V_{i+1}$ and $\Sigma$
  which are $\tilde\cH_i$-measurable.
  Similarly, $\{\tau = n\} = \{H_n \succeq 0\} \in \tilde\cH_{n-1}$.
  Let $\tilde V_i = V_i \I\{i\leq\tau\}$ for
  $1\leq i\leq n$ and
  $\tilde V_i = H_\tau/m$ for $n+1\leq i\leq n+m$.
  Note that $\tilde X_i$ is $\tilde \cH_i$-measurable
  and $\tilde V_i$ is $\tilde \cH_{i-1}$-measurable.
  Further, $\E \big[ \tilde X_i \mid \tilde \cH_{i-1} \big] = 0$ and
  $\E \big[ \tilde X_i \tilde X_i^\T \mid \tilde \cH_{i-1} \big]
  = \tilde V_i$.

  \proofparagraph{bounding the difference between the original and
  modified martingales}

  By the triangle inequality,
  \begin{align*}
    \|S - \tilde S \|_p
    &\leq
    \left\| \sum_{i=\tau+1}^n  X_i \right\|_p
    + \left\| \frac{1}{\sqrt{m}} \sum_{i=n+1}^{n+m} H_\tau^{1/2} Z_i \right\|_p.
  \end{align*}
  The first term on the right vanishes on
  $\{\tau = n\} = \{H_n \succeq 0\} = \{\Omega \preceq M\}$.
  For the second term, note that
  $\tfrac{1}{\sqrt{m}}\sum_{i=n+1}^{n+m} H_\tau^{1/2} Z_i$
  is distributed as $H_\tau^{1/2}Z$,
  where $Z$ is an independent standard Gaussian. Also
  $\P\big( \| H_\tau^{1/2} Z \|_p > \eta \big)
  \leq \P\big( \| H_n^{1/2} Z \|_p > \eta,\, \Omega \preceq M)
  + \P\big( \Omega \npreceq M \big)$.
  Therefore
  \begin{align}%
    \nonumber
    \P\big( \| S - \tilde S \|_p > \eta\big)
    &\leq
    2 \P\big(\Omega \npreceq M \big)
    + \P\big( \| (M-\Omega)^{1/2}Z \|_p > \eta,\,
    \Omega \preceq M \big) \\
    \label{eq:approx_modified_original}
    &= 2 \gamma(M) + \varepsilon_p(M, \eta).
  \end{align}

  \proofparagraph{strong approximation of the modified martingale}

  Let $\tilde Z_1, \ldots, \tilde Z_{n+m}$ be i.i.d.\ $\cN(0, I_d)$
  and independent of $\tilde \cH_{n+m}$.
  Define $\check X_i = \tilde V_i^{1/2} \tilde Z_i$
  and $\check S = \sum_{i=1}^{n+m} \check X_i$.
  Fix a Borel set $A \subseteq \R^d$ and $\sigma, \eta > 0$ and
  let $f = f_{A\eta\sigma}$ be the function defined in
  Lemma~\ref{lem:smooth_approximation}.
  By the Lindeberg method, write the telescoping sum
  \begin{align*}
    \E\Big[f\big(\tilde S\big) - f\big(\check S\big)
    \mid \cH_0 \Big]
    &=
    \sum_{i=1}^{n+m}
    \E\Big[ f\big(Y_i + \tilde X_i\big)
      - f\big(Y_i + \check X_i\big)
    \mid \cH_0 \Big]
  \end{align*}
  where
  $Y_i = \sum_{j=1}^{i-1} \tilde X_j + \sum_{j=i+1}^{n+m} \check X_j$.
  By Lemma~\ref{lem:smooth_approximation} we have for $k \geq 0$
  \begin{align*}
    &\Bigg|
    \E\big[
      f(Y_i + \tilde X_i)
      - f(Y_i + \check X_i)
      \mid \cH_0
    \big]
    - \sum_{|\kappa| = 0}^k
    \frac{1}{\kappa!}
    \E \left[
      \partial^\kappa f(Y_i)
      \left( \tilde X_i^\kappa - \check X_i^\kappa \right)
      \bigm| \cH_0
    \right]
    \Bigg| \\
    &\quad\leq
    \frac{1}{\sigma^k \eta \sqrt{k!}}
    \E \left[
      \|\tilde X_i\|_p \|\tilde X_i\|_2^k
      + \|\check X_i\|_p \|\check X_i\|_2^k
      \bigm| \cH_0
    \right].
  \end{align*}
  With $k \in \{2, 3\}$, we bound each summand.
  With $|\kappa| = 0$ we have
  $\tilde X_i^\kappa = \check X_i^\kappa$,
  so consider $|\kappa| = 1$.
  Noting that $\sum_{i=1}^{n+m} \tilde V_i = \Sigma + M$, define
  \begin{align*}
    \tilde Y_i
    &=
    \sum_{j=1}^{i-1} \tilde X_j
    + \tilde Z_i
    \Bigg(\sum_{j=i+1}^{n+m} \tilde V_j\Bigg)^{1/2}
    =
    \sum_{j=1}^{i-1} \tilde X_j
    + \tilde Z_i
    \Bigg(\Sigma + M - \sum_{j=1}^{i} \tilde V_j\Bigg)^{1/2}
  \end{align*}
  and let $\check \cH_i$ be the $\sigma$-algebra generated by
  $\tilde \cH_{i-1}$ and $\tilde Z_i$.
  Note that $\tilde Y_i$ is $\check \cH_i$-measurable
  and that $Y_i$ and $\tilde Y_i$
  have the same distribution conditional on $\tilde \cH_{n+m}$.
  So
  \begin{align*}
    &\sum_{|\kappa| = 1}
    \frac{1}{\kappa!}
    \E\hspace*{-0.5mm}\left[
      \partial^\kappa f(Y_i)
      \big( \tilde X_i^\kappa - \check X_i^\kappa \big)
      \bigm| \cH_0
    \right]
    =
    \E \left[
      \nabla f(Y_i)^\T
      \big( \tilde X_i - \tilde V_i^{1/2} \tilde Z_i \big)
      \bigm| \cH_0
    \right] \\[-3mm]
    &\quad=
    \E \left[
      \nabla f(\tilde Y_i)^\T \tilde X_i
      \bigm| \cH_0
    \right]
    - \E \left[
      \nabla f(Y_i)^\T \tilde V_i^{1/2} \tilde Z_i
      \bigm| \cH_0
    \right] \\
    &\quad=
    \E \left[
      \nabla f(\tilde Y_i)^\T
      \E \left[
        \tilde X_i
        \mid \check \cH_i
      \right]
      \bigm| \cH_0
    \right]
    - \E \left[
      \tilde Z_i
    \right]
    \E \left[
      \nabla f(Y_i)^\T \tilde V_i^{1/2}
      \bigm| \cH_0
    \right] \\
    &\quad=
    \E \left[
      \nabla f(\tilde Y_i)^\T
      \E \left[
        \tilde X_i
        \mid \tilde \cH_{i-1}
      \right]
      \bigm| \cH_0
    \right]
    - 0
    = 0.
  \end{align*}
  Next, if $|\kappa| = 2$ then
  \begin{align*}
    &\sum_{|\kappa| = 2}
    \frac{1}{\kappa!}
    \E \left[
      \partial^\kappa f(Y_i)
      \left( \tilde X_i^\kappa - \check X_i^\kappa \right)
      \bigm| \cH_0
    \right] \\
    &\quad=
    \frac{1}{2}
    \E \left[
      \tilde X_i^\T \nabla^2 f(Y_i) \tilde X_i
      - \tilde Z_i^\T \tilde V_i^{1/2} \nabla^2 f(Y_i)
      \tilde V_i^{1/2} \tilde Z_i
      \bigm| \cH_0
    \right] \\
    &\quad=
    \frac{1}{2}
    \E \left[
      \E \left[
        \Tr \nabla^2 f(\tilde Y_i) \tilde X_i \tilde X_i^\T
        \bigm| \check \cH_i
      \right]
      \bigm| \cH_0
    \right]
    - \frac{1}{2}
    \E \left[
      \Tr \tilde V_i^{1/2} \nabla^2 f(Y_i) \tilde V_i^{1/2}
      \bigm| \cH_0
    \right]
    \E \left[
      \tilde Z_i \tilde Z_i^\T
    \right] \\
    &\quad=
    \frac{1}{2}
    \E \left[
      \Tr \nabla^2 f(Y_i)
      \E \left[
        \tilde X_i \tilde X_i^\T
        \bigm| \tilde \cH_{i-1}
      \right]
      \bigm| \cH_0
    \right]
    - \frac{1}{2}
    \E \left[
      \Tr \nabla^2 f(Y_i) \tilde V_i
      \bigm| \cH_0
    \right]
    = 0.
  \end{align*}
  Finally if $|\kappa| = 3$, then since
  $\check X_i \sim \cN(0, \tilde V_i)$
  conditional on $\tilde \cH_{n+m}$, we have by symmetry of the Gaussian
  distribution and Lemma~\ref{lem:smooth_approximation},
  \begin{align*}
    &
    \left|
    \sum_{|\kappa| = 3}
    \frac{1}{\kappa!}
    \E \left[
      \partial^\kappa f(Y_i)
      \left( \tilde X_i^\kappa - \check X_i^\kappa \right)
      \bigm| \cH_0
    \right]
    \right|
    \\
    &\quad=
    \left|
    \sum_{|\kappa| = 3}
    \frac{1}{\kappa!}
    \left(
      \E \left[
        \partial^\kappa f(\tilde Y_i)
        \E \left[ \tilde X_i^\kappa \mid \check \cH_i \right]
        \bigm| \cH_0
      \right]
      - \E \left[
        \partial^\kappa f(Y_i) \,
        \E \left[
          \check X_i^\kappa
          \bigm| \tilde \cH_{n+m}
        \right]
        \bigm| \cH_0
      \right]
    \right)
    \right|
    \\
    &\quad=
    \left|
    \sum_{|\kappa| = 3}
    \frac{1}{\kappa!}
    \E \left[
      \partial^\kappa f(Y_i) \,
      \E \left[ \tilde X_i^\kappa \mid \tilde \cH_{i-1} \right]
      \bigm| \cH_0
    \right]
    \right|
    \leq
    \frac{1}{\sigma^3}
    \sum_{|\kappa| = 3}
    \E \left[
      \left|
      \E \left[ \tilde X_i^\kappa \mid \tilde \cH_{i-1} \right]
      \right|
      \bigm| \cH_0
    \right].
  \end{align*}
  Combining these and summing over $i$ with $k=2$ shows
  \begin{align*}
    \E\left[
      f\big(\tilde S\big) - f\big(\check S\big)
      \bigm| \cH_0
    \right]
    &\leq
    \frac{1}{\sigma^2 \eta \sqrt{2}}
    \sum_{i=1}^{n+m}
    \E \left[
      \|\tilde X_i\|_p \|\tilde X_i\|_2^2
      + \|\check X_i\|_p \|\check X_i\|_2^2
      \bigm| \cH_0
    \right]
  \end{align*}
  On the other hand, taking $k = 3$ gives
  \begin{align*}
    \E\left[
      f\big(\tilde S\big) - f\big(\check S\big)
      \bigm| \cH_0
    \right]
    &\leq
    \frac{1}{\sigma^3 \eta \sqrt{6}}
    \sum_{i=1}^{n+m}
    \E \left[
      \|\tilde X_i\|_p \|\tilde X_i\|_2^3
      + \|\check X_i\|_p \|\check X_i\|_2^3
      \bigm| \cH_0
    \right] \\
    &\quad+
    \frac{1}{\sigma^3}
    \sum_{i=1}^{n+m}
    \sum_{|\kappa| = 3}
    \E \left[
      \left|
      \E \left[ \tilde X_i^\kappa \mid \tilde \cH_{i-1} \right]
      \right|
      \bigm| \cH_0
    \right].
  \end{align*}
  For $1 \leq i \leq n$ we have
  $\|\tilde X_i\| \leq \|X_i\|$
  and $\|\check X_i\| \leq \|V_i^{1/2} \tilde Z_i\|$.
  For $n+1 \leq i \leq n+m$ we have
  $\tilde X_i = H_\tau^{1/2} Z_i / \sqrt m$
  and $\check X_i = H_\tau^{1/2} \tilde Z_i / \sqrt m$
  which are equal in distribution given $\cH_0$.
  Therefore with
  \begin{align*}
    \tilde \beta_{p,k}
    &=
    \sum_{i=1}^{n}
    \E \left[
      \|X_i\|_p \|X_i\|_2^k
      + \|V_i^{1/2} Z_i\|_p \|V_i^{1/2} Z_i\|_2^k
      \bigm| \cH_0
    \right],
  \end{align*}
  we have, since $k \in \{2,3\}$,
  \begin{align*}
    &\sum_{i=1}^{n+m}
    \E \left[
      \|\tilde X_i\|_p \|\tilde X_i\|_2^k
      + \|\check X_i\|_p \|\check X_i\|_2^k
      \bigm| \cH_0
    \right]
    \leq
    \tilde\beta_{p,k}
    + \frac{2}{\sqrt m}
    \E \left[
      \|H_\tau^{1/2} Z\|_p \|H_\tau^{1/2} Z\|_2^k
      \bigm| \cH_0
    \right].
  \end{align*}
  Since $H_i$ is weakly decreasing under the
  semi-definite partial order, we have
  $H_\tau \preceq H_0 = \Sigma + M$
  implying that $|(H_\tau)_{j j}| \leq \|\Sigma + M\|_{\max}$ and
  $\E\big[|(H_\tau^{1/2} Z)_j|^3 \mid \cH_0 \big]
  \leq \sqrt{8/\pi}\, \|\Sigma + M\|_{\max}^{3/2}$.
  Hence as $p \geq 1$ and $k \in \{2,3\}$,
  \begin{align*}
    \E\left[
      \|H_\tau^{1/2}Z\|_p
      \|H_\tau^{1/2}Z\|_2^k
      \bigm| \cH_0
    \right]
    &\leq
    \E\left[\|H_\tau^{1/2} Z\|_1^{k+1}
      \bigm| \cH_0
    \right] \\
    &\leq
    d^{k+1} \max_{1\leq j\leq d}
    \E\left[|(H_\tau^{1/2} Z)_j|^{k+1}
      \bigm| \cH_0
    \right] \\
    &\leq 3 d^4 \,
    \|\Sigma + M\|_{\max}^{(k+1)/2}
    \leq 6 d^4 \,
    \|\Sigma \|_{\max}^{(k+1)/2}
    + 6 d^4 \|M\|.
  \end{align*}
  Assuming some $X_i$ is not identically zero so
  the result is non-trivial,
  and supposing that $\Sigma$ is bounded a.s.\
  (replacing $\Sigma$ by $\Sigma \cdot \I\{\|\Sigma\|_{\max} \leq C\}$
  for an appropriately large $C$ if necessary),
  take $m$ large enough that
  \begin{align}
    \label{eq:bound_extra_terms}
    \frac{2}{\sqrt m}
    \E \left[
      \|H_\tau^{1/2} Z\|_p \|H_\tau^{1/2} Z\|_2^k
      \bigm| \cH_0
    \right]
    \leq
    \frac{1}{4}
    \beta_{p,k}.
  \end{align}
  Further, if $|\kappa| = 3$ then
  $\big|\E \big[
  \tilde X_i^\kappa \mid \tilde \cH_{i-1} \big]\big|
  \leq \big| \E \left[ X_i^\kappa \mid \cH_{i-1} \right]\big|$
  for $1 \leq i \leq n$
  while by symmetry of the Gaussian distribution
  $\E \left[ \tilde X_i^\kappa \mid \tilde \cH_{i-1} \right] = 0$
  for $n+1 \leq i \leq n+m$.
  Hence with
  \begin{align*}
    \tilde \pi_3
    &=
    \sum_{i=1}^{n+m}
    \sum_{|\kappa| = 3}
    \E \Big[ \big|
      \E \left[ X_i^\kappa \mid \cH_{i-1} \right]
    \big| \mid \cH_0 \Big],
  \end{align*}
  we have
  \begin{align*}
    \E\left[
      f\big(\tilde S\big) - f\big(\check S\big)
      \bigm| \cH_0
    \right]
    &\leq
    \min \left\{
      \frac{3 \tilde \beta_{p,2}}{4 \sigma^2 \eta}
      + \frac{\beta_{p,2}}{4 \sigma^2 \eta},
      \frac{3 \tilde \beta_{p,3}}{4 \sigma^3 \eta}
      + \frac{\beta_{p,3}}{4 \sigma^3 \eta}
      + \frac{\tilde \pi_3}{\sigma^3}
    \right\}.
  \end{align*}
  Along with Lemma~\ref{lem:smooth_approximation}, and with
  $\sigma = \eta / t$ and $\varepsilon = \P(\|Z\|_p > t)$,
  we conclude that
  \begin{align*}
    &\P(\tilde S \in A \mid \cH_0)
    =
    \E\big[\I\{\tilde S \in A\} - f(\tilde S)
      \mid \cH_0
    \big]
    + \E\big[f(\tilde S) - f\big(\check S\big)
      \mid \cH_0
    \big]
    + \E \big[f\big(\check S\big)
      \mid \cH_0
    \big] \\
    &\,\leq
    \varepsilon\,\P(\tilde S \in A
    \mid \cH_0)
    + \min \left\{
      \frac{3 \tilde \beta_{p,2}}{4 \sigma^2 \eta}
      + \frac{\beta_{p,2}}{4 \sigma^2 \eta},
      \frac{3 \tilde \beta_{p,3}}{4 \sigma^3 \eta}
      + \frac{\beta_{p,3}}{4 \sigma^3 \eta}
      + \frac{\tilde \pi_3}{\sigma^3}
    \right\} \\
    &\quad+
    \varepsilon
    + (1 - \varepsilon) \P\big(\check S \in A_p^{3\eta}
      \mid \cH_0
    \big) \\
    &\,\leq
    \P\big( \check S \in A_p^{3\eta}
      \mid \cH_0
    \big)
    + 2 \P(\|Z\|_p > t)
    + \min\!\left\{
      \frac{3 \tilde \beta_{p,2} t^2}{4 \eta^3}
      + \frac{\beta_{p,2} t^2}{4 \eta^3},
      \frac{3 \tilde \beta_{p,3} t^3}{4 \eta^4}
      + \frac{\beta_{p,3} t^3}{4 \eta^4}
      + \frac{\tilde \pi_3 t^3}{\eta^3}
    \right\}.
  \end{align*}
  Taking a supremum and an outer expectation yields
  with $\beta_{p,k} = \E\big[\tilde \beta_{p,k}\big]$
  and $\pi_3 = \E[\tilde \pi_3]$,
  \begin{align*}
    &\E^* \left[
      \sup_{A \in \cB(\R^d)}
      \left\{
        \P(\tilde S \in A \mid \cH_0)
        - \P\big( \check S \in A_p^{3\eta} \mid \cH_0 \big)
      \right\}
    \right] \\
    &\quad\leq
    2 \P(\|Z\|_p > t)
    + \min \left\{
      \frac{\beta_{p,2} t^2}{\eta^3},
      \frac{\beta_{p,3} t^3}{\eta^4}
      + \frac{\pi_3 t^3}{\eta^3}
    \right\}.
  \end{align*}
  Finally, since
  $\check S = \sum_{i=1}^n \tilde V_i^{1/2} \tilde Z_i
  \sim \cN(0,\Sigma + M)$ conditional on $\cH_0$,
  the conditional Strassen theorem
  in Lemma~\ref{lem:strassen}
  ensures the existence of $\tilde S$ and
  $\tilde T \mid \cH_0 \sim \cN(0, \Sigma + M)$
  such that
  \begin{align}
    \label{eq:approx_modified_martingale}
    \P\left(\|\tilde S-\tilde T\|_p>3\eta\right)
    &\leq
    \inf_{t>0}
    \left\{
      2 \P(\|Z\|_p > t)
      + \min \left\{
        \frac{\beta_{p,2} t^2}{\eta^3},
        \frac{\beta_{p,3} t^3}{\eta^4} + \frac{\pi_3 t^3}{\eta^3}
      \right\}
    \right\},
  \end{align}
  since the infimum is attained by continuity of $\|Z\|_p$.

  \proofparagraph{conclusion}

  We show how to write
  $\tilde T = (\Sigma + M)^{1/2} W$
  where $W \sim \cN(0,I_d)$
  and use this representation to construct
  $T \mid \cH_0 \sim \cN(0, \Sigma)$.
  By the spectral theorem, let $\Sigma + M = U \Lambda U^\T$
  where $U$ is a $d \times d$ orthogonal random matrix
  and $\Lambda$ is a diagonal $d \times d$ random matrix with
  diagonal entries satisfying
  $\lambda_1 \geq \cdots \geq \lambda_r > 0$
  and $\lambda_{r+1} = \cdots = \lambda_d = 0$
  where $r = \rank (\Sigma + M)$.
  Let $\Lambda^+$ be the Moore--Penrose pseudo-inverse of $\Lambda$
  (obtained by inverting its non-zero elements) and define
  $W = U (\Lambda^+)^{1/2} U^\T \tilde T + U \tilde W$, where
  the first $r$ elements of $\tilde W$ are zero
  and the last $d-r$ elements are i.i.d.\ $\cN(0,1)$
  independent from $\tilde T$.
  Then, it is easy to check that
  $W \sim \cN(0, I_d)$ and that
  $\tilde T = (\Sigma + M)^{1/2} W$.
  Now define $T = \Sigma^{1/2} W$ so
  \begin{equation}%
    \label{eq:approx_target}
    \P\big(\|T - \tilde T\|_p > \eta\big)
    = \P\big(\big\|\big((\Sigma + M)^{1/2}
    - \Sigma^{1/2} \big) W \big\|_p>\eta \big)
    = \delta_p(M, \eta).
  \end{equation}
  Finally
  \eqref{eq:approx_modified_original},
  \eqref{eq:approx_modified_martingale},
  \eqref{eq:approx_target},
  the triangle inequality
  and a union bound conclude the proof since
  by taking an infimum over $M \succeq 0$,
  and by possibly reducing the constant of $1/4$ in
  \eqref{eq:bound_extra_terms} to account for
  this infimum being potentially unattainable,
  \begin{align*}
    \P\big(\|S-T\|_p > 5\eta\big)
    &\leq
    \P\big(\|\tilde S - \tilde T \|_p > 3\eta \big)
    +\P\big(\|S - \tilde S \|_p > \eta\big)
    +\P\big(\|T - \tilde T \|_p > \eta\big) \\
    &\leq
    \inf_{t>0}
    \left\{
      2 \P\big( \|Z\|_p > t \big)
      + \min\left\{
        \frac{\beta_{p,2} t^2}{\eta^3},
        \frac{\beta_{p,3} t^3}{\eta^4}
        + \frac{\pi_3 t^3}{\eta^3}
      \right\}
    \right\} \\
    &\quad+
    \inf_{M \succeq 0}
    \big\{ 2\gamma(M) + \delta_p(M,\eta)
    + \varepsilon_p(M, \eta)\big\}.
  \end{align*}
\end{myproof}

Applying Lemma~\ref{lem:sa_martingale}
and the martingale approximation
immediately yields
Theorem~\myref{thm:sa_dependent}.

\begin{myproof}[Theorem~\myref{thm:sa_dependent}]
  Apply Lemma~\ref{lem:sa_martingale} to
  the martingale $\sum_{i=1}^{n} \tilde X_i$,
  noting that $S - \sum_{i=1}^{n} \tilde X_i = U$.
\end{myproof}

Bounding the quantities
in Theorem~\myref{thm:sa_dependent} gives a
user-friendly version as Proposition~\myref{pro:sa_simplified}.

\begin{myproof}[Proposition~\myref{pro:sa_simplified}]

  We set $M = \nu^2 I_d$ and
  bound each term appearing on the right-hand side of
  the main inequality in Proposition~\myref{pro:sa_simplified}

  \proofparagraph{bounding $\P( \|Z\|_p > t )$}

  By Markov's inequality and Lemma~\ref{lem:gaussian_pnorm},
  we have
  $\P( \|Z\|_p > t ) \leq \E[\|Z\|_p] / t \leq \phi_p(d) / t$.

  \proofparagraph{bounding $\gamma(M)$}

  With $M = \nu^2 I_d$
  and by Markov,
  $\gamma(M) = \P\big(\Omega \npreceq M\big)
  = \P\big(\|\Omega\|_2 > \nu^2 \big)
  \leq \nu^{-2} \E[\|\Omega\|_2]$.

  \proofparagraph{bounding $\delta(M, \eta)$}

  By Markov's inequality and Lemma~\ref{lem:gaussian_pnorm},
  using
  $\max_j |M_{j j}| \leq \|M\|_2$
  for $M \succeq 0$,
  \begin{align*}
    \delta_{p}(M,\eta)
    &= \P\left(
      \big\|\big((\Sigma +M)^{1/2}- \Sigma^{1/2}\big) Z\big\|_p
      \geq \eta
    \right)
    \leq \frac{\phi_p(d)} {\eta}
    \E \left[
      \big\|(\Sigma +M)^{1/2}- \Sigma^{1/2}\big\|_2
    \right].
  \end{align*}
  For semi-definite matrices
  the eigenvalue operator commutes with smooth matrix functions so
  \begin{align*}
    \|(\Sigma +M)^{1/2}- \Sigma^{1/2}\|_2
    &=
    \max_{1 \leq j \leq d}
    \left|
    \sqrt{\lambda_j(\Sigma) + \nu^2} - \sqrt{\lambda_j(\Sigma)}
    \right|
    \leq \nu
  \end{align*}
  and hence $\delta_{p}(M,\eta) \leq \phi_p(d)\nu / \eta$.

  \proofparagraph{bounding $\varepsilon(M, \eta)$}

  Note that $(M -\Omega)^{1/2}Z$ is a centered Gaussian
  conditional on $\cH_n$,
  on the event $\{\Omega \preceq M\}$.
  We thus have by Markov's inequality,
  Lemma~\ref{lem:gaussian_pnorm}
  and Jensen's inequality that
  \begin{align*}
    \varepsilon_p(M, \eta)
    &= \P\left(\big\| (M - \Omega)^{1/2} Z \big\|_p\geq \eta, \
    \Omega \preceq M\right) \\
    &\leq
    \frac{1}{\eta}
    \E\left[
      \I\{\Omega \preceq M\}
      \E\left[
        \big\| (M - \Omega)^{1/2} Z \big\|_p
        \mid \cH_n
      \right]
    \right] \\
    &\leq
    \frac{\phi_p(d)}{\eta}
    \E\left[
      \I\{\Omega \preceq M\}
      \max_{1 \leq j \leq d}
      \sqrt{(M - \Omega)_{j j}}
    \right]
    \leq
    \frac{\phi_p(d)}{\eta}
    \E\left[
      \sqrt{\|M - \Omega\|_2}
    \right] \\
    &\leq
    \frac{\phi_p(d)}{\eta}
    \E\left[
      \sqrt{\|\Omega\|_2} + \nu
    \right]
    \leq
    \frac{\phi_p(d)}{\eta}
    \left(\sqrt{\E[\|\Omega\|_2]} + \nu \right).
  \end{align*}
  Thus by Theorem~\myref{thm:sa_dependent} and the previous parts,
  \begin{align*}
    \P\big(\|S-T\|_p > 6\eta\big)
    &\leq
    \inf_{t>0}
    \left\{
      2 \P\big(\|Z\|_p>t\big)
      + \min\left\{
        \frac{\beta_{p,2} t^2}{\eta^3},
        \frac{\beta_{p,3} t^3}{\eta^4}
        + \frac{\pi_3 t^3}{\eta^3}
      \right\}
    \right\} \\
    &\quad+
    \inf_{M \succeq 0}
    \big\{ 2\gamma(M) + \delta_p(M,\eta)
    + \varepsilon_p(M, \eta)\big\}
    +\P\big(\|U\|_p>\eta\big) \\
    &\leq
    \inf_{t>0}
    \left\{
      \frac{2 \phi_p(d)}{t}
      + \min\left\{
        \frac{\beta_{p,2} t^2}{\eta^3},
        \frac{\beta_{p,3} t^3}{\eta^4}
        + \frac{\pi_3 t^3}{\eta^3}
      \right\}
    \right\} \\
    &\quad+
    \inf_{\nu > 0}
    \left\{ \frac{2\E \left[ \|\Omega\|_2 \right]}{\nu^2}
      + \frac{2 \phi_p(d) \nu}{\eta}
    \right\}
    + \frac{\phi_p(d) \sqrt{\E \left[ \|\Omega\|_2 \right]}}{\eta}
    +\P\big(\|U\|_p>\eta\big).
  \end{align*}
  In general, set
  $t = 2^{1/3} \phi_p(d)^{1/3} \beta_{p,2}^{-1/3} \eta$
  and $\nu = \E[\|\Omega\|_2]^{1/3} \phi_p(d)^{-1/3} \eta^{1/3}$,
  replacing $\eta$ with $\eta / 6$ to see
  \begin{align*}
    \P\big(\|S-T\|_p > 6\eta\big)
    &\leq
    24 \left(
      \frac{\beta_{p,2} \phi_p(d)^2}{\eta^3}
    \right)^{1/3}
    + 17 \left(
      \frac{\E \left[ \|\Omega\|_2 \right] \phi_p(d)^2}{\eta^2}
    \right)^{1/3}
    +\P\left(\|U\|_p>\frac{\eta}{6}\right).
  \end{align*}
  Whenever $\pi_3 = 0$ we can set
  $t = 2^{1/4} \phi_p(d)^{1/4} \beta_{p,3}^{-1/4} \eta$,
  and with $\nu$ as above we obtain
  \begin{align*}
    \P\big(\|S-T\|_p > \eta\big)
    &\leq
    24 \left(
      \frac{\beta_{p,3} \phi_p(d)^3}{\eta^4}
    \right)^{1/4}
    + 17 \left(
      \frac{\E \left[ \|\Omega\|_2 \right] \phi_p(d)^2}{\eta^2}
    \right)^{1/3}
    +\P\left(\|U\|_p>\frac{\eta}{6}\right).
  \end{align*}
\end{myproof}

After establishing Proposition~\myref{pro:sa_simplified},
Corollaries~\myref{cor:sa_mixingale}, \myref{cor:sa_martingale}
and \myref{cor:sa_indep} follow as in the main text.

\begin{myproof}[Corollary~\myref{cor:sa_mixingale}]
  Proposition~\myref{pro:sa_simplified} with
  $\P ( \|U\|_p > \frac{\eta}{6} )
  \leq \frac{6}{\eta} \sum_{i=1}^{n} c_i (\zeta_{i} + \zeta_{n-i+1})$.
\end{myproof}

\begin{myproof}[Corollary~\myref{cor:sa_martingale}]
  By Proposition~\myref{pro:sa_simplified}
  with $U=0$ a.s.
\end{myproof}

\begin{myproof}[Corollary~\myref{cor:sa_indep}]
  By Corollary~\myref{cor:sa_martingale}
  with $\Omega=0$ a.s.
\end{myproof}

We conclude this section with a discussion expanding on the comments made
in Remark~\myref{rem:coupling_bounds_probability} on deriving bounds in
probability from Yurinskii's coupling. Consider for illustration the
independent data second-order result given in
Corollary~\myref{cor:sa_indep}: for each $\eta > 0$,
there exists $T_n \mid \cH_0 \sim \cN(0, \Sigma)$ satisfying
\begin{align*}
  \P\big(\|S_n-T_n\|_p > \eta\big)
  &\leq
  24 \left(
    \frac{\beta_{p,2} \phi_p(d)^2}{\eta^3}
  \right)^{1/3},
\end{align*}
where here we make explicit the dependence on the sample size $n$ for clarity.
The naive approach to converting this into a probability bound for
$\|S_n-T_n\|_p$ is to select $\eta$ to ensure the right-hand side is
of order $1$, arguing that the probability can then be made arbitrarily
small by taking, in this case, $\eta$ to be a large enough multiple of
$\beta_{p,2}^{1/3} \phi_p(d)^{2/3}$. However, the somewhat subtle mistake is
in neglecting the fact that the realization of the coupling variable $T_n$
will in general depend on $\eta$, rendering the resulting
bound invalid.
As an explicit example of this phenomenon, take $\eta > 1$ and suppose
$\|S_n - T_n(\eta)\| = \eta$ with probability $1 - 1/\eta$ and
$\|S_n - T_n(\eta)\| = n$ with probability $1/\eta$.
Then $\P\big(\|S_n - T_n(\eta)\| > \eta\big) = 1/\eta$
but it is not true for any $\eta$ that $\|S_n - T_n(\eta)\| \lesssim_\P 1$.

We propose in Remark~\myref{rem:coupling_bounds_probability} the following fix.
Instead of selecting $\eta$ to ensure the right-hand side is of order $1$,
we instead choose it so the bound converges (slowly) to zero. This is
easily achieved by taking the naive and incorrect bound and multiplying
by some divergent sequence $R_n$. The resulting inequality reads,
in the case of Corollary~\myref{cor:sa_indep} with
$\eta = \beta_{p,2}^{1/3} \phi_p(d)^{2/3} R_n$,
\begin{align*}
  \P\Big(\|S_n-T_n\|_p >
    \beta_{p,2}^{1/3} \phi_p(d)^{2/3} R_n
  \Big)
  &\leq
  \frac{24}{R_n}
  \to 0.
\end{align*}
We thus recover, for the price of a rate which is slower by an arbitrarily
small amount, a valid upper bound in probability, as we can immediately
conclude that
\begin{align*}
  \|S_n-T_n\|_p
  \lesssim_\P
  \beta_{p,2}^{1/3} \phi_p(d)^{2/3} R_n.
\end{align*}

\subsection{Strong approximation for martingale empirical processes}

We begin by presenting some calculations omitted from the main text
relating to the motivating example of kernel density estimation with
i.i.d.\ data.
First, the bias of this estimator is bounded as
\begin{align*}
  \big| \E \big[ \hat g(x) \big] - g(x) \big|
  &=
  \left|
  \int_{\frac{-x}{h}}^{\frac{1-x}{h}}
  K(\xi)
  \diff \xi
  - 1
  \right|
  \leq
  2 \int_{\frac{a}{h}}^\infty
  \frac{1}{\sqrt{2 \pi}}
  e^{-\frac{\xi^2}{2}}
  \diff \xi
  \leq
  \frac{h}{a}
  \sqrt{\frac{2}{\pi}}
  e^{-\frac{a^2}{2 h^2}}.
\end{align*}
Next, we do the calculations necessary to apply
Corollary~\myref{cor:sa_indep}.
Define
$k_{i j} = \frac{1}{n h} K \left( \frac{X_i - x_j}{h} \right)$ and
$k_i = (k_{i j} : 1 \leq j \leq N)$.
Then $\|k_i\|_\infty \leq \frac{1}{n h \sqrt{2 \pi}}$ a.s.\ and
$\E[\|k_i\|_2^2] \leq \frac{N}{n^2 h} \int_{-\infty}^\infty K(\xi)^2 \diff \xi
\leq \frac{N}{2 n^2 h \sqrt{\pi}}$.
Let $V = \Var[k_i] \in \R^{N \times N}$,
so assuming that $1/h \geq \log 2 N$,
by Lemma~\ref{lem:gaussian_useful},
\begin{align*}
  \beta_{\infty,2}
  &=
  n \E\left[\| k_i \|^2_2 \| k_i \|_\infty
  \right]
  + n \E \left[ \|V^{1/2} Z \|^2_2 \|V^{1/2} Z \|_\infty \right] \\
  &\leq
  \frac{N}{\sqrt{8} n^2 h^2 \pi}
  + \frac{4 N \sqrt{\log 2 N}}{\sqrt{8} n^2 h^{3/2} \pi^{3/4}}
  \leq
  \frac{N}{n^2 h^2}.
\end{align*}
Finally, we verify the stochastic continuity bounds.
By the Lipschitz property of $K$, it is easy to show that
for $x,x' \in \cX$ we have
$\left|\frac{1}{h} K \left( \frac{X_i - x}{h} \right)
- \frac{1}{h} K \left( \frac{X_i - x'}{h} \right)\right|
\lesssim \frac{|x-x'|}{h^2}$ almost surely, and also that
$\E \Big[ \left|\frac{1}{h} K \left( \frac{X_i - x}{h} \right)
- \frac{1}{h} K \left( \frac{X_i - x'}{h} \right)\right|^2 \Big]
\lesssim \frac{|x-x'|^2}{h^3}$.
By chaining with the Bernstein--Orlicz norm and polynomial covering numbers,
\begin{align*}
  \sup_{|x-x'| \leq \delta}
  \big\|S(x) - S(x')\big\|_\infty
  \lesssim_\P
  \delta
  \sqrt{\frac{\log n}{n h^3}}
\end{align*}
whenever $\log(N/h) \lesssim \log n$
and $n h \gtrsim \log n$.
By a Gaussian process maximal inequality
\citep[Corollary~2.2.8]{van1996weak}
the same bound holds for $T(x)$ with
\begin{align*}
  \sup_{|x-x'| \leq \delta}
  \big\|T(x) - T(x')\big\|_\infty
  \lesssim_\P
  \delta
  \sqrt{\frac{\log n}{n h^3}}.
\end{align*}

\iftoggle{aos}{}{\pagebreak}

\begin{myproof}[Lemma~\myref{lem:kde_eigenvalue}]
  For $x, x' \in [a, 1-a]$, the scaled covariance function
  of this nonparametric estimator is
  \begin{align*}
    n h\, \Cov\big[\hat g(x), \hat g(x')\big]
    &=
    \frac{1}{h}
    \E \left[
      K \left( \frac{X_i - x}{h} \right)
      K \left( \frac{X_i - x'}{h} \right)
    \right] \\
    &\quad-
    \frac{1}{h}
    \E \left[
      K \left( \frac{X_i - x}{h} \right)
    \right]
    \E \left[
      K \left( \frac{X_i - x'}{h} \right)
    \right] \\
    &=
    \frac{1}{2 \pi}
    \int_{\frac{-x}{h}}^{\frac{1-x}{h}}
    \exp \left( - \frac{t^2}{2} \right)
    \exp \left( - \frac{1}{2} \left( t + \frac{x - x'}{h} \right)^2 \right)
    \diff t
    - h I(x) I(x')
  \end{align*}
  where
  $I(x) = \frac{1}{\sqrt 2 \pi} \int_{-x/h}^{(1-x)/h} e^{-t^2/2} \diff t$.
  Completing the square and a substitution gives
  \begin{align*}
    n h\, \Cov\big[\hat g(x), \hat g(x')\big]
    &=
    \frac{1}{2 \pi}
    \exp \left( - \frac{1}{4} \left( \frac{x-x'}{h} \right)^2 \right)
    \int_{\frac{-x-x'}{2h}}^{\frac{2-x-x'}{2h}}
    \exp \left(-t^2\right)
    \diff t
    - h I(x) I(x').
  \end{align*}
  Now we show that since $x, x'$ are not too close to the boundary
  of $[0,1]$,
  the limits in the above integral can be replaced by $\pm \infty$.
  Note that $\frac{-x-x'}{2h} \leq \frac{-a}{h}$
  and $\frac{2-x-x'}{2h} \geq \frac{a}{h}$ so
  \begin{align*}
    \int_{-\infty}^{\infty}
    \exp \left(-t^2\right)
    \diff t
    - \int_{\frac{-x-x'}{2h}}^{\frac{2-x-x'}{2h}}
    \exp \left(-t^2\right)
    \diff t
    \leq
    2 \int_{a/h}^\infty
    \exp \left(-t^2\right)
    \diff t
    \leq
    \frac{h}{a}
    \exp \left(- \frac{a^2}{h^2}\right).
  \end{align*}
  Therefore since
  $\int_{-\infty}^{\infty} e^{-t^2} \diff t = \sqrt \pi$,
  \begin{align*}
    \left|
    n h\, \Cov\big[\hat g(x), \hat g(x')\big]
    - \frac{1}{2 \sqrt \pi}
    \exp \left( - \frac{1}{4} \left( \frac{x-x'}{h} \right)^2 \right)
    + h I(x) I(x')
    \right|
    \leq
    \frac{h}{2 \pi a}
    \exp \left(- \frac{a^2}{h^2}\right).
  \end{align*}
  Define the $N \times N$ matrix
  $\tilde\Sigma_{i j} = \frac{1}{2 \sqrt \pi}
  \exp \left( - \frac{1}{4} \left( \frac{x_i-x_j}{h} \right)^2 \right)$.
  By \citet[Proposition~2.4,
  Proposition~2.5 and Equation~2.10]{baxter1994norm},
  with
  $\cB_k = \big\{b \in \R^\Z :
  \sum_{i \in \Z} \I\{b_i \neq 0\} \leq k \big\}$,
  \begin{align*}
    \inf_{k \in \N}
    \inf_{b \in \R^k}
    \frac{\sum_{i=1}^k \sum_{j=1}^k b_i b_j \, e^{-\lambda(i-j)^2}}
    {\sum_{i=1}^k b_i^2}
    =
    \sqrt{\frac{\pi}{\lambda}}
    \sum_{i=-\infty}^{\infty}
    \exp \left( - \frac{(\pi e + 2 \pi i)^2}{4 \lambda} \right).
  \end{align*}
  We use Riemann sums,
  noting that $\pi e + 2 \pi x = 0$ at
  $x = -e/2 \approx -1.359$.
  Consider the substitutions
  $\Z \cap (-\infty, -3] \mapsto (-\infty, -2]$,
  $\{-2, -1\} \mapsto \{-2, -1\}$ and
  $\Z \cap [0, \infty) \mapsto [-1, \infty)$.
  \begin{align*}
    \sum_{i \in \Z}
    e^{-(\pi e + 2 \pi i)^2 / 4 \lambda}
    &\leq
    \int_{-\infty}^{-2}
    e^{ - (\pi e + 2 \pi x)^2/4 \lambda}
    \diff x
    + e^{- (\pi e - 4 \pi)^2/4 \lambda} \\
    &\quad+
    e^{ - (\pi e - 2 \pi)^2 / 4 \lambda}
    + \hspace*{-1mm} \int_{-1}^{\infty}
    e^{ -(\pi e + 2 \pi x)^2 / 4 \lambda}
    \diff x.
  \end{align*}
  Now use the substitution $t = \frac{\pi e + 2 \pi x}{2 \sqrt \lambda}$
  and suppose $\lambda < 1$, yielding
  \begin{align*}
    \sum_{i \in \Z}
    e^{-(\pi e + 2 \pi i)^2 / 4 \lambda}
    &\leq
    \frac{\sqrt \lambda}{\pi}
    \int_{-\infty}^{\frac{\pi e - 4 \pi}{2 \sqrt \lambda}}
    e^{-t^2}
    \diff t
    + e^{- (\pi e - 4 \pi)^2/4 \lambda} \\
    &\quad+
    e^{ - (\pi e - 2 \pi)^2 / 4 \lambda}
    + \frac{\sqrt \lambda}{\pi}
    \int_{\frac{\pi e - 2 \pi}{2 \sqrt \lambda}}^{\infty}
    e^{-t^2}
    \diff t \\
    &\leq
    \left( 1 + \frac{1}{\pi} \frac{\lambda}{4 \pi - \pi e} \right)
    e^{-(\pi e - 4 \pi)^2 / 4 \lambda}
    +
    \left( 1 + \frac{1}{\pi} \frac{\lambda}{\pi e - 2 \pi} \right)
    e^{- (\pi e - 2 \pi)^2 / 4 \lambda} \\
    &\leq
    \frac{13}{12}
    e^{-(\pi e - 4 \pi)^2 / 4 \lambda}
    +
    \frac{8}{7}
    e^{- (\pi e - 2 \pi)^2 / 4 \lambda}
    \leq
    \frac{9}{4}
    \exp \left( - \frac{5}{4 \lambda} \right).
  \end{align*}
  Therefore
  \begin{align*}
    \inf_{k \in \N}
    \inf_{b \in \cB_k}
    \frac{\sum_{i \in \Z} \sum_{j \in \Z} b_i b_j \, e^{-\lambda(i-j)^2}}
    {\sum_{i \in \Z} b_i^2}
    < \frac{4}{\sqrt \lambda}
    \exp \left( - \frac{5}{4 \lambda} \right)
    < 4 e^{-1/\lambda}.
  \end{align*}
  From this and since
  $\tilde\Sigma_{i j} = \frac{1}{2 \sqrt \pi} e^{-\lambda(i-j)^2}$
  with $\lambda = \frac{1}{4(N-1)^2 h^2} \leq \frac{\delta^2}{h^2}$,
  for each $h$ and some $\delta \leq h$,
  \begin{align*}
    \lambda_{\min}(\tilde\Sigma)
    &\leq
    2 e^{-h^2/\delta^2}.
  \end{align*}
  Recall that
  \begin{align*}
    \left|
    \Sigma_{i j}
    - \tilde\Sigma_{i j}
    + h I(x_i) I(x_j)
    \right|
    \leq
    \frac{h}{2 \pi a}
    \exp \left(- \frac{a^2}{h^2}\right).
  \end{align*}
  Now for any positive semi-definite $N \times N$ matrices $A$ and $B$
  and vector $v$ we have $\lambda_{\min}(A - v v^\T) \leq \lambda_{\min}(A)$
  and $\lambda_{\min}(B) \leq \lambda_{\min}(A) + \|B-A\|_2
  \leq \lambda_{\min}(A) + N \|B-A\|_{\max}$.
  Hence with $I_i = I(x_i)$,
  \begin{align*}
    \lambda_{\min}(\Sigma)
    &\leq
    \lambda_{\min}(\tilde\Sigma - h I I^\T)
    + \frac{N h}{2 \pi a}
    \exp \left(- \frac{a^2}{h^2}\right)
    \leq
    2 e^{-h^2/\delta^2}
    + \frac{h}{\pi a \delta}
    e^{-a^2 / h^2}.
  \end{align*}
\end{myproof}

\begin{myproof}[Proposition~\myref{pro:emp_proc}]

  Let $\cF_\delta$ be a $\delta$-cover of $(\cF, d)$.
  Using a union bound, we can write
  \begin{align*}
    &\P\left(\sup_{f \in \cF}
      \big| S(f) - T(f) \big|
    \geq 2t + \eta \right)
    \leq
    \P\left(\sup_{f \in \cF_\delta}
      \big| S(f) - T(f) \big|
    \geq \eta \right) \\
    &\qquad\qquad+
    \P\left(\sup_{d(f,f') \leq \delta}
      \big| S(f) - S(f') \big|
    \geq t \right)
    + \P\left(\sup_{d(f,f') \leq \delta}
      \big| T(f) - T(f') \big|
    \geq t \right).
  \end{align*}

  \proofparagraph{bounding the difference on $\cF_\delta$}

  We apply Corollary~\myref{cor:sa_martingale}
  with $p = \infty$ to the
  martingale difference sequence
  $\cF_\delta(X_i) = \big(f(X_i) : f \in \cF_\delta\big)$
  which takes values in $\R^{|\cF_\delta|}$.
  Square integrability can be assumed otherwise
  $\beta_\delta = \infty$.
  Note
  $\sum_{i=1}^n \cF_\delta(X_i) = S(\cF_\delta)$
  and $\phi_\infty(\cF_\delta) \leq \sqrt{2 \log 2 |\cF_\delta|}$.
  Therefore there exists a conditionally Gaussian vector $T(\cF_\delta)$
  with the same covariance structure as $S(\cF_\delta)$
  conditional on $\cH_0$ satisfying
  \begin{align*}
    \P\left(
      \sup_{f \in \cF_\delta}
      \big| S(f) - T(f) \big|
      \geq \eta
    \right)
    &\leq
    \frac{24\beta_\delta^{\frac{1}{3}}
    (2\log 2 |\cF_\delta|)^{\frac{1}{3}}}{\eta}
    + 17\left(\frac{\sqrt{2 \log 2 |\cF_\delta|}
    \sqrt{\E\left[\|\Omega_\delta\|_2\right]}}{\eta }\right)^{\frac{2}{3}}.
  \end{align*}

  \proofparagraph{bounding the fluctuations in $S(f)$}

  Since $\bigvvvert S(f) - S(f') \bigvvvert_\psi
  \leq L d(f,f')$,
  by Theorem~2.2.4 in \citet{van1996weak}
  \begin{align*}
    \biggvvvert
    \sup_{d(f,f') \leq \delta}
    \big| S(f) - S(f') \big|
    \biggvvvert_\psi
    &\leq
    C_\psi L
    \left(
      \int_0^\delta
      \psi^{-1}(N_\varepsilon) \diff{\varepsilon}
      + \delta \psi^{-1}(N_\delta^2)
    \right)
    = C_\psi L J_\psi(\delta).
  \end{align*}
  Then, by Markov's inequality and the definition of the Orlicz norm,
  \begin{align*}
    \P\left(
      \sup_{d(f,f') \leq \delta}
      \big| S(f) - S(f') \big|
      \geq t
    \right)
    &\leq
    \psi\left(\frac{t}{C_\psi L J_\psi(\delta)} \right)^{-1}.
  \end{align*}

  \proofparagraph{bounding the fluctuations in $T(f)$}

  By the Vorob'ev--Berkes--Philipp theorem
  \citep{dudley1999uniform},
  $T(\cF_\delta)$ extends to a conditionally Gaussian process $T(f)$.
  Firstly since
  $\bigvvvert T(f) - T(f') \bigvvvert_2 \leq L d(f,f')$
  conditionally on $\cH_0$,
  and $T(f)$ is a conditional Gaussian process, we have
  $\bigvvvert T(f) - T(f') \bigvvvert_{\psi_2}
  \leq 2 L d(f,f')$
  conditional on $\cH_0$
  by \citet[Chapter~2.2, Complement~1]{van1996weak},
  where $\psi_2(x) = \exp(x^2) - 1$.
  Thus again by Theorem~2.2.4 in \citet{van1996weak},
  again conditioning on $\cH_0$,
  \begin{align*}
    \biggvvvert
    \sup_{d(f,f') \leq \delta}
    \big| T(f) - T(f') \big|
    \biggvvvert_{\psi_2}
    &\leq
    C_1 L
    \int_0^\delta
    \sqrt{\log N_\varepsilon} \diff{\varepsilon}
    = C_1 L J_2(\delta)
  \end{align*}
  for some universal constant $C_1 > 0$,
  where we used $\psi_2^{-1}(x) = \sqrt{\log(1+x)}$
  and monotonicity of covering numbers.
  Then by Markov's inequality and the definition of the Orlicz norm,
  \begin{align*}
    \P\left(
      \sup_{d(f,f') \leq \delta}
      \big| T(f) - T(f') \big|
      \geq t
    \right)
    &\leq
    \left(
      \exp\left(
        \frac{t^2}{C_1^2 L^2 J_2(\delta)^2}
      \right) - 1
    \right)^{-1}
    \vee 1 \\
    &\leq
    2 \exp\left(
      \frac{-t^2}{C_1^2 L^2 J_2(\delta)^2}
    \right).
  \end{align*}

  \proofparagraph{conclusion}

  The result follows by scaling $t$ and $\eta$
  and enlarging constants if necessary.
\end{myproof}

\subsection{Applications to nonparametric regression}

\begin{myproof}[Proposition~\myref{pro:series}]

  We proceed according to the decomposition given in
  Section~\myref{sec:series}.
  By stationarity and Lemma~SA-2.1 in
  \citet{cattaneo2020large},
  we have $\sup_w \|p(w)\|_1 \lesssim 1$
  and also $\|H\|_1 \lesssim n/k$
  and $\|H^{-1}\|_1 \lesssim k/n$.

  \proofparagraph{bounding $\beta_{\infty,2}$ and $\beta_{\infty,3}$}

  Set $X_i = p(W_i) \varepsilon_i$
  so $S = \sum_{i=1}^n X_i$
  and set $\sigma^2_i = \sigma^2(W_i)$ and
  $V_i = \Var[X_i \mid \cH_{i-1}]
  = \sigma_i^2 p(W_i) p(W_i)^\T$.
  Recall from
  Corollary~\myref{cor:sa_martingale} that for $r \in \{2,3\}$,
  \begin{align*}
    \beta_{\infty,r}
    = \sum_{i=1}^n \E\left[\| X_i \|^r_2 \| X_i \|_\infty
    + \|V_i^{1/2} Z_i \|^r_2 \|V_i^{1/2} Z_i \|_\infty \right]
  \end{align*}
  with $Z_i \sim \cN(0,1)$ i.i.d.\ and independent of $V_i$.
  For the first term, we use
  $\sup_w \|p(w)\|_2 \lesssim 1$
  and bounded third moments of $\varepsilon_i$:
  \begin{align*}
    \E\left[ \| X_i \|^r_2 \| X_i \|_\infty \right]
    &\leq
    \E\left[ |\varepsilon_i|^3 \| p(W_i) \|^{r+1}_2 \right]
    \lesssim 1.
  \end{align*}
  For the second term, apply Lemma~\ref{lem:gaussian_useful} conditionally on
  $\cH_n$ with $\sup_w \|p(w)\|_2 \lesssim 1$ to see
  \begin{align*}
    \E\left[ \|V_i^{1/2} Z_i \|^r_2 \|V_i^{1/2} Z_i \|_\infty \right]
    &\lesssim
    \sqrt{\log 2k} \
    \E\left[
      \max_{1 \leq j \leq k}
      (V_i)_{j j}^{1/2}
      \bigg( \sum_{j=1}^k (V_i)_{j j} \bigg)^{r/2}
    \right] \\
    &\lesssim
    \sqrt{\log 2k} \
    \E\left[
      \sigma_i^{r+1}
      \max_{1 \leq j \leq k}
      p(W_i)_j
      \bigg(
        \sum_{j=1}^k
        p(W_i)_{j}^2
      \bigg)^{r/2}
    \right] \\
    &\lesssim
    \sqrt{\log 2k} \
    \E\left[
      \sigma_i^{r+1}
    \right]
    \lesssim
    \sqrt{\log 2k}.
  \end{align*}
  Putting these together yields
  $\beta_{\infty,2} \lesssim n \sqrt{\log 2k}$
  and $\beta_{\infty,3} \lesssim n \sqrt{\log 2k}$.

  \proofparagraph{bounding $\Omega$}

  Set $\Omega = \sum_{i=1}^n \big(V_i - \E[V_i] \big)$ as in
  Lemma~\ref{lem:sa_martingale} so
  \begin{align*}
    \Omega
    &= \sum_{i=1}^n
    \big(\sigma_i^2 p(W_i)p(W_i)^\T - \E\left[ \sigma_i^2 p(W_i)p(W_i)^\T
    \right]\big).
  \end{align*}
  Observe that $\Omega_{j l}$ is the sum of a zero-mean
  strictly stationary $\alpha$-mixing sequence and so $\E[\Omega_{j l}^2]
  \lesssim n$ by
  Lemma~\ref{lem:variance_mixing}\ref{eq:variance_mixing_bounded}.
  Since the basis functions
  satisfy Assumption~3 in \citet{cattaneo2020large}, $\Omega$ has a bounded
  number of non-zero entries in each row, and so by Jensen's inequality
  \begin{align*}
    \E\left[
      \|\Omega\|_2
    \right]
    &\leq
    \E\left[
      \|\Omega\|_\rF
    \right]
    \leq
    \left(
      \sum_{j=1}^k
      \sum_{l=1}^k
      \E\left[
        \Omega_{j l}^2
      \right]
    \right)^{1/2}
    \lesssim \sqrt{n k}.
  \end{align*}

  \proofparagraph{strong approximation}

  By Corollary~\myref{cor:sa_martingale} and the previous parts,
  with any sequence $R_n \to \infty$,
  \begin{align*}
    \|S  - T \|_\infty
    &\lesssim_\P
    \beta_{\infty,2}^{1/3} (\log 2k)^{1/3} R_n
    + \sqrt{\log 2k} \sqrt{\E[\|\Omega\|_2]} R_n \\
    &\lesssim_\P
    n^{1/3} \sqrt{\log 2k} R_n
    + (n k)^{1/4} \sqrt{\log 2k} R_n.
  \end{align*}
  If further $\E \left[ \varepsilon_i^3 \mid \cH_{i-1} \right] = 0$ then
  the third-order version of Corollary~\myref{cor:sa_martingale}
  applies since
  \begin{align*}
    \pi_3
    &=
    \sum_{i=1}^{n}
    \sum_{|\kappa| = 3}
    \E \Big[ \big|
      \E [ X_i^\kappa \mid \cH_{i-1} ]
    \big| \Big]
    =
    \sum_{i=1}^{n}
    \sum_{|\kappa| = 3}
    \E \Big[ \big|
      p(W_i)^\kappa \,
      \E [ \varepsilon_i^3 \mid \cH_{i-1} ]
    \big| \Big]
    = 0,
  \end{align*}
  giving
  \begin{align*}
    \|S  - T \|_\infty
    &\lesssim_\P
    \beta_{\infty,3}^{1/4} (\log 2k)^{3/8} R_n
    + \sqrt{\log 2k} \sqrt{\E[\|\Omega\|_2]} R_n
    \lesssim_\P
    (n k)^{1/4} \sqrt{\log 2k} R_n.
  \end{align*}
  By H{\"o}lder's inequality and with
  $\|H^{-1}\|_1 \lesssim k/n$ we have
  \begin{align*}
    \sup_{w \in \cW}
    \left|
    p(w)^\T H^{-1} S
    - p(w)^\T H^{-1} T
    \right|
    &\leq
    \sup_{w \in \cW}
    \|p(w)\|_1
    \|H^{-1}\|_1
    \| S - T \|_\infty
    \lesssim
    n^{-1} k
    \| S - T \|_\infty.
  \end{align*}

  \proofparagraph{convergence of $\hat H$}

  We have
  $\hat H - H = \sum_{i=1}^n \big(p(W_i)p(W_i)^\T - \E\left[
  p(W_i)p(W_i)^\T \right]\big)$.
  Observe that $(\hat H - H)_{j l}$ is the sum of
  a zero-mean strictly stationary $\alpha$-mixing sequence and so
  $\E[(\hat H - H)_{j l}^2] \lesssim n$ by
  Lemma~\ref{lem:variance_mixing}\ref{eq:variance_mixing_bounded}.
  Since the basis
  functions satisfy Assumption~3 in \citet{cattaneo2020large},
  $\hat H-H$ has a
  bounded number of non-zero entries in each row and so by Jensen's inequality
  \begin{align*}
    \E\left[
      \|\hat H-H\|_1
    \right]
    &=
    \E\left[
      \max_{1 \leq i \leq k}
      \sum_{j=1}^k
      \big|(\hat H-H)_{i j}\big|
    \right]
    \leq
    \E\left[
      \sum_{1 \leq i \leq k}
      \Bigg(
        \sum_{j=1}^k
        |(\hat H-H)_{i j}|
      \Bigg)^2
    \right]^{\frac{1}{2}}
    \lesssim \sqrt{n k}.
  \end{align*}

  \proofparagraph{bounding the matrix term}

  Note $\|\hat H^{-1}\|_1 \leq \|H^{-1}\|_1
  + \|\hat H^{-1}\|_1 \|\hat H-H\|_1 \|H^{-1}\|_1$
  so by the previous part, we deduce
  \begin{align*}
    \|\hat H^{-1}\|_1
    \leq
    \frac{\|H^{-1}\|_1}
    {1 - \|\hat H-H\|_1 \|H^{-1}\|_1}
    \lesssim_\P
    \frac{k/n}
    {1 - \sqrt{n k}\, k/n}
    \lesssim_\P
    \frac{k}{n}
  \end{align*}
  as $k^3 / n \to 0$. Also, note that by the martingale structure, since
  $p(W_i)$ is bounded and supported on a region with volume at most of the order
  $1/k$, and as $W_i$ has a Lebesgue density,
  \begin{align*}
    \Var[T_j]
    &=
    \Var[S_j]
    =
    \Var\left[
      \sum_{i=1}^n \varepsilon_i p(W_i)_j
    \right]
    =
    \sum_{i=1}^n
    \E\left[
      \sigma_i^2 p(W_i)_j^2
    \right]
    \lesssim
    \frac{n}{k}.
  \end{align*}
  So by the Gaussian maximal inequality in Lemma~\ref{lem:gaussian_pnorm},
  $\|T\|_\infty \lesssim_\P \sqrt{\frac{n \log 2k}{k}}$.
  Since $k^3/n \to 0$,
  \begin{align*}
    \sup_{w \in \cW}
    \left|
    p(w)^\T (\hat H^{-1} - H^{-1}) S
    \right|
    &\leq
    \sup_{w \in \cW}
    \|p(w)^\T\|_1
    \|\hat H^{-1}\|_1
    \|\hat H - H\|_1
    \|H^{-1}\|_1
    \|S - T\|_\infty \\
    &\quad+
    \sup_{w \in \cW}
    \|p(w)^\T\|_1
    \|\hat H^{-1}\|_1
    \|\hat H - H\|_1
    \|H^{-1}\|_1
    \|T\|_\infty \\
    &\lesssim_\P
    \frac{k}{n}
    \sqrt{n k}
    \frac{k}{n}
    \left(
      n^{1/3} \sqrt{\log 2k}
      + (n k)^{1/4} \sqrt{\log 2k}
    \right) \\
    &\quad+
    \frac{k}{n}
    \sqrt{n k}
    \frac{k}{n}
    \sqrt{\frac{n \log 2k}{k}}
    \lesssim_\P
    \frac{k^2}{n}
    \sqrt{\log 2k}.
  \end{align*}

  \proofparagraph{conclusion of the main result}

  By the previous parts,
  with $G(w) = p(w)^\T H^{-1} T$,
  \begin{align*}
    &\sup_{w \in \cW}
    \left|
    \hat\mu(w) - \mu(w)
    - p(w)^\T H^{-1} T
    \right| \\
    &\quad=
    \sup_{w \in \cW}
    \left|
    p(w)^\T H^{-1} (S - T)
    + p(w)^\T (\hat H^{-1} - H^{-1}) S
    + \Bias(w)
    \right| \\
    &\quad\lesssim_\P
    \frac{k}{n}
    \|S - T\|_\infty
    + \frac{k^2}{n} \sqrt{\log 2k}
    + \sup_{w \in \cW} |\Bias(w)| \\
    &\quad\lesssim_\P
    \frac{k}{n}
    \left( n^{1/3} \sqrt{\log 2k} + (n k)^{1/4} \sqrt{\log 2k} \right) R_n
    + \frac{k^2}{n} \sqrt{\log 2k}
    + \sup_{w \in \cW} |\Bias(w)| \\
    &\quad\lesssim_\P
    n^{-2/3} k \sqrt{\log 2k} R_n
    + n^{-3/4} k^{5/4} \sqrt{\log 2k} R_n
    + \frac{k^2}{n} \sqrt{\log 2k}
    + \sup_{w \in \cW} |\Bias(w)| \\
    &\quad\lesssim_\P
    n^{-2/3} k \sqrt{\log 2k} R_n
    + \sup_{w \in \cW} |\Bias(w)|
  \end{align*}
  since $k^3/n \to 0$.
  If further $\E \left[ \varepsilon_i^3 \mid \cH_{i-1} \right] = 0$ then
  \begin{align*}
    \sup_{w \in \cW}
    \left|
    \hat\mu(w) - \mu(w)
    - p(w)^\T H^{-1} T
    \right|
    &\lesssim_\P
    \frac{k}{n}
    \|S - T\|_\infty
    + \frac{k^2}{n} \sqrt{\log 2k}
    + \sup_{w \in \cW} |\Bias(w)| \\
    &\lesssim_\P
    n^{-3/4} k^{5/4} \sqrt{\log 2k} R_n
    + \sup_{w \in \cW} |\Bias(w)|.
  \end{align*}
  Finally, we verify the variance bounds for the Gaussian process.
  Since $\sigma^2(w)$ is bounded,
  \begin{align*}
    \Var[G(w)]
    &=
    p(w)^\T H^{-1}
    \Var\left[ \sum_{i=1}^n p(W_i) \varepsilon_i \right]
    H^{-1} p(w) \\
    &=
    p(w)^\T H^{-1}
    \E\left[\sum_{i=1}^n p(W_i) p(W_i)^\T \sigma^2(W_i) \right]
    H^{-1} p(w) \\
    &\lesssim
    \|p(w)\|_2^2 \|H^{-1}\|_2^2
    \|H\|_2
    \lesssim
    k/n.
  \end{align*}
  Similarly, since $\sigma^2(w)$ is bounded away from zero,
  \begin{align*}
    \Var[G(w)]
    &\gtrsim
    \|p(w)\|_2^2 \|H^{-1}\|_2^2
    \|H^{-1}\|_2^{-1}
    \gtrsim
    k/n.
  \end{align*}

  \proofparagraph{bounding the bias}

  We delegate the task of deriving bounds on the bias to
  \citet{cattaneo2020large}, who provide a high-level assumption on the
  approximation error in Assumption~4 and then use it to derive bias bounds in
  Section~3 of the form $\sup_{w \in \cW} |\Bias(w)| \lesssim_\P k^{-\gamma}$.
  This assumption is verified for B-splines, wavelets and piecewise
  polynomials in their supplemental appendix.
\end{myproof}

\begin{myproof}[Proposition~\myref{pro:series_feasible}]
  \iftoggle{aos}{\phantom{a}}{}
  \proofparagraph{infeasible supremum approximation}

  Provided that the bias is negligible,
  for all $s > 0$ we have
  \begin{align*}
    &\sup_{t \in \R}
    \left|
    \P\left(
      \sup_{w \in \cW}
      \left|
      \frac{\hat\mu(w)-\mu(w)}{\sqrt{\rho(w,w)}}
      \right| \leq t
    \right)
    -
    \P\left(
      \sup_{w \in \cW}
      \left|
      \frac{G(w)}{\sqrt{\rho(w,w)}}
      \right| \leq t
    \right)
    \right| \\
    &\quad\leq
    \sup_{t \in \R}
    \P\left(
      t \leq
      \sup_{w \in \cW}
      \left|
      \frac{G(w)}{\sqrt{\rho(w,w)}}
      \right|
      \leq t + s
    \right)
    +
    \P\left(
      \sup_{w \in \cW}
      \left|
      \frac{\hat\mu(w)-\mu(w)-G(w)}{\sqrt{\rho(w,w)}}
      \right| > s
    \right).
  \end{align*}
  By the Gaussian anti-concentration result given as Corollary~2.1 in
  \citet{chernozhukov2014anti} applied to a discretization of $\cW$, the first
  term is at most $s \sqrt{\log n}$ up to a constant factor, and the second
  term converges to zero whenever
  $\frac{1}{s} \left( \frac{k^3 (\log k)^3}{n} \right)^{1/6} \to 0$.
  Thus a suitable value of $s$ exists whenever $\frac{k^3(\log n)^6}{n} \to 0$.

  \proofparagraph{feasible supremum approximation}

  By \citet[Lemma~3.1]{chernozhukov2013gaussian},
  with $\rho(w,w') = \E[\hat\rho(w,w')]$,
  \begin{align*}
    &\sup_{t \in \R}
    \left|
    \P\left(
      \sup_{w \in \cW}
      \left|
      \frac{\hat G(w)}{\sqrt{\hat\rho(w,w)}}
      \right|
      \leq t \biggm| \bW, \bY
    \right)
    - \P\left(
      \left|
      \frac{G(w)}{\sqrt{\rho(w,w)}}
      \right|
      \leq t
    \right)
    \right| \\
    &\quad\lesssim_\P
    \sup_{w,w' \in \cW}
    \left|
    \frac{\hat\rho(w,w')}
    {\sqrt{\hat\rho(w,w)\hat\rho(w',w')}}
    - \frac{\rho(w,w')}
    {\sqrt{\rho(w,w)\rho(w',w')}}
    \right|^{1/3}
    (\log n)^{2/3} \\
    &\quad\lesssim_\P
    \left(\frac n k \right)^{1/3}
    \sup_{w,w' \in \cW} |\hat\rho(w,w') - \rho(w,w')|^{1/3}
    (\log n)^{2/3} \\
    &\quad\lesssim_\P
    \left( \frac{n (\log n)^2}{k} \right)^{1/3}
    \sup_{w,w' \in \cW}
    \left|
    p(w)^\T \hat H^{-1}
    \left(
      \hat{\Var}[S]
      - \Var[S]
    \right)
    \hat H^{-1} p(w')
    \right|^{1/3} \\
    &\quad\lesssim_\P
    \left( \frac{k (\log n)^2}{n} \right)^{1/3}
    \left\|
    \hat{\Var}[S]
    - \Var[S]
    \right\|_2^{1/3},
  \end{align*}
  and vanishes in probability when
  $\frac{k (\log n)^2}{n}
  \big\| \hat{\Var}[S] - \Var[S] \big\|_2 \to_\P 0$.
  For the plug-in estimator,
  \begin{align*}
    &\left\|
    \hat{\Var}[S]
    - \Var[S]
    \right\|_2
    =
    \left\|
    \sum_{i=1}^n
    p(W_i) p(W_i^\T)
    \hat\sigma^2(W_i)
    - n \E\left[
      p(W_i) p(W_i^\T)
      \sigma^2(W_i)
    \right]
    \right\|_2 \\
    &\quad\lesssim_\P
    \sup_{w \in \cW}
    |\hat{\sigma}^2(w)-\sigma^2(w)|
    \, \big\| \hat H \big\|_2 \\
    &\qquad+
    \left\|
    \sum_{i=1}^n
    p(W_i) p(W_i^\T)
    \sigma^2(W_i)
    - n \E\left[
      p(W_i) p(W_i^\T)
      \sigma^2(W_i)
    \right]
    \right\|_2 \\
    &\quad\lesssim_\P
    \frac{n}{k}
    \sup_{w \in \cW}
    |\hat{\sigma}^2(w)-\sigma^2(w)|
    + \sqrt{n k},
  \end{align*}
  where the second term is bounded by the same argument
  used to bound $\|\hat H - H\|_1$.
  Thus, the feasible approximation is valid whenever
  $(\log n)^2 \sup_{w \in \cW}
  |\hat{\sigma}^2(w)-\sigma^2(w)| \to_\P 0$
  and $\frac{k^3 (\log n)^4}{n} \to 0$.
  The validity of the uniform confidence band follows immediately.
\end{myproof}

\begin{myproof}[Proposition~\myref{pro:local_poly}]

  We apply Proposition~\myref{pro:emp_proc}
  with the metric $d(f_w, f_{w'}) = \|w-w'\|_2$
  and the function class
  \begin{align*}
    \cF
    &=
    \left\{
      (W_i, \varepsilon_i) \mapsto
      e_1^\T H(w)^{-1} K_h(W_i-w) p_h(W_i-w)
      \varepsilon_i
      :\ w \in \cW
    \right\},
  \end{align*}
  with $\psi$ chosen as a suitable Bernstein--Orlicz function.

  \proofparagraph{bounding $H(w)^{-1}$}

  Recall that
  $H(w) = \sum_{i=1}^n \E[K_h(W_i-w) p_h(W_i-w)p_h(W_i-w)^\T]$
  and let $a(w) \in \R^k$ with $\|a(w)\|_2 = 1$.
  Since the density of $W_i$ is bounded away from zero on $\cW$,
  \begin{align*}
    a(w)^\T H(w) a(w)
    &=
    n \E\left[
      \big( a(w)^\T p_h(W_i-w) \big)^2
      K_h(W_i-w)
    \right] \\
    &\gtrsim
    n \int_\cW
    \big( a(w)^\T p_h(u-w) \big)^2
    K_h(u-w)
    \diff{u} \\
    &\gtrsim
    n \int_{\frac{\cW-w}{h}}
    \big( a(w)^\T p(u) \big)^2
    K(u)
    \diff{u}.
  \end{align*}
  This is continuous in $a(w)$ on the compact set
  $\|a(w)\|_2 = 1$
  and $p(u)$ forms a polynomial basis so
  $a(w)^\T p(u)$ has finitely many zeroes.
  Since $K(u)$ is compactly supported
  and $h \to 0$,
  the above integral is eventually strictly positive
  for all $x \in \cW$,
  and hence is bounded below uniformly in $w \in \cW$
  by a positive constant.
  Therefore
  $\sup_{w \in \cW} \|H(w)^{-1}\|_2 \lesssim 1/n$.

  \proofparagraph{bounding $\beta_\delta$}

  Let $\cF_\delta$ be a $\delta$-cover of $(\cF, d)$
  with cardinality $|\cF_\delta| \asymp \delta^{-m}$
  and let
  $\cF_\delta(W_i, \varepsilon_i)
  = \big(f(W_i, \varepsilon_i) : f\in \cF_\delta\big)$.
  Define the truncated errors
  $\tilde\varepsilon_i =
  \varepsilon_i\I\{-a \log n \leq \varepsilon_i \leq b \log n\}$
  and note that
  $\E\big[e^{|\varepsilon_i|/C_\varepsilon}\big] < \infty$
  implies that
  $\P(\exists i: \tilde\varepsilon_i \neq \varepsilon_i)
  \lesssim n^{1-(a \vee b)/C_\varepsilon}$.
  Hence, by choosing $a$ and $b$ large enough,
  with high probability, we can replace all
  $\varepsilon_i$ by $\tilde\varepsilon_i$.
  Further, it is always possible to increase either $a$ or $b$
  along with some randomization to ensure that
  $\E[\tilde\varepsilon_i] = 0$.
  Since $K$ is bounded and compactly supported,
  $W_i$ has a bounded density and
  $|\tilde\varepsilon_i| \lesssim \log n$,
  \begin{align*}
    \bigvvvert
    f(W_i, \tilde\varepsilon_i)
    \bigvvvert_2
    &=
    \E\left[
      \left|
      e_1^\T H(w)^{-1} K_h(W_i-w) p_h(W_i-w)
      \tilde\varepsilon_i
      \right|^2
    \right]^{1/2} \\
    &\leq
    \E\left[
      \|H(w)^{-1}\|_2^2
      K_h(W_i-w)^2
      \|p_h(W_i-w)\|_2^2
      \sigma^2(W_i)
    \right]^{1/2} \\
    &\lesssim
    n^{-1}
    \E\left[
      K_h(W_i-w)^2
    \right]^{1/2}
    \lesssim
    n^{-1}
    h^{-m / 2}, \\
    \bigvvvert
    f(W_i, \tilde\varepsilon_i)
    \bigvvvert_\infty
    &\leq
    \bigvvvert
    \|H(w)^{-1}\|_2
    K_h(W_i-w)
    \|p_h(W_i-w)\|_2
    |\tilde\varepsilon_i|
    \bigvvvert_\infty \\
    &\lesssim
    n^{-1}
    \bigvvvert
    K_h(W_i-w)
    \bigvvvert_\infty
    \log n
    \lesssim
    n^{-1}
    h^{-m}
    \log n.
  \end{align*}
  Therefore
  \begin{align*}
    \E\left[
      \|\cF_\delta(W_i, \tilde\varepsilon_i)\|_2^2
      \|\cF_\delta(W_i, \tilde\varepsilon_i)\|_\infty
    \right]
    &\leq
    \sum_{f\in\cF_\delta}
    \bigvvvert f(W_i, \tilde\varepsilon_i) \bigvvvert_2^2
    \, \max_{f\in\cF_\delta}
    \bigvvvert f(W_i, \tilde\varepsilon_i) \bigvvvert_\infty \\
    &\lesssim
    n^{-3} \delta^{-m} h^{-2m} \log n.
  \end{align*}
  Let
  $V_i(\cF_\delta) =
  \E\big[\cF_\delta(W_i, \tilde\varepsilon_i)
    \cF_\delta(W_i, \tilde\varepsilon_i)^\T
  \mid \cH_{i-1}\big]$
  and $Z_i \sim \cN(0, I_d)$ be i.i.d.\ and
  independent of $\cH_n$.
  Note that
  $V_i(f,f) = \E[f(W_i, \tilde\varepsilon_i)^2 \mid W_i]
  \lesssim n^{-2} h^{-2m}$
  and
  $\E[V_i(f,f)] = \E[f(W_i, \tilde\varepsilon_i)^2]
  \lesssim n^{-2} h^{-m}$.
  Thus by Lemma~\ref{lem:gaussian_useful},
  \begin{align*}
    \E\left[
      \big\| V_i(\cF_\delta)^{1/2} Z_i \big\|^2_2
      \big\| V_i(\cF_\delta)^{1/2} Z_i \big\|_\infty
    \right]
    &=
    \E\left[
      \E\left[
        \big\| V_i(\cF_\delta)^{1/2} Z_i \big\|^2_2
        \big\| V_i(\cF_\delta)^{1/2} Z_i \big\|_\infty
        \mid \cH_n
      \right]
    \right] \\
    &\leq
    4 \sqrt{\log 2|\cF_\delta|}
    \,\E\Bigg[
      \max_{f \in \cF_\delta} \sqrt{V_i(f,f)}
      \sum_{f \in \cF_\delta} V_i(f,f)
    \Bigg] \\
    &\lesssim
    n^{-3}
    h^{-2m}
    \delta^{-m}
    \sqrt{\log(1/\delta)}.
  \end{align*}
  Thus since $\log(1/\delta) \asymp \log(1/h) \asymp\log n$,
  \begin{align*}
    \beta_\delta
    &=
    \sum_{i=1}^n
    \E\left[
      \|\cF_\delta(W_i, \tilde\varepsilon_i)\|_2^2
      \|\cF_\delta(W_i, \tilde\varepsilon_i)\|_\infty
      + \big\| V_i(\cF_\delta)^{1/2} Z_i \big\|^2_2
      \big\| V_i(\cF_\delta)^{1/2} Z_i \big\|_\infty
    \right] \\
    &\lesssim
    \frac{\log n}
    {n^2 h^{2m} \delta^m}.
  \end{align*}

  \proofparagraph{bounding $\Omega_\delta$}

  Let $C_K>0$ be the radius of a $\ell_2$-ball
  containing the support of $K$
  and note that
  \begin{align*}
    \left|
    V_i(f,f')
    \right|
    &=
    \Big|
    \E\Big[
      e_1^\T H(w)^{-1}
      p_h(W_i-w)
      e_1^\T H(w')^{-1}
      p_h(W_i-w') \\
      &\qquad\times
      K_h(W_i-w)
      K_h(W_i-w')
      \tilde\varepsilon_i^2
      \Bigm| \cH_{i-1}
    \Big]
    \Big| \\
    &\lesssim
    n^{-2}
    K_h(W_i-w)
    K_h(W_i-w') \\
    &\lesssim
    n^{-2}
    h^{-m}
    K_h(W_i-w)
    \I\{\|w-w'\|_2 \leq 2 C_K h\}.
  \end{align*}
  Since $W_i$ are $\alpha$-mixing
  with $\alpha(j) < e^{-2j / C_\alpha}$,
  Lemma~\ref{lem:variance_mixing}\ref{eq:variance_mixing_exponential}
  with $r=3$ gives
  \begin{align*}
    &\Var\left[
      \sum_{i=1}^n V_i(f,f')
    \right] \\
    &\quad\lesssim
    \sum_{i=1}^n
    \E\left[
      |V_i(f,f')|^3
    \right] ^{2/3}
    \lesssim
    n^{-3} h^{-2m}
    \E\left[
      K_h(W_i-w)^3
    \right] ^{2/3}
    \I\{\|w-w'\|_2 \leq 2 C_K h\} \\
    &\quad\lesssim
    n^{-3} h^{-2m}
    (h^{-2m})^{2/3}
    \I\{\|w-w'\|_2 \leq 2 C_K h\} \\
    &\quad\lesssim
    n^{-3} h^{-10m/3}
    \I\{\|w-w'\|_2 \leq 2 C_K h\}.
  \end{align*}
  Therefore, by Jensen's inequality,
  \begin{align*}
    \E\big[ \|\Omega_\delta\|_2 \big]
    &\leq
    \E\big[ \|\Omega_\delta\|_\rF \big]
    \leq
    \E\Bigg[
      \sum_{f,f' \in \cF_\delta}
      (\Omega_\delta)_{f,f'}^2
    \Bigg]^{1/2}
    \leq
    \Bigg(
      \sum_{f,f' \in \cF_\delta}
      \Var\left[
        \sum_{i=1}^n V_i(f,f')
      \right]
    \Bigg)^{1/2} \\
    &\lesssim
    n^{-3/2} h^{-5m/3}
    \Bigg(
      \sum_{f,f' \in \cF_\delta}
      \I\{\|w-w'\|_2 \leq 2 C_K h\}
    \Bigg)^{1/2} \\
    &\lesssim
    n^{-3/2} h^{-5m/3}
    \big(h^{m} \delta^{-2m} \big)^{1/2}
    \lesssim
    n^{-3/2}
    h^{-7m/6}
    \delta^{-m}.
  \end{align*}
  Note that we could have used
  $\|\cdot\|_1$ rather than $\|\cdot\|_\rF$,
  but this term is negligible either way.

  \proofparagraph{regularity of the stochastic processes}

  For each $f, f' \in \cF$,
  define the zero-mean and $\alpha$-mixing random variables
  \begin{align*}
    u_i(f,f')
    &=
    e_1^\T
    \big(
      H(w)^{-1} K_h(W_i-w) p_h(W_i-w)
      - H(w')^{-1} K_h(W_i-w') p_h(W_i-w')
    \big)
    \tilde\varepsilon_i.
  \end{align*}
  To bound this we use that for all $1 \leq j \leq k$,
  by the Lipschitz property of the kernel and monomials,
  \begin{align*}
    &\left|
    K_h(W_i-w) - K_h(W_i-w')
    \right| \\
    &\quad\lesssim
    h^{-m-1}
    \|w-w'\|_2
    \big(
      \I\{\|W_i-w\| \leq C_K h\}
      + \I\{\|W_i-w'\| \leq C_K h\}
    \big), \\
    &\left|
    p_h(W_i-w)_j - p_h(W_i-w')_j
    \right|
    \lesssim
    h^{-1}
    \|w-w'\|_2,
  \end{align*}
  to deduce that for any $1 \leq j,l \leq k$,
  \begin{align*}
    \big| H(w)_{j l} - H(w')_{j l} \big|
    &=
    \big|
    n \E\big[
      K_h(W_i-w) p_h(W_i-w)_j p_h(W_i-w)_l \\
      &\qquad-
      K_h(W_i-w') p_h(W_i-w')_j p_h(W_i-w')_l
    \big]
    \big| \\
    &\leq
    n\E\left[
      \left|
      K_h(W_i-w) - K_h(W_i-w')
      \right|
      \left|
      p_h(W_i-w)_j
      p_h(W_i-w)_l
      \right|
    \right] \\
    &\quad+
    n\E\left[
      \left|
      p_h(W_i-w)_j - p_h(W_i-w')_j
      \right|
      \left|
      K_h(W_i-w')
      p_h(W_i-w)_l
      \right|
    \right] \\
    &\quad+
    n\E\left[
      \left|
      p_h(W_i-w)_l - p_h(W_i-w')_l
      \right|
      \left|
      K_h(W_i-w')
      p_h(W_i-w')_j
      \right|
    \right] \\
    &\lesssim
    n h^{-1}\|w-w'\|_2.
  \end{align*}
  Therefore as the dimension of the matrix $H(w)$ is fixed,
  \begin{align*}
    \big\| H(w)^{-1} - H(w')^{-1} \big\|_2
    &\leq
    \big\| H(w)^{-1}\big\|_2
    \big\| H(w')^{-1}\big\|_2
    \big\| H(w) - H(w') \big\|_2
    \lesssim
    \frac{\|w-w'\|_2}{n h}.
  \end{align*}
  Hence
  \begin{align*}
    \big| u_i(f,f') \big|
    &\leq
    \big\|
    H(w)^{-1} K_h(W_i-w) p_h(W_i-w)
    - H(w')^{-1} K_h(W_i-w') p_h(W_i-w')
    \tilde\varepsilon_i
    \big\|_2 \\
    &\leq
    \big\| H(w)^{-1} - H(w')^{-1} \big\|_2
    \big\| K_h(W_i-w) p_h(W_i-w)
    \tilde\varepsilon_i
    \big\|_2 \\
    &\quad+
    \big| K_h(W_i-w) - K_h(W_i-w') \big|
    \big\| H(w')^{-1} p_h(W_i-w)
    \tilde\varepsilon_i
    \big\|_2 \\
    &\quad+
    \big\| p_h(W_i-w) - p_h(W_i-w') \big\|_2
    \big\| H(w')^{-1} K_h(W_i-w')
    \tilde\varepsilon_i \big\|_2 \\
    &\lesssim
    \frac{\|w-w'\|_2}{n h}
    \big| K_h(W_i-w) \tilde\varepsilon_i \big|
    + \frac{1}{n}
    \big| K_h(W_i-w) - K_h(W_i-w') \big|
    \,|\tilde\varepsilon_i| \\
    &\lesssim
    \frac{\|w-w'\|_2 \log n}{n h^{m+1}},
  \end{align*}
  and from the penultimate line, we also deduce that
  \begin{align*}
    \Var[u_i(f,f')]
    &\lesssim
    \frac{\|w-w'\|_2^2}{n^2h^2}
    \E\left[
      K_h(W_i-w)^2 \sigma^2(X_i)
    \right] \\
    &\quad+
    \frac{1}{n^2}
    \E\left[
      \big( K_h(W_i-w) - K_h(W_i-w') \big)^2
      \sigma^2(X_i)
    \right]
    \lesssim
    \frac{\|w-w'\|_2^2}{n^2h^{m+2}}.
  \end{align*}
  Further, $\E[u_i(f,f') u_j(f,f')] = 0$ for $i \neq j$ so
  by Lemma~\ref{lem:exponential_mixing}\ref{eq:exponential_mixing_bernstein},
  for a constant $C_1>0$,
  \begin{align*}
    \P\left(
      \Big| \sum_{i=1}^n u_i(f,f') \Big|
      \geq \frac{C_1 \|w-w'\|_2}{\sqrt n h^{m/2+1}}
      \left(
        \sqrt{t}
        + \sqrt{\frac{(\log n)^2}{n h^m}} \sqrt t
        + \sqrt{\frac{(\log n)^6}{n h^m}} t
      \right)
    \right)
    &\leq
    C_1 e^{-t}.
  \end{align*}
  Therefore, adjusting the constant if necessary
  and since $n h^{m} \gtrsim (\log n)^7$,
  \begin{align*}
    \P\left(
      \Big| \sum_{i=1}^n u_i(f,f') \Big|
      \geq
      \frac{C_1 \|w-w'\|_2}{\sqrt{n} h^{m/2+1}}
      \left(
        \sqrt{t} + \frac{t}{\sqrt{\log n}}
      \right)
    \right)
    &\leq
    C_1 e^{-t}.
  \end{align*}
  By Lemma~2 in \citet{van2013bernstein} with
  $\psi(x) =
  \exp\Big(\big(\sqrt{1+2 x / \sqrt{\log n}}-1 \big)^2
  \log n \Big)-1$,
  \begin{align*}
    \Bigvvvert \sum_{i=1}^n u_i(f,f') \Bigvvvert_\psi
    &\lesssim
    \frac{\|w-w'\|_2}{\sqrt{n} h^{m/2+1}}
  \end{align*}
  so we take $L = \frac{1}{\sqrt{n} h^{m/2+1}}$.
  Noting
  $\psi^{-1}(t) = \sqrt{\log(1+t)} + \frac{\log(1+t)}{2\sqrt{\log n}}$
  and $N_\delta \lesssim \delta^{-m}$,
  \begin{align*}
    J_\psi(\delta)
    &=
    \int_0^\delta
    \psi^{-1}\big( N_\varepsilon \big)
    \diff{\varepsilon}
    + \delta
    \psi^{-1} \big( N_\delta \big)
    \lesssim
    \frac{\delta \log(1/\delta)}{\sqrt{\log n}}
    + \delta \sqrt{\log(1/\delta)}
    \lesssim
    \delta \sqrt{\log n}, \\
    J_2(\delta)
    &=
    \int_0^\delta
    \sqrt{\log N_\varepsilon}
    \diff{\varepsilon}
    \lesssim
    \delta \sqrt{\log(1/\delta)}
    \lesssim
    \delta \sqrt{\log n}.
  \end{align*}

  \proofparagraph{strong approximation}

  Recalling that
  $\tilde\varepsilon_i = \varepsilon_i$
  for all $i$ with high probability,
  by Proposition~\myref{pro:emp_proc},
  for all $t, \eta > 0$ there exists a
  zero-mean Gaussian process $T(w)$ satisfying
  \begin{align*}
    \E\left[
      \left(\sum_{i=1}^n f_w(W_i, \varepsilon_i)\right)
      \left(\sum_{i=1}^n f_{w'}(W_i, \varepsilon_i)\right)
    \right]
    &= \E\big[ T(w) T(w')
    \big]
  \end{align*}
  for all $w, w' \in \cW$ and
  \begin{align*}
    &\P\left(
      \sup_{w \in \cW}
      \left| \sum_{i=1}^n f_{w}(W_i, \varepsilon_i)
      - T(w) \right|
      \geq C_\psi(t + \eta)
    \right) \\
    &\quad\leq
    C_\psi
    \inf_{\delta > 0}
    \inf_{\cF_\delta}
    \Bigg\{
      \frac{\beta_\delta^{1/3} (\log 2 |\cF_\delta|)^{1/3}}{\eta }
      + \left(\frac{\sqrt{\log 2 |\cF_\delta|}
      \sqrt{\E\left[\|\Omega_\delta\|_2\right]}}{\eta }\right)^{2/3} \\
      &\qquad+
      \psi\left(\frac{t}{L J_\psi(\delta)}\right)^{-1}
      + \exp\left(\frac{-t^2}{L^2 J_2(\delta)^2}\right)
    \Bigg\} \\
    &\quad\leq
    C_\psi
    \Bigg\{
      \frac{
        \left(\frac{\log n} {n^2 h^{2m} \delta^{m}} \right)^{1/3}
      (\log n)^{1/3}}{\eta }
      + \left(\frac{\sqrt{\log n}
          \sqrt{n^{-3/2} h^{-7m/6} \delta^{-m}}
      }{\eta }\right)^{2/3} \\
      &\qquad+
      \psi\left(\frac{t}{\frac{1}{\sqrt{n} h^{m/2+1}}
      J_\psi(\delta)}\right)^{-1}
      + \exp\left(\frac{-t^2}{
          \left( \frac{1}{\sqrt{n} h^{m/2+1}} \right)^2
      J_2(\delta)^2}\right)
    \Bigg\} \\
    &\quad\leq
    C_\psi
    \Bigg\{
      \frac{
      (\log n)^{2/3}}{n^{2/3} h^{2m/3} \delta^{m/3} \eta}
      + \left(\frac{
        n^{-3/4} h^{-7m/12} \delta^{-m/2} \sqrt{\log n}}
      {\eta }\right)^{2/3} \\
      &\qquad+
      \psi\left(\frac{t\sqrt{n} h^{m/2+1}}
      {\delta \sqrt{\log n}}\right)^{-1}
      + \exp\left(\frac{-t^2n h^{m+2}}
      {\delta^2 \log n}\right)
    \Bigg\}.
  \end{align*}
  Noting $\psi(x) \geq e^{x^2/4}$ for $x \leq 4 \sqrt{\log n}$,
  any $R_n \to \infty$ gives the probability bound
  \begin{align*}
    \sup_{w \in \cW}
    \left| \sum_{i=1}^n f_{w}(W_i, \varepsilon_i)
    - T(w) \right|
    &\lesssim_\P
    \frac{(\log n)^{2/3}}{n^{2/3} h^{2m/3} \delta^{m/3}} R_n
    + \frac{\sqrt{\log n}}{n^{3/4} h^{7m/12} \delta^{m/2}} R_n
    + \frac{\delta \sqrt{\log n}} {\sqrt{n} h^{m/2+1}}.
  \end{align*}
  Optimizing over $\delta$ gives
  $\delta \asymp \left(\frac{\log n}{n h^{m-6}}\right)^{\frac{1}{2m+6}}
  = h \left( \frac{\log n}{n h^{3m}} \right)^{\frac{1}{2m+6}}$
  and so
  \begin{align*}
    \sup_{w \in \cW}
    \left| \sum_{i=1}^n f_{w}(W_i, \varepsilon_i)
    - T(w) \right|
    &\lesssim_\P
    \left(
      \frac{(\log n)^{m+4}}{n^{m+4}h^{m(m+6)}}
    \right)^{\frac{1}{2m+6}} R_n.
  \end{align*}

  \proofparagraph{convergence of $\hat H(w)$}

  For $1 \leq j,l \leq k$
  define the zero-mean random variables
  \begin{align*}
    u_{i j l}(w)
    &=
    K_h(W_i-w) p_h(W_i-w)_j p_h(W_i-w)_l \\
    &\quad-
    \E\big[K_h(W_i-w) p_h(W_i-w)_j p_h(W_i-w)_l \big]
  \end{align*}
  and note that
  $|u_{i j l}(w)| \lesssim h^{-m}$.
  By Lemma~\ref{lem:exponential_mixing}\ref{eq:exponential_mixing_bounded}
  for a constant $C_2 > 0$ and all $t > 0$,
  \begin{align*}
    \P\left(
      \left|
      \sum_{i=1}^n
      u_{i j l}(w)
      \right|
      > C_2 h^{-m} \big( \sqrt{n t}
      + (\log n)(\log \log n) t \big)
    \right)
    &\leq
    C_2 e^{-t}.
  \end{align*}
  Further, note that by Lipschitz properties,
  \begin{align*}
    \left|
    \sum_{i=1}^n u_{i j l}(w)
    - \sum_{i=1}^n  u_{i j l}(w')
    \right|
    &\lesssim
    h^{-m-1} \|w-w'\|_2
  \end{align*}
  so there is a $\delta$-cover of $(\cW, \|\cdot\|_2)$
  with size at most $n^a \delta^{-a}$ for some $a > 0$.
  Adjusting $C_2$,
  \begin{align*}
    \P\left(
      \sup_{w \in \cW}
      \left|
      \sum_{i=1}^n
      u_{i j l}(w)
      \right|
      > C_2 h^{-m} \big( \sqrt{n t}
      + (\log n)(\log \log n) t \big)
      + C_2 h^{-m-1} \delta
    \right)
    &\leq
    C_2 n^a \delta^{-a}
    e^{-t}
  \end{align*}
  and hence
  \begin{align*}
    \sup_{w \in \cW}
    \left|
    \sum_{i=1}^n
    u_{i j l}(w)
    \right|
    &\lesssim_\P
    h^{-m} \sqrt{n \log n}
    + h^{-m} (\log n)^3
    \lesssim_\P
    \sqrt{\frac{n \log n}{h^{2m}}}.
  \end{align*}
  Therefore
  \begin{align*}
    \sup_{w\in\cW} \|\hat H(w)-H(w)\|_2
    &\lesssim_\P
    \sqrt{\frac{n \log n}{h^{2m}}}.
  \end{align*}

  \proofparagraph{bounding the matrix term}

  Firstly note that,
  since $\sqrt{\frac{\log n}{n h^{2m}}} \to 0$,
  we have that uniformly in $w \in \cW$
  \begin{align*}
    \|\hat H(w)^{-1}\|_2
    \leq
    \frac{\|H(w)^{-1}\|_2}
    {1 - \|\hat H(w)-H(w)\|_2 \|H(w)^{-1}\|_2}
    &\lesssim_\P
    \frac{1/n}
    {1 - \sqrt{\frac{n \log n}{h^{2m}}} \frac{1}{n}}
    \lesssim_\P
    \frac{1}{n}.
  \end{align*}
  Therefore
  \begin{align*}
    &\sup_{w \in \cW}
    \big|
    e_1^\T \big(\hat H(w)^{-1} - H(w)^{-1}\big)
    S(w)
    \big|
    \leq
    \sup_{w \in \cW}
    \big\|\hat H(w)^{-1} - H(w)^{-1}\big\|_2
    \|S(w)\|_2 \\
    &\quad\leq
    \sup_{w \in \cW}
    \big\|\hat H(w)^{-1}\big\|_2
    \big\|H(w)^{-1}\big\|_2
    \big\|\hat H(w) - H(w)\big\|_2
    \|S(w)\|_2
    \lesssim_\P
    \sqrt{\frac{\log n}{n^3 h^{2m}}}
    \sup_{w \in \cW}
    \|S(w)\|_2.
  \end{align*}
  Now for $1 \leq j \leq k$ write
  $u_{i j}(w) = K_h(W_i-w) p_h(W_i-w)_j \tilde \varepsilon_i$
  so that $S(w)_j = \sum_{i=1}^n u_{i j}(w)$ with high probability.
  Note that $u_{i j}(w)$ are zero-mean with
  $\Cov[u_{i j}(w), u_{i' j}(w)] = 0$ for $ i \neq i'$.
  Also $|u_{i j}(w)| \lesssim h^{-m} \log n$
  and $\Var[u_{i j}(w)] \lesssim h^{-m}$.
  Thus by
  Lemma~\ref{lem:exponential_mixing}\ref{eq:exponential_mixing_bernstein}
  for a constant $C_3>0$,
  \begin{align*}
    \P\left(
      \Big| \sum_{i=1}^n u_{i j}(w) \Big|
      \geq C_3 \big( (h^{-m/2} \sqrt n + h^{-m} \log n) \sqrt t
      + h^{-m} (\log n)^3 t \big)
    \right)
    &\leq
    C_3 e^{-t}, \\
    \P\left(
      \Big| \sum_{i=1}^n u_{i j}(w) \Big|
      >
      C_3 \left(
        \sqrt{\frac{tn}{h^{m}}}
        + \frac{t(\log n)^3}{h^{m}}
      \right)
    \right)
    &\leq
    C_3 e^{-t},
  \end{align*}
  where we used $n h^{m} \gtrsim (\log n)^2$
  and adjusted the constant if necessary.
  As before,
  $u_{i j}(w)$ is Lipschitz in $w$ with a constant which is at most
  polynomial in $n$,
  so for some $a>0$
  \begin{align*}
    \P\left(
      \sup_{w \in \cW}
      \Big| \sum_{i=1}^n u_{i j}(w) \Big|
      >
      C_3 \left(
        \sqrt{\frac{tn}{h^{m}}}
        + \frac{t(\log n)^3}{h^{m}}
      \right)
    \right)
    &\leq
    C_3 n^a e^{-t}, \\
    \sup_{w \in \cW}
    \|S(w)\|_2
    \lesssim_\P
    \sqrt{\frac{n \log n}{h^{m}}}
    + \frac{(\log n)^4}{h^{m}}
    &\lesssim_\P
    \sqrt{\frac{n \log n}{h^{m}}}
  \end{align*}
  as $n h^m \gtrsim (\log n)^7$.
  Finally
  \begin{align*}
    \sup_{w \in \cW}
    \big|
    e_1^\T \big(\hat H(w)^{-1} - H(w)^{-1}\big)
    S(w)
    \big|
    &\lesssim_\P
    \sqrt{\frac{\log n}{n^3 h^{2m}}}
    \sqrt{\frac{n \log n}{h^{m}}}
    \lesssim_\P
    \frac{\log n}{\sqrt{n^2 h^{3m}}}.
  \end{align*}

  \proofparagraph{bounding the bias}

  Since $\mu \in \cC^\gamma$, we have, by the multivariate version of Taylor's
  theorem,
  \begin{align*}
    \mu(W_i)
    &=
    \sum_{|\kappa|=0}^{\gamma-1}
    \frac{1}{\kappa!}
    \partial^{\kappa} \mu(w)
    (W_i-w)^\kappa
    + \sum_{|\kappa|=\gamma}
    \frac{1}{\kappa!}
    \partial^{\kappa} \mu(w')
    (W_i-w)^\kappa
  \end{align*}
  for some $w'$ on the line segment connecting
  $w$ and $W_i$.
  Now since $p_h(W_i-w)_1 = 1$,
  \begin{align*}
    &e_1^\T \hat H(w)^{-1}
    \sum_{i=1}^n K_h(W_i-w) p_h(W_i-w) \mu(w) \\
    &\quad=
    e_1^\T \hat H(w)^{-1}
    \sum_{i=1}^n K_h(W_i-w) p_h(W_i-w) p_h(W_i-w)^\T e_1 \mu(w)
    = e_1^\T e_1 \mu(w) = \mu(w).
  \end{align*}
  Therefore
  \begin{align*}
    \Bias(w)
    &=
    e_1^\T \hat H(w)^{-1}
    \sum_{i=1}^n K_h(W_i-w) p_h(W_i-w) \mu(W_i)
    - \mu(w) \\
    &=
    e_1^\T \hat H(w)^{-1}
    \sum_{i=1}^n K_h(W_i-w) p_h(W_i-w) \\
    &\quad\times
    \Bigg(
      \sum_{|\kappa|=0}^{\gamma-1}
      \frac{1}{\kappa!}
      \partial^{\kappa} \mu(w)
      (W_i-w)^\kappa
      + \sum_{|\kappa|=\gamma}
      \frac{1}{\kappa!}
      \partial^{\kappa} \mu(w')
      (W_i-w)^\kappa
      - \mu(w)
    \Bigg) \\
    &=
    \sum_{|\kappa|=1}^{\gamma-1}
    \frac{1}{\kappa!}
    \partial^{\kappa} \mu(w)
    e_1^\T \hat H(w)^{-1}
    \sum_{i=1}^n K_h(W_i-w) p_h(W_i-w)
    (W_i-w)^\kappa \\
    &\quad+
    \sum_{|\kappa|=\gamma}
    \frac{1}{\kappa!}
    \partial^{\kappa} \mu(w')
    e_1^\T \hat H(w)^{-1}
    \sum_{i=1}^n K_h(W_i-w) p_h(W_i-w)
    (W_i-w)^\kappa \\
    &=
    \sum_{|\kappa|=\gamma}
    \frac{1}{\kappa!}
    \partial^{\kappa} \mu(w')
    e_1^\T \hat H(w)^{-1}
    \sum_{i=1}^n K_h(W_i-w) p_h(W_i-w)
    (W_i-w)^\kappa,
  \end{align*}
  where we used that
  $p_h(W_i-w)$ is a vector containing monomials
  in $W_i-w$ of order up to $\gamma$, so
  $e_1^\T \hat H(w)^{-1}
  \sum_{i=1}^n K_h(W_i-w) p_h(W_i-w)
  (W_i-w)^\kappa = 0$
  whenever $1 \leq |\kappa| \leq \gamma$.
  Finally
  \begin{align*}
    &\sup_{w\in\cW}
    |\Bias(w)| \\
    &\quad=
    \sup_{w\in\cW}
    \left|
    \sum_{|\kappa|=\gamma}
    \frac{1}{\kappa!}
    \partial^{\kappa} \mu(w')
    e_1^\T \hat H(w)^{-1}
    \sum_{i=1}^n K_h(W_i-w) p_h(W_i-w)
    (W_i-w)^\kappa
    \right| \\
    &\quad\lesssim_\P
    \sup_{w\in\cW}
    \max_{|\kappa| = \gamma}
    \left|
    \partial^{\kappa} \mu(w')
    \right|
    \|\hat H(w)^{-1}\|_2
    \left\|
    \sum_{i=1}^n K_h(W_i-w) p_h(W_i-w)
    \right\|_2
    h^\gamma \\
    &\quad\lesssim_\P
    \frac{h^\gamma}{n}
    \sup_{w\in\cW}
    \left\|
    \sum_{i=1}^n K_h(W_i-w) p_h(W_i-w)
    \right\|_2.
  \end{align*}
  Now write
  $\tilde u_{i j}(w) = K_h(W_i-w)p_h(W_i-w)_j$
  and note that
  $|\tilde u_{i j}(w)| \lesssim h^{-m}$
  and
  $\E[\tilde u_{i j}(w)] \lesssim 1$.
  By Lemma~\ref{lem:exponential_mixing}%
  \ref{eq:exponential_mixing_bounded},
  for a constant $C_4$,
  \begin{align*}
    \P\left(
      \left|
      \sum_{i=1}^n \tilde u_{i j}(w)
      - \E\left[
        \sum_{i=1}^n \tilde u_{i j}(w)
      \right]
      \right|
      > C_4 h^{-m} \big( \sqrt{n t}
      + (\log n)(\log \log n) t \big)
    \right)
    &\leq
    C_4 e^{-t}.
  \end{align*}
  As in previous parts, by Lipschitz properties,
  this implies
  \begin{align*}
    \sup_{w \in \cW}
    \left|
    \sum_{i=1}^n \tilde u_{i j}(w)
    \right|
    &\lesssim_\P
    n
    \left(
      1 + \sqrt{\frac{\log n}{n h^{2m}}}
    \right)
    \lesssim_\P
    n.
  \end{align*}
  Therefore
  $\sup_{w\in\cW} |\Bias(w)|
  \lesssim_\P n h^\gamma / n
  \lesssim_\P h^\gamma$.

  \proofparagraph{conclusion}

  By the previous parts,
  \begin{align*}
    \sup_{w \in \cW}
    \left|\hat \mu(w) - \mu(w) - T(w) \right|
    &\leq
    \sup_{w \in \cW}
    \left|e_1^\T H(w)^{-1} S(w) - T(w) \right| \\
    &\quad+
    \sup_{w \in \cW}
    \left| e_1^\T \big(\hat H(w)^{-1} - H(w)^{-1}\big) S(w) \right|
    + \sup_{w \in \cW}
    |\Bias(w)| \\
    &\lesssim_\P
    \left(
      \frac{(\log n)^{m+4}}{n^{m+4}h^{m(m+6)}}
    \right)^{\frac{1}{2m+6}} R_n
    + \frac{\log n}{\sqrt{n^2 h^{3m}}}
    + h^\gamma \\
    &\lesssim_\P
    \frac{R_n}{\sqrt{n h^m}}
    \left(
      \frac{(\log n)^{m+4}}{n h^{3m}}
    \right)^{\frac{1}{2m+6}}
    + h^\gamma,
  \end{align*}
  where the last inequality follows because
  $n h^{3m} \to \infty$
  and $\frac{1}{2m+6} \leq \frac{1}{2}$.
  Finally, we verify the upper and lower bounds
  on the variance of the Gaussian process.
  Since the spectrum of $H(w)^{-1}$
  is bounded above and below by $1/n$,
  \begin{align*}
    \Var[T(w)]
    &=
    \Var\left[
      e_1^\T H(w)^{-1}
      \sum_{i=1}^n K_h(W_i-w) p_h(W_i-w) \varepsilon_i
    \right] \\
    &=
    e_1^\T H(w)^{-1}
    \Var\left[
      \sum_{i=1}^n K_h(W_i-w) p_h(W_i-w) \varepsilon_i
    \right]
    H(w)^{-1} e_1^\T \\
    &\lesssim
    \|H(w)^{-1}\|_2^2
    \max_{1 \leq j \leq k}
    \sum_{i=1}^n
    \Var\big[
      K_h(W_i-w) p_h(W_i-w)_j \sigma(W_i)
    \big] \\
    &\lesssim
    \frac{1}{n^2} n
    \frac{1}{h^m}
    \lesssim
    \frac{1}{n h^m}.
  \end{align*}
  $\Var[T(w)] \gtrsim \frac{1}{n h^m}$
  by the same argument given to bound the eigenvalues of
  $H(w)^{-1}$.
\end{myproof}
 
\end{document}